\documentclass[11pt]{article}
\usepackage{srcltx}
\usepackage{epsfig}
\usepackage{url}
\usepackage[colorlinks=true]{hyperref}
\def\myP{{\Pi}}
\def\mycO{{\cal SO}}
\def\cO{{\cal O}}

\def\bZ{{\mathbf{Z}}}

\def\cC{{\cal C}}

\newcommand{\wh}[1]{{\widehat#1}}
\newcommand{\ov}[1]{{\overline#1}}
\newcommand{\wt}[1]{{\widetilde#1}}

\def\r{{\hbox{\tiny\rm r}}}
\def\b{{\hbox{\tiny\rm b}}}

\def\RiskOpt{\mathrm{RiskOpt}}
\usepackage{psfrag}
\usepackage{graphicx}
\usepackage[utf8]{inputenc}
\usepackage{subfigure}
\usepackage{hyperref}
\usepackage{mathtools}
\usepackage{amsfonts,amsmath,amssymb,amsthm}
\usepackage{xcolor}
\usepackage{booktabs}

\usepackage{float,enumitem}
\setlist{itemsep=3pt,parsep=0pt,topsep=3pt}

\newcommand{\eps}{\varepsilon}

\floatstyle{ruled}
\newfloat{algorithm}{tbp}{loa}
\floatname{algorithm}{Algorithm}

\DeclareMathOperator{\rank}{rank}

\DeclareMathOperator*{\argmin}{Arg\,min}

\newtheorem{thm}{Theorem}

\newtheorem{lemma}[thm]{Lemma}

\newcommand{\cP}{{\cal P}}
\def\qed{\ \hfill$\square$\par\smallskip}

\def\Opt{{\mathop{\hbox{Opt}}}}

\def\rank{{\mathop{\hbox{\rm Rank}}}}

\newcommand{\cN}{I\!\! N}

\newcommand{\E}{{\cal E}}

\newcommand{\N}{{\cal N}}
\newcommand{\I}{{\cal I}}

  \newcommand{\X}{{\cal X}}

\newcommand{\half}{ \mbox{\small$\frac{1}{2}$}}

\def\mR{{\bar{\varrho}^*}}
\def\mr{{\mathfrak{r}}}
\def\myr{\mathfrak{r}}
\def\myrg{\mathfrak{g}}

\def\eps{\varepsilon}

\def\e{\epsilon}
\def\ex{\mathrm{e}}

\def\inter{{\hbox{\rm int}\,}}

\def\cM{{\cal M}}
\def\Prob{{\hbox{\rm Prob}}}

\newcommand{\be}{\begin{eqnarray}}
\newcommand{\ee}[1]{\label{#1}\end{eqnarray}}
\newcommand{\nn}{\nonumber \\}
\newcommand{\ese}{\end{eqnarray*}}
\newcommand{\bse}{\begin{eqnarray*}}
\newcommand{\rf}[1]{~(\ref{#1})}

\newtheorem{proposition}{Proposition}
\newtheorem{corollary}{Corollary}
\newtheorem{theorem}{Theorem}
\def\argmin{\mathop{\hbox{\rm argmin$\,$}}}
\def\Argmin{\mathop{\hbox{\rm Argmin$\,$}}}

\def\Risk{\hbox{\rm Risk}}

\def\cU{{\cal U}}
\def\cA{{\cal A}}
\def\cS{{\cal S}}

\def\cH{{\cal H}}
\def\cI{{\cal I}}
\def\cR{{\cal R}}

\def\cN{{\cal N}}

\def\bR{{\mathbf{R}}}

\def\cB{{\cal B}}

\def\cF{{\cal F}}
\def\cG{{\cal G}}
\def\cX{{\cal X}}
\def\T{{\cal T}}
\def\cJ{{\cal J}}

\definecolor{darkmagenta}{RGB}{125,38,205}

\def\myC{{\mathfrak{C}}}
\begin{document}

\title{Aggregating estimates by convex optimization}
\author{
Anatoli Juditsky
\thanks{LJK, Universit\'e Grenoble Alpes, 700 Avenue Centrale,  38401 Domaine Universitaire de Saint-Martin-d'Hères, France,	
{\tt anatoli.juditsky@univ-grenoble-alpes.fr}}
\and Arkadi Nemirovski
\thanks{ISyE, Georgia Institute
 of Technology, Atlanta, Georgia
30332, USA, {\tt arkadi.nemirovski@isye.gatech.edu}.}
\footnote{Research of the authors was supported by Multidisciplinary Institute in Artificial intelligence MIAI {@} Grenoble Alpes (ANR-19-P3IA-0003).
}}
\date{}
\maketitle
\begin{abstract}
We discuss the approach to estimate aggregation and adaptive estimation based upon (nearly optimal) testing of convex hypotheses.
We show that in the situation where the observations stem from {\em simple observation schemes} \cite{juditsky2020statistical} and where set of unknown signals is a finite union of convex and compact sets, the proposed approach leads to aggregation and adaptation routines with nearly optimal performance. As an illustration, we consider application of the proposed estimates to the problem of recovery of unknown signal known to belong to a union of ellitopes \cite{l2estimation,juditsky2020statistical} in Gaussian observation scheme. The proposed approach can be implemented efficiently when the number of sets in the union is  ``not very large.'' We conclude the paper with a small simulation study illustrating practical performance of the proposed procedures in the problem of signal estimation in the single-index model.
\end{abstract}

\section{Introduction}
We address the problem of data-driven selection of estimators from a given collection. A simplified version of the problem considered in this paper is as follows.
\begin{quotation}
\noindent{\bf Problem I.} We are given in advance $N$ nonempty convex compact signal sets $X_j\subset\bR^n$ and $m\times n$ sensing matrices $A_j$, $1\leq j\leq N$. Given access to $M$ independent observations
\be
\omega^M=(\omega_1,...,\omega_M): \omega_k=Ax + \sigma\xi_k,\,1\leq k\leq M,\quad \xi_k\sim\cN(0,I_m)
\ee{eq0}
we want to recover the signal $x\in \bR^n$ in the situation when it is known {\em a priori} that $x\in X_j$ and $A=A_j$ for some (unknown!) $j\leq N$. Given reliability tolerance $\epsilon\in(0,1)$,  we quantify the performance of a
candidate estimate $\widehat{x}(\omega^M):\bR^{mM}\to\bR^n$ by its worst-case $\epsilon$-risk---the radius of the smallest ball, in a  given seminorm\footnote{Recall that a seminorm on $\bR^n$
satisfies exactly the same requirements as a norm, with positivity outside the origin replaced with nonnegativity. A standard example of a seminorm is $\|x\|=\pi(Bx)$ where $\pi(\cdot)$ is a norm on some $\bR^m$
and  $B\in \bR^{m\times n}$ has a nontrivial kernel.} $\|\cdot\|$, around $\wh x$
which contains the signal $x$ underlying observations with probability at least $1-\epsilon$, that is, by the quantity
$$
{\Risk^{\overline{1,N}}_{\epsilon,M}}[\widehat{x}|X]:=\min\left\{\rho: \Prob_{{\omega^M\sim p_j^M}}\{\|\widehat{x}(\omega^K)-x\|>\rho\}\leq\epsilon\,\,\forall (j\leq N,x\in X_j)\right\}
$$
where $p_j$  is the normal distribution $\cN(A_jx, \sigma^2I_m)$ and $p_j^M=\overbrace{p_j\times...\times p_j}^{M}$.
We intend to estimate the signal by aggregating $N$ selected in advance ``preliminary'' estimates $\wt x_j$, $j=1,...,N$, $j$th of them associated with $j$th observation model in which $x$ is known to belong to $X_j$  and $A=A_j$. Specifically, we split $M$ available observations into ``pilot sample'' $\omega_1,...,\omega_{\overline{K}}$ used to build points $x_j=\wt x_j(\omega^{\ov K})$, and use the remaining $K=M-\ov K$ observations to ``assemble'' $x_j$ into the resulting estimate $\widehat{x}$ of the signal.
\end{quotation}
A related problem is that of constructing an  estimate which is  {\em adaptive}---such that its risk is ``as close as possible'' to the maximal risk of the $j$th estimate under the $j$th observation model, $1\leq j\leq N$.

In this work, our focus is on the aggregation step, thus, for the most of the exposition below, estimates $x_j=\wt x_j$, $j=1,...,N$, are regarded as known fixed points in $\bR^n$. The above problem is closely related to another fundamental statistical problem, that of aggregation and, in particular, to ``model selection'' version of the
problem in which
 the objective is to select the ``nearly $\|\cdot\|$-closest'' to $x$ point among given points $x_1,...,x_N$.
 Both problems have received a lot of attention in the statistical literature. The adaptive estimation problem, in its general form which is relevant for us, has been stated in O. Lepski's pioneering works \cite{Lepski1991,Lepski1991a,Lepski1992a,Lepski1992}
 (for the setting in which $\{X_j\}$ is an injected family of sets), then substantially generalized in \cite{goldenshluger2008universal,goldenshluger2009structural,goldenshluger2011bandwidth,lepski2015adaptive},
 giving rise to the celebrated Lepski's and Goldenshluger-Lepski's adaptation schemes put to use in various contexts and by various authors.
 A remarkable progress has also been achieved when solving the aggregation problem, in particular, in the context of $L_2$-estimation in the white noise model where exact oracle inequalities were derived for collections of arbitrary estimators. Specifically, the notion of optimal rates of aggregation has been introduced in \cite{tsybakov2003optimal}, and aggregation procedures attaining the risk which approaches the risk of the best point among $x_i$ with the smallest possible, in the minimax sense, remainder term have been introduced (see also \cite{audibert2004aggregated,yang2004aggregating,bunea2007aggregation,rigollet2007linear,juditsky2008learning,dai2012deviation} and references therein).
Aggregation of estimators with respect to other
loss functions has also been studied extensively. The problem of aggregating estimates with the Kullback–Leibler divergence as a loss function has been studied in \cite{catoni2004statistical,yang2000mixing} in the problem of density estimation and in \cite{rigollet2012kullback} for  generalized  linear models. Aggregation w.r.t. $L_1$-risk in the context of density estimation has been studied in \cite{yatracos1985rates,devroye1996universally,devroye2012combinatorial,mahalanabis2007density}; that approach has been extended to the regression setup in \cite{hengartner2001estimation}. Finally, one of our principal motivations comes from \cite{goldenshluger2009universal} where a general aggregation scheme which applies to wide variety of the risk measures have been proposed.
In this paper, we aim at extending adaptive estimation and estimate aggregation framework in several directions. Specifically, we propose adaptive estimation and aggregation routines for problems where indirect observations are available under general convex constraints on unknown signal.\footnote{
We should mention here a special status of the problem of adaptive estimation of general linear functionals of unknown signal: in a separate line of research \cite{cailow2004,cailow2005} the minimax affine estimator was used as  ``working horse'' to build the near-optimal estimator  of a linear functional over a finite union $X$ of convex compact sets in  Gaussian observation scheme.  A different general construction for nearly minimax optimal estimation of linear functionals over union of convex sets in simple observation schemes has been developed in \cite{juditsky2020near}.}
\paragraph{The underlying idea} of the proposed routines is that of pairwise comparison of candidate estimates: to decide if estimate $\wt x_i$ is better/worse than $\wt x_j$, $i\neq j$, we replace the relation  ``risk of $\wt x_i$ is less  than risk of $\wt x_j$'' with a pair of convex hypotheses about $x$. To see how this reduction operates, consider the situation of Problem I with $N=2$, where we want to choose between just two estimates $\wt x_1$ and $\wt x_2$, associated with models indexed by $j\in \{1,2\}$, assuming that $\epsilon$-risk of $\wt x_j$ is bounded with $\myr_j$ under the $j$th model.
For the sake of definiteness, assume that $\myr_1\leq \myr_2$, and that the realization of  noise in the preliminary observation
belongs to the subset of the corresponding probability space of probability $1-\epsilon$ such that $\|\wt x_{j_*}-x\|\leq \myr_{j_*}$,
where $j_*$ is the index of the ``true'' observation model, that is, $x\in X_{j_*}$ and $A=A_{j_*}$. In this case,
 if $j_*=1$, we have $x\in X_1$ with $\|x-\wt x_1\|\leq \myr_1$ with probability $1-\epsilon$ and
 $A=A_1$, i.e., \[
 x\in B_1:=\{x\in X_1,\,\|x-\wt x_1\|\leq \myr_1\},
 \] so that $A_1x$ belongs to the convex and compact set $Y_1=A_1B_1$. When $j_*=2$, we have
 \[x\in B_2:=\{x\in X_2,\,\|x-\wt x_2\|\leq \myr_2\}
  \]
  and $A_2x\in Y_2=A_2B_2$. Now, assuming that we do have $x\in B_j$ when the actual model if the $j$th one, $j=1,2$, given observation $\omega^K$, consider the problem of testing ``convex hypotheses''
\[
H_1:\,A_1x\in Y_1\;\;\mbox{and}\;\;H_2:\,A_2x\in Y_2.
\]
As it is well known (see, e.g., \cite{chernoff1952,Burnashev1979,Burnashev1982}), when $Y_1$ and $Y_2$ do not intersect, the optimal test (that with the smallest maximal risk) deciding on $H_1$ against $H_2$ in Gaussian o.s. is the likelihood ratio test of simple hypotheses
\[
\ov H_1:\,A_1x=\ov y_1\;\;\mbox{and}\;\;\ov H_2:\,A_2x=\ov y_2
\]
 where
\[
   (\ov y_1,\ov y_2)\in \Argmin_{y_1\in Y_1,y_2\in Y_2}\|y_1-y_2\|_2.
\]
Thus, assuming that two hypotheses can be separated with maximal risk $\leq \epsilon$, when the first model is true, $x\in B_1$ and $H_1$ holds,  the test will accept it (and reject $H_2$) with probability $1-\epsilon$, implying that the $2\epsilon$-risk of the estimate $\wh x(\omega^K)=\wt x_1$ is bounded with $\myr_1$, and ``symmetric'' bound holds when the second model is true. On the other hand, in the case the hypotheses cannot be separated $(1-\epsilon)$-reliably, selecting $\wh x=\wt x_1$ results in the $\epsilon$-risk of $\wh x$ bounded with $\myr_1$ when the first model is true, and with  $\myr_1+2\myr_2+2\myr_{12}$ where
\[
\myr_{12}=\min\big\{\half \|x_1-x_2\|: x_1\in B_1, x_2\in B_2\big\}
\]
in the case of the true second model. A simple calculation shows (cf. e.g., Theorem \ref{thm:genthmG} in Section \ref{sec:ellall}) that in the latter case the quantity $\myr_{12}$ is upper-bounded by the maximal risk of estimation over $X=X_1\cup X_2$. Note that if ``separation'' $\myr_{12}$ is majorated by $\myr_2$, estimate $\wh x$ is adaptive in the sense of \cite{Lepski1991}---when $x\in X_1$ the $\epsilon$-risk of $\wh x$ is bounded with $\myr_1$, and when $x\in X_2$ its risk is bounded with $\myr_1+2\myr_2+2\myr_{12}$ which is the same as $\myr_2$, up to a moderate absolute factor. On the other hand, if $\myr_{12}\gg\myr_2$,  the corresponding bound is the best one can achieve under the circumstances.
More generally, reducing the problem of risk minimization to that of pairwise testing of convex hypotheses makes the problem amenable to the machinery of nearly-optimal testing of convex hypotheses developed in \cite{GJN2015}.

The proposed approach  shares its motivation with another construction of estimates based on testing  multiple hypotheses---the $T$-estimators developed in \cite{birge2006model,birge2007model,birge2013robust}. When applied to Problem I, the latter approach amounts to building a net of
points $\{x_\tau\},\,\tau\in \cal T,$ in $X$ and selecting the estimate by applying pairwise tests to small Euclidean balls around images of $x_\tau$, $x_{\tau'}$ in the observation space. Note that, typically, $T$-estimators  cannot be obtained in a computationally efficient fashion and are usually considered as a theoretical tool to explore the properties of statistical problems. Despite obvious similarities with $T$-estimates (e.g., great flexibility shared by the both approaches), adaptive and aggregation estimates we discuss in this paper are of a different nature. Our approach can be seen as an ``operational counterpart'' to that of \cite{birge2006model}
leading to adaptive estimates which are efficiently computable provided the data of the problem---sets $X_j$ and norm $\|\cdot\|$---are computationally tractable, and $N$ is moderate (hundreds, perhaps, thousands). As the price to pay for generality of the proposed constructions, our  estimates and their risks (provably near-optimal, as we shall see, under natural assumptions) are given by efficient computation rather than in a closed analytic form.\footnote{We believe that in our setting, allowing for arbitrary sensing matrices and general convex parameter sets, closed form results are just impossible.} This is hardly a problem in application where efficient computation  usually is not inferior to a formula.

\paragraph{What is ahead.}
In what follows we discuss two adaptive estimates: a ``generic'' selection procedure in the situation where $\|\cdot\|$ is an arbitrary seminorm, and a special aggregation routine for the problem setting in which $\|\cdot\|$ is a Euclidean seminorm.
Our principal contribution (cf., e.g., Theorem \ref{thm:genthm} and Corollary \ref{cor:if} of Section \ref{sec:genest} in the case of general seminorm), as applied to Problem I above, may be summarized as follows.

Let a real $\theta\geq 1$, an integer $\overline K\geq 1$ and $\epsilon\in (0,1/2)$ be fixed.
Assume that we are given preliminary $\ov K$-observation estimates $\wt x_j(\cdot)$ along with reals $\myr_j$ such that
\[
{\Risk^{\{j\}}_{\epsilon,\ov K}}[\wt{x}_j|X_j]\leq\myr_j{\leq \theta} {\RiskOpt^{\{j\}}_{\epsilon,\ov K}}[X_j],\;\;j=1,...,N,
\]
where
for a nonempty $\cJ\subset\overline{1,N}$ and a $ K$-observation estimate $\widehat{x}(\cdot)$,
$$
\Risk^{\cJ}_{\epsilon,K}[\widehat{x}|\cup_{j\in\cJ}X_j]=\min\left\{\rho:\Prob_{\omega^{K}\sim p_j^{K}}\left\{\|\widehat{x}(\omega^{ K})-x\|>\rho\right\}\leq\epsilon\,\,\forall (j\in\cJ,x\in X_j)\right\}
$$
is the risk of estimate ``on the union of models with indexes from $\cJ$,'' and
\[
{\RiskOpt^{\cJ}_{\epsilon,K}[\cup_{j\in\cJ}X_j]}=\inf_{\wh x}{\Risk^{\cJ}_{\epsilon,K}[\widehat{x}|\cup_{j\in\cJ}X_j]}
\]
is the corresponding minimax risk.

{Now, suppose that given} $M\geq O(1) {\ln (N/\epsilon)\over \ln(1/\epsilon) } \ov K$ independent observations \rf{eq0}, we utilize the first $\ov K$ observations to build ``preliminary'' estimates $x_j=\wt x_j(\omega^{\ov K})$, $j=1,...,N$, and then proceed with selection procedure of Section \ref{sec:genest} using $K=M-\ov K$ remaining observations to aggregate points $x_j$ into an adaptive estimate $\wh x^{(a)}(\omega^M)$.
Then $\wh x^{(a)}(\omega^M)$ satisfies
\[
{\Risk^{\overline{1,N}}_{2\epsilon,M}}[\wh{x}^{(a)}|X]\leq {O(1)\theta\RiskOpt^{\overline{1,N}}_{\epsilon,\ov K}[X]}
\]
In other words, modulo logarithmic in $N$ increase in observation count and reliability parameter $\epsilon$ of the risk replaced with $2\epsilon$, estimate  $\wh x^{(a)}(\omega^M)$ is minimax optimal on $X$ within factor $O(1)\theta$.
\par Furthermore, the ``overall'' minimax risk $\RiskOpt^{\overline{1,N}}$ is nearly upper-bounded by the maximum of pairwise minimax risks $\RiskOpt^{\{i,j\}}$. Specifically, with $M$ as above,
$$
\RiskOpt^{{\overline{1,N}}}_{2\epsilon,M}[\cup_{j\leq N}X_j]\leq O(1)\max\limits_{i\neq j}\RiskOpt^{\{i,j\}}_{\epsilon,\ov K}[X_i\cup X_j].
$$
Finally,
suppose that $N$ models in Problem I are ordered, so that bounds $\myr_i$ for partial risks of estimates $\wt x_j(\omega^{\ov K})$ satisfy
\[
\myr_1\leq \myr_2\leq ...\leq \myr_N,
\]
and that minimax risks of estimation over ``pairwise unions'' ${\RiskOpt^{\{i,j\}}_{\epsilon,\ov K}}[X_i\cup X_j]$, $1\leq i,j\leq N$, are dominated by the pairwise maxima of the corresponding partial risks, i.e.,
\[
\max_{1\leq j\leq i}\RiskOpt^{\{i,j\}}_{\epsilon,\ov K}[X_i\cup X_j]\leq O(1)\RiskOpt^{\{i\}}_{\epsilon,\ov K}[X_i],\;\;i=1,...,N.
\]
Then we also have
\[
\Risk^{\{i\}}_{2\epsilon,M}[\wh{x}^{(a)}|X_i]\leq {O(1)\theta} \RiskOpt^{\{i\}}_{\epsilon,\ov K}[X_i],\;\;\forall i=1,...,N,
 \]
 i.e., estimate $\wh{x}^{(a)}$ is (again, up to logarithmic in $N$ increase in observation count and reliability parameter $\epsilon$ of the risk replaced with $2\epsilon$) minimax adaptive, within  factor $O(1)\theta$, in the sense of \cite{Lepski1991,Lepski1991a} over considered family of observation models.

Our results are not restricted to the Gaussian observation scheme \rf{eq0} and deal with {\sl simple observation schemes}\footnote{Our results can be easily extended to the more general case of {\em simple families}---families of distributions specified in terms of upper bounds on their moment-generating functions, see \cite{juditsky2020statistical} for details. Restricting the framework to the case of simple observation schemes is aimed at  streamlining the presentation.} (o.s.'s),
as defined in \cite{GJN2015,juditsky2020statistical}. Aside of Gaussian o.s., important examples of simple o.s. are
\begin{itemize}
\item {\em Poisson o.s.}, where $\omega_k$ are independent across $k$ identically distributed vectors with independent across $i\leq m$ entries $[\omega_k]_i\sim\hbox{\rm Poisson}(a_i^Tx)$, and
\item {\em Discrete o.s.}, where $\omega_k$ are independent across $k$ realizations of discrete random variable taking values $1,...,m$ with probabilities affinely parameterized by $x$.
\end{itemize}
The presentation is organized in two parts. In the first part, we consider the problem of adaptive and minimax estimation over the sets which are unions of convex sets---a generalization of the setting of Problem I to the case of simple o.s.. We start with stating the general estimation problem and provide an ``operational summary'' of results on testing in simple observation schemes in Section \ref{sec:simpleos}.  Section \ref{sec:genest} deals with adaptation in the case of a general seminorm $\|\cdot\|$, and Section \ref{sec:l2est} with the special case of Euclidean seminorm.
The second part of the paper deals with the problem of model selection aggregation. Although closely related to the problem of adaptive estimation,
this problem calls for different notion of optimality with respect to which estimation routines discussed in Section \ref{sec:adest} may be heavily suboptimal. The second part begins with a description of two ``abstract'' aggregation routines utilizing pairwise tests in Section \ref{sect:general}, which we specify for aggregation in simple o.s. in Section \ref{sec:simplebs5}. We consider next the application of these routines to signal recovery in the situation described in Section \ref{sect:setting11}.
We conclude the paper (Section \ref{sec:ellall}) detailing how the proposed approach can be used to build nearly minimax estimates in the problem of signal recovery in Gaussian o.s. when the signal set $X$ is a union of ellitopes (cf. \cite{l2estimation}); these results are accompanied by a small simulation study illustrating numerical performance of the proposed estimates in that problem.

Proofs of the results are postponed till the appendix.

\section{Preliminaries: testing convex hypotheses in simple observation schemes}\label{sec:simpleos}
\subsection{Simple observation schemes: definitions}\label{sec:sos}
All  developments {to follow} make use of the notion of a simple observation scheme, see \cite{juditsky2020statistical}. To make the presentation self-contained we start with explaining this notion here.
\par
Formally, a {\sl simple observation scheme} (o.s.) is a collection $\mycO=\left((\Omega,\myP ),\{p_\mu(\cdot):\mu\in\cM\},\cF\right)$, where
\begin{itemize}
\item $(\Omega,\myP )$ is an {\sl observation space}: $\Omega$ is a Polish (complete metric separable) space, and $\myP $ is a $\sigma$-finite $\sigma$-additive Borel reference measure on $\Omega$, such that $\Omega$ is the support of $\myP $;
\item $\{p_\mu(\cdot):\mu\in\cM\}$ is a parametric family of probability densities, specifically, $\cM$ is a convex relatively open set in some $\bR^M$, and for $\mu\in \cM$, $p_\mu(\cdot)$ is a probability density, taken w.r.t. $\myP $, on $\Omega$. We assume that the function $p_\nu(\omega)$ is positive and continuous in $(\mu,\omega)\in\cM\times\Omega$;
\item $\cF$ is a finite-dimensional linear subspace in the space of continuous functions on $\Omega$. We assume that $\cF$ contains constants and all functions of the form $\ln(p_\mu(\cdot)/p_\nu(\cdot))$, $\mu,\nu\in \cM$, and that the function
    \begin{equation}\label{PhicalO}
    \Phi_{\mycO}(\phi;\mu)=\ln\left(\int_\Omega \ex^{\phi(\omega)} p_\mu(\omega)\myP (d\omega)\right)
    \end{equation}
    is  real-valued on $\cF\times\cM$ and is {\sl concave} in $\mu\in \cM$; note that this function is automatically convex in $\phi\in\cF$. From real-valuedness, convexity-concavity and the fact that both $\cF$ and $\cM$ are convex and relatively open, it follows that $\Phi$ is continuous on $\cF\times\cM$.
\end{itemize}
\subsubsection{Examples of simple observation schemes}
As shown in \cite{juditsky2020statistical} (and can be immediately verified), the following o.s.'s are simple:
\begin{enumerate}
\item {\sl Gaussian o.s.}, where $\myP $ is the Lebesgue measure on $\Omega=\bR^d$, $\cM=\bR^d$, $p_\mu(\omega)$ is the density of the Gaussian distribution $\cN(\mu,I_d)$ (mean $\mu$, unit covariance), and $\cF$ is the family of affine functions on $\bR^d$. Gaussian o.s. with $\mu$ linearly parameterized by signal $x$ underlying observations is the standard observation model in signal processing;
\item {\sl Poisson o.s.}, where $\myP $ is the counting measure on the nonnegative integer $d$-dimensional lattice $\Omega=\bZ^d_+$,  $\cM=\bR^d_{++}=\{\mu=[\mu_1;...;\mu_d]>0\}$, $p_{\mu}$ is the density, taken w.r.t. $\myP $, of random $d$-dimensional vector with independent $\hbox{\rm Poisson}(\mu_i)$ entries, $i=1,..,d$, and $\cF$ is the family of all affine functions on $\Omega$. Poisson o.s. with $\mu$ affinely parameterized by signal $x$ underlying observation is the standard observation model in {\sl Poisson imaging};
\item {\sl Discrete o.s.}, where $\myP$ is the counting measure on the finite set $\Omega=\{1,2,...,d\}$, $\cM$ is the set of positive $d$-dimensional probabilistic vectors $\mu=[\mu_1;...;\mu_d]$, $p_\mu(\omega)=\mu_\omega$, $\omega\in\Omega$, is the density, taken w.r.t. $\myP$, of a probability distribution $\mu$ on $\Omega$, and $\cF=\bR^d$ is the space of all real-valued functions on   $\Omega$;
\item {\sl Direct product of simple o.s.'s}. Given $K$ simple o.s.'s $\mycO_t=((\Omega_t,\myP_t),\{p_{t,\mu}:\mu\in\cM_r\},\cF_t)$, $t=1,...,K$, we can build from them a new (direct product) o.s. $\mycO_1\times....\times\mycO_K$  with  observation space $\Omega_1\times...\times\Omega_K$, reference measure $\myP_1\times...\times \myP_K$, family of probability densities $\{p_\mu(\omega_1,...,\omega_K)=\prod_{t=1}^Kp_{t,\mu_t}(\omega_t):\mu=[\mu_1;...;\mu_K]\in\cM_1\times...\times\cM_K\}$, and $\cF=\{\phi(\omega_1,...,\omega_K)=\sum_{t=1}^K\phi_t(\omega_t):\phi_t\in\cF_t,t\leq K\}$. In other words, the direct product of o.s.'s $\mycO_t$ is the observation scheme in which we observe collections $\omega^K=(\omega_1,...,\omega_K)$ with independent across $t$ components $\omega_t$ yielded by o.s.'s $\mycO_t$.
    \par
    When all factors $\mycO_t$, $t=1,...,K$, are identical to each other, we can reduce the direct product $\mycO_1\times...\times\mycO_K$ to its ``diagonal,''
    referred to as {\sl $K$th power $\mycO^K$}, or {\sl stationary $K$-repeated version}, of $\mycO=\mycO_1=...=\mycO_K$. Just as in the direct product case, the observation space and reference measure in $\mycO^K$ are $\Omega^K=\underbrace{\Omega\times...\times\Omega}_{K}$ and $\myP^K=\underbrace{\myP\times...\times\myP}_{K}$, the family of densities is $\{p^K_\mu(\omega^K)=\prod_{t=1}^Kp_\mu(\omega_t):\mu\in\cM\}$, and the family $\cF$ is $\{\phi^{(K)}(\omega_1,...,\omega_K)=\sum_{t=1}^K\phi(\omega_t):\phi\in\cF\}$. Informally, $\mycO^K$ is the observation scheme we arrive at when passing from a single observation drawn from a distribution $p_\mu$, $\mu\in\cM$,  to $K$ independent observations drawn from the same distribution $p_\mu$.\par
    It is immediately seen that direct product of simple o.s.'s, same as power of simple o.s., are themselves simple o.s.
\end{enumerate}
\subsection{Testing pairs of convex hypotheses in simple o.s.}\label{sectpairs}
What follows is a summary of  results of \cite{juditsky2020statistical} which are relevant to our current needs.\\
Assume that $\omega^K=(\omega_1,...,\omega_K)$ is a stationary $K$-repeated observation in a simple o.s. $\mycO=((\Omega,\myP),\{p_\mu:\mu\in\cM\},\cF)$, so that $\omega_1,...,\omega_K$ are, independently of each other, drawn from a distribution $p_\mu$ with some $\mu\in\cM$. Given $\omega^K$ we want to decide on the hypotheses $H_1$ and $H_2$,
 with $H_\chi$, $\chi=1,2$,  stating that $\omega_t\sim p_\mu$ for some $\mu\in M_\chi$,  where $M_\chi$ is a nonempty convex compact subset of $\cM$.
 In the sequel, we refer to hypotheses of this type, parameterized  by nonempty convex compact subsets of $\cM$,
 as to {\sl convex} hypotheses in the simple o.s. in question.\par
The principal ``building block'' of our subsequent constructions is a {\em simple} test\footnote{A test deciding on a pair of hypotheses is called simple, if given an observation, it always accepts exactly one of the hypotheses and rejects the other one.} $\T^K$ for this problem which is as follows:
\begin{itemize}
\item Given convex compact sets $M_\chi$, $\chi=1,2$, we solve the optimization problem
\begin{equation}\label{eq1}
\Opt=\max_{\mu\in M_1,\nu\in M_2}\ln\Big(\underbrace{\int_\Omega\sqrt{p_\mu(\omega)p_\nu(\omega)}\myP(d\omega)}_{=:\varrho(\mu,\nu)}\Big)
\end{equation}
It is shown in \cite{GJN2015} that in the case of simple o.s., problem (\ref{eq1}) is a convex problem (convexity meaning that the objective to be maximized is a concave continuous function of $\mu,\nu$) and an optimal solution exists.
\begin{quote}
\noindent Note that for basic simple o.s.'s problem (\ref{eq1}) reads
\begin{equation}\label{eq11}
\Opt=\max_{\mu\in M_1,\nu\in M_2}\left\{\begin{array}{ll}
-\hbox{\small $1\over 8$}\|\mu-\nu\|_2^2,&\hbox{Gaussian o.s.}\\
-{\half}\sum_{i=1}^d [\sqrt{\mu_i}-\sqrt{\nu_i}]^2,&\hbox{Poisson o.s.}\\
\ln\left(\sum_{i=1}^d\sqrt{\mu_i\nu_i}\right),&\hbox{Discrete o.s.}\\
\end{array}\right.
\end{equation}
\end{quote}
\item An optimal solution $\mu_*$, $\nu_*$ to (\ref{eq1}) induces  {\sl detectors}
\begin{equation}\label{eq2}
\begin{array}{rcl}\phi_*(\omega)&=&{\half}\ln(p_{\mu_*}(\omega)/p_{\nu_*}(\omega)):\;\Omega\to\bR,
\\
\phi_*^{(K)}(\omega^K)&=&\sum_{t=1}^K\phi_*(\omega_t):\;\; \Omega\times...\times\Omega\to\bR
\end{array}\end{equation}
Given a stationary $K$-repeated observation $\omega^K$, the test $\T^K$ accepts hypothesis $H_1$ and rejects hypothesis $H_2$ whenever $\phi_*^{(K)}(\omega^K) \geq0$, otherwise the test rejects $H_1$ and accepts $H_2$. The {\sl risk} of $\T^K$ -- the maximal probability to reject a hypothesis when it is true -- does not exceed $\epsilon_\star^K,$ where
$$\epsilon_\star=\exp(\Opt).
$$
In other words, whenever observation $\omega^K$ stems from a distribution $p_\mu$ with $\mu\in M_1\cup M_2$,
\begin{itemize}
\item the $p_\mu$-probability to reject $H_1$ when the hypothesis is true (i.e., when $\mu\in M_1$) is at most $\epsilon_\star^K$,     and
\item the $p_\mu$-probability to reject $H_2$ when the hypothesis is true (i.e., when $\mu\in M_2$) is at most $\epsilon_\star^K$.
\end{itemize}
\end{itemize}
The test $\T^K$ possesses the following optimality properties:
\begin{itemize}
\item[{\bf A.}] The associated detector {\sl $\phi_*^{(K)} $ and the risk $\epsilon_\star^K$ form an optimal solution and the optimal value in the optimization problem
\[
\begin{array}{c}
\min\limits_\phi\max\left[\max_{\mu\in M_1}\int_{\Omega^K} \e^{-\phi(\omega^K)}p^{K}_\mu(\omega^K)\myP^K(d\omega^K),
\max_{\nu\in M_2}\int_{\Omega^K} \e^{\phi(\omega^K)}p^{K}_\nu(\omega^K)\myP^K(d\omega^K)\right],\\
\big[\Omega^K=\underbrace{\Omega\times...\times\Omega}_{K},\;\;p_\mu^{K}(\omega^K)=\prod_{t=1}^Kp_\mu(\omega_t),\,
\big]\\
\end{array}
\]
where the minimum is taken w.r.t. all Borel functions $\phi(\cdot):\Omega^K\to\bR$;}
\item[{\bf B.}] {\sl Let $\epsilon\in (0,1/2)$, and suppose that there exists a test which, using a stationary $\overline{K}$-repeated observation, decides on the hypotheses $H_1$, $H_2$ with risk $\leq\epsilon$. Then
    \begin{equation}\label{epsilonstar}
    \epsilon_\star\leq [2\sqrt{\epsilon(1-\epsilon)}]^{1/\overline{K}}
    \end{equation}
    and the test $\T^K$ with\footnote{Here $\lceil a\rceil$ stands for the ``upper'' integer part---the smallest integer greater or equal to $a$.}
    $$
    K=\left\lceil {2\ln(1/\epsilon)\over \ln\left( [4\epsilon(1-\epsilon)]^{-1}\right)} \overline{K}\right\rceil
    $$
    decides on the hypotheses $H_1,H_2$ with risk $\leq\epsilon$ as well. Note that $K=2(1+o(1))\overline{K}$ as $\epsilon\to+0$.}\footnote{It is worth mentioning that in the Gaussian o.s. test $\T^K$ optimal---it is the test minimizing the maximal risk of testing of $H_1$ vs $H_2$ among all tests; the corresponding optimal risk is $\epsilon =1-\Phi\big(\half \|\mu_*-\nu_*\|_2\sqrt{K}\big)
    $ where $\Phi$ is the standard normal c.d.f. and $[\mu_*;\nu_*]$ is an optimal solution to \rf{eq1}.}
\end{itemize}
In what follows we augment the test $\T^K$ to address the situation where one {or both hypotheses are empty. When one of the hypotheses is empty, $\T^K$, by convention, accepts the nonempty} hypothesis. When both hypotheses are empty, $\T^K$ accepts, say, the first of them. Because the true hypothesis cannot be empty, the risk of $\T^K$ vanishes in this case.

\subsection{Testing multiple hypotheses in simple o.s.}\label{sectinfcol}
As shown in \cite{GJN2015}, near-optimal pairwise tests deciding on pairs of convex hypotheses in simple o.s.'s outlined in Section \ref{sectpairs} can be used as building blocks when constructing near-optimal tests deciding on multiple convex hypotheses. In the sequel, we  use one of these constructions, namely, as follows.
\par
Assume that we are given a simple o.s. $\mycO=((\Omega,P),\{p_\mu:\mu\in\cM\},\cF)$ and two finite collections of nonempty convex compact subsets  $B_1,...,B_b$ (``blue sets'') and $R_1,...,R_r$ (``red sets'')   of $\cM$. Our objective is, given a stationary $K$-repeated observation $\omega^K$ stemming from a distribution $p_\mu$, $\mu\in\cM$, to infer the color of $\mu$, that is, to decide on the hypothesis $H_B:\,\mu\in B:=B_1\cup...\cup B_b$    vs. the alternative $H_R:\,\mu\in R:=R_1\cup...\cup R_r$. To this end we act as follows:
\begin{enumerate}
\item For every pair $i,j$ with $i\leq b$ and $j\leq r$, we solve the problem (\ref{eq11})  with $B_i$ in the role of $M_1$ and $R_j$ in the role of $M_2$; we denote $\Opt_{ij}$ the associated optimal values. The corresponding optimal solutions $\mu_{ij}$ and $\nu_{ij}$ give rise to the detectors
\begin{equation}\label{eqij}
\begin{array}{c}
\phi_{ij}(\omega)=\half \ln\left(p_{\mu_{ij}}(\omega)/p_{\nu_{ij}}(\omega)\right):\Omega\to\bR,\,\,
\phi_{ij}^{(K)}(\omega^K)=\sum_{t=1}^K\phi_{ij}(\omega_t):\Omega^K\to\bR\\
\end{array}
\end{equation}
(cf. (\ref{eq2})) and {\sl risks}
\begin{equation}\label{epsij}
\epsilon_{ij}=\exp(\Opt_{ij})=\int_\Omega \sqrt{p_{\mu_{ij}}(\omega)p_{\nu_{ij}}(\omega)}P(d\omega).
\end{equation}
\item We build the entrywise positive $b\times r$ matrix $E^{(K)}=[\epsilon_{ij}^K]_{{1\leq i\leq b\atop 1\leq j\leq r}}$ and symmetric entrywise nonnegative $(b+r)\times (b+r)$ matrix $E_K=\hbox{\tiny$\left[\begin{array}{c|c}&E^{(K)}\cr\hline [E^{(K)}]^T&\cr\end{array}\right]$}$. Let $\epsilon_{K}$ be the spectral norm of the matrix $E^{(K)}$ (equivalently, spectral norm of $E_K$), and let $e=[g;h]$\footnote{We use ``Matlab notation'' $[a;b]$ for vertical and $[a,b]$ for horizontal concatenation of matrices $a, b$ of appropriate dimensions.} be the Perron-Frobenius eigenvector of $E_K$, so that $e$ is a nontrivial nonnegative vector such that $E_Ke=\epsilon_Ke$. Note that from entrywise positivity of $E^{(K)}$ it immediately follows that $e>0$, so that the quantities
    $$
    \alpha_{ij}=\ln(h_j/g_i),\;\;1\leq i\leq b,\; 1\leq j\leq r
    $$
    are well defined. We set
    \begin{equation}\label{eq12}
    \psi_{ij}^{(K)}(\omega^K)=\phi_{ij}^{(K)}(\omega^K)-\alpha_{ij}=\sum_{t=1}^K\phi_{ij}(\omega_t)-\alpha_{ij}:\Omega^K\to\bR,\,1\leq i\leq b,1\leq j\leq r
\end{equation}
\item Let now $\T^K$ be the test which given observation $\omega^K\in\Omega^K$ with $\omega_t$, $t=1,...,K$, drawn, independently of each other, from a distribution $p_\mu$, claims that $\mu$ is blue (equivalently, $\mu\in B$), if there exists $i\leq b$ such that $\psi_{ij}(\omega^K)\geq0$ for all $j=1,...,r$, and claims that $\mu$ is red (equivalently, $\mu\in R$) otherwise.
\end{enumerate}
The main result about the just described ``color inferring'' test is as follows
\begin{proposition}\label{basicprop} {\rm \cite[Propositions 3.2]{GJN2015}}  Let the components  $\omega_t$ of $\omega^K$ be drawn, independently of each other, from distribution $p_\mu$, $\mu\in B\cup R$. Then the just defined test for every $\omega^K$ assigns $\mu$ with exactly one color, blue or red, depending on the observation. Moreover,
\begin{itemize}
\item when $\mu$ is blue (i.e., $\mu\in B$), the test makes correct inference ``$\mu$ is blue'' with $p_\mu^{K}$-probability at least $1-\epsilon_K$;
\item similarly, when $\mu$ is red (i.e., $\mu\in R$), the test makes correct inference ``$\mu$ is red'' with $p_\mu^{K}$-probability at least $1-\epsilon_K$.
     \end{itemize}
\end{proposition}
Now, suppose that $\overline \T$ is some color inferring test with maximal risk $\leq \epsilon\in(0,\half)$. Obviously, $\overline \T$ gives rise to a straightforward test of hypotheses $H_{B_i}:\mu\in B_i$, $i\leq b$ vs $H_{R_j}:\mu\in R_j$, $j\leq r$ with maximal risk bounded with $\epsilon$. This simple observation implies the following corollary of Proposition \ref{basicprop} (cf. \cite[Proposition 3.4]{GJN2015}).
\begin{proposition}\label{col:basic}
In the just described situation, given $\epsilon\in (0,\half)$, assume that in nature there exists test $\overline{\T}$,
based on
$\overline{K}$-repeated observation $\omega^{\overline K}\sim p_\mu^{\overline K}$ and deciding on blue and red hypotheses, and
such that $\overline \T$ never accepts more than one hypothesis and has risk $\leq\epsilon$, meaning that whenever $\mu\in B$ (whenever
$\mu\in R$), $H_b$ (resp., $H_r$) will be accepted with $p_\mu^{\overline{K}}$- probability $\geq1-\epsilon$. Then risk of detector-based test $\T^K$ utilizing $K$-repeated observation $\omega^K$ does not exceed
$\varepsilon\in(0,1)$ provided that\footnote{The case of unique observation may be of interest when the considered o.s. is Gaussian. The corresponding near-optimality result admits the following reformulation in this case: suppose that in a Gaussian o.s. in nature there exists test $\overline \T$ deciding with risk $\leq \epsilon$ on
hypotheses $H_B$ and $H_R$ using (unique) observation $\omega\sim  \cN(\mu, \bar \sigma^2 I_d)$. Then detector-based coloring inference $\T$ utilizing (unique) observation $\omega$, $\omega\sim  \cN(\mu, \sigma^2 I_d)$
with
\[
{\sigma\leq {q_\N(1-\epsilon)\over q_{\N}\left(1-
\mbox{\small$\varepsilon\over\max[b,r]$}\right)
}
\bar\sigma}
\]
has its risk bounded with $\varepsilon$. Here $q_{\N}(p)$ is the $p$-quantile of $\N(0,1)$: $\Prob_{s\sim\N(0,1)}\{s\leq q_{\N}(p)\}=p$, $0\leq p\leq 1$.
}
\[
K\geq \left\lceil {2\ln\left(\max[b,r]\varepsilon^{-1}\right)\over \ln\left( [4\epsilon(1-\epsilon)]^{-1}\right)} \overline{K}\right\rceil.
\]
\end{proposition}
\section{Adaptive estimation by testing}\label{sec:adest}
\subsection{Estimation over unions of convex sets in simple o.s.: problem setting}\label{sect:setting11}
{\bf Problem setup.}
In the sequel, we deal with the situation as follows. Given are:
\begin{itemize}
\item[1.] simple o.s.  $\mycO=((\Omega,\myP),\{p_\mu(\cdot):\mu\in\cM\},\cF)$,
\item[2.] a collection of $N\geq 2$ convex compact sets $X_j\subset \bR^n$, giving rise to the set $X=\bigcup_{j=1}^NX_j$,
\item[3.] affine mappings $x\mapsto \cA_j(x)$ such that $\cA_j(X_j)\subset\cM$, $j=1,...,N$,
\item[4.] a seminorm $\|\cdot\|$ on $\bR^n$,
\item[5.] reliability tolerance $\epsilon\in(0,1/2)$,
\end{itemize}
\paragraph{Risks.}
Given a nonempty subset $\cJ=\{j_1<...<j_s\}$ of $\{1,2,...,N\}$, set $Y\subset \bigcup_{j=1}^N {X}_j$ and $\varepsilon\in(0,1)$, we define the $\varepsilon$-risk of an $M$-observation estimate $\widehat{x}(\omega^M):\,\Omega^M\to\bR^n$ on $Y$ as
$$
\Risk^{\cJ}_{\varepsilon,M}[\widehat{x}|Y]=\min\left\{\rho: \,\Prob_{\omega^M\sim p^M_{\cA_j(x)} }
\left\{\|\widehat{x}(\omega^M)-x\|>\rho\right\}\leq\varepsilon\,\;\forall (j\in \cJ,x\in Y\cap {X}_j)\right\},
$$
and the associated minimax risk as
\[
\RiskOpt^{\cJ}_{\varepsilon,M}[Y]
=\inf_{\widehat{x}(\cdot)}\Risk^\cJ_{\varepsilon,M}[\widehat{x}|Y]
\]
where the infimum  is taken over all estimates utilizing $M$-repeated observation $\omega^M$.

We assume that, in addition to the above setup,  we are given
 \begin{itemize}
 \item[6.]
   positive integers $\ov K$ and $K$ such that  $M=\ov K+K$ and  $N$ ``preliminary'' $\ov K$-observation estimates $\widetilde{x}_i(\cdot):\,\Omega^{\ov K}\to\bR^n$, along with reals $\myr_i=\myr^{\ov K}_i(\epsilon),\, 1\leq i\leq N$---upper bounds for the partial $\epsilon$-risks of $\widetilde{x}_i(\cdot)$:
\begin{equation}\label{tolya1}
\Risk^{\{i\}}_{\epsilon,{\ov K}}[\widetilde{x}_i|X_i]\leq \myr^{\ov K}_i(\epsilon), \;\;1\leq i\leq N.
\end{equation}
\end{itemize}
\paragraph{Goal and strategy.} Assume that we are given $M$ independent across $k$ observations
\[
\omega_k\sim p_{\cA_{\ell_*}(x_*)},\;1\leq k\leq M
\]
(using the terminology of Section \ref{sec:sos}---a stationary $M$-repeated observation $\omega^M=(\omega_1,...,\omega_M)$),
stemming from an unknown pair $(\ell_*,x_*)$ with $1\leq \ell_*\leq N$ and $x_*\in {X}_{\ell_*}$. Our goal is to build an estimate $\wh x$ of $x_*$ with the least possible risk. To this end
we intend to use collection $\omega^{\ov K}=(\omega_1,...,\omega_{\ov K})$ of the first ${\ov K}$  observations to compute points
$${x}_i=\widetilde{x}_i(\omega^{\ov K}).$$
 Our goal is to use the remaining---secondary---$K$ observations
$\omega^K=(\omega_{{\ov K}+1},...,\omega_{{\ov K}+K})$ to ``aggregate'' these points into an estimate $\wh x$ of $x_*\in X$.
We are going to achieve this goal via techniques for convex hypothesis testing developed in \cite{GJN2015,juditsky2020statistical}.
\paragraph{Notational conventions.}   We denote by $\cO$ and $\cU$ the sets of all ordered pairs $(i,j)$ (resp., unordered pairs $\{i,j\}$)  with $1\leq i,j\leq N$ and $j\neq i$.\par
In the sequel we fix $\ell_*$ and $x_*\in {X}_{\ell_*}$ and, in accordance with what was said {above}, deal with repeated observations with i.i.d. components $\omega_k\sim p_{\cA_{\ell_*}(x_*)}$. We denote by $\wt \Omega^{\ov K}$ the set of all realizations of the ``preliminary'' (pilot) observation $\omega^{\ov K}$ such that
\begin{equation}\label{neweq555}
\|x_*-\widetilde{x}_{\ell_*}(\omega^{\ov K})\|\leq \myr_{\ell_*}:=\myr^{\ov K}_{\ell_*}.
\end{equation}  Due to (\ref{tolya1}) the $p^L_{\cA_{\ell_*}(x_*)}$-probability of $\wt\Omega^{\ov K}$ is at least $1-\epsilon$.
\paragraph{Note:} {\sl From now on we fix a realization $\wt\omega^{\ov K}\in{\wt\Omega^{\ov K}}$ of the preliminary observation $\omega^{\ov K}$};
 in what follows,   $\omega^K$ is the secondary (post-pilot) $K$-repeated observation, $\omega^K=(\omega_{\overline{K}+1},...,\omega_{\overline{K}+K})$. For notational convenience, in the sequel, we suppress explicit reference to $\wt\omega^{\ov K}$ when defining/denoting subsequent entities which in fact depend on
$\wt\omega^{\ov K}$ as parameter.

\subsection{Case of general seminorm}\label{sec:genest}
{\subsubsection{Construction}}
The aggregation routine is as follows.
\begin{enumerate}
\item For $1\leq i\neq j\leq N$ we put
\be
\begin{aligned}
x_i&=\widetilde{x}_i(\wt\omega^{\ov K}),\label{emeq101aa}\\
B_i&=B_i(\wt\omega^{\ov K})=\big\{x\in X_i:\,\|x-x_i\|\leq \myr_i:=\myr_i^{\overline{K}}(\epsilon)\big\},\label{emeq101ab}\\
\delta_{ij}&=\delta_{ij}(\wt\omega^{\ov K})=\half \min\limits_{x\in B_i,y\in B_j}\|x-y\|,
\end{aligned}
\ee{emeq101a}
with the standard convention that minimum over an empty set is $+\infty$.

We specify hypotheses $H_i=H_i\big(B_i(\wt\omega^{\ov K})\big)$ ``the  observations stem from a pair $(i,x)$ with $x\in\ B_i(\wt\omega^{\ov K})$'' (equivalently: $H_i$ states that the distribution of independent across ${k} \leq K$ observations $\omega_{{\ov K}+k}$  belongs to the set $M_i=\{\cA_i(x):\,x\in B_i\}$). Note that sets $M_i=M_i(\widetilde{\omega}^{\ov K})$  are convex and compact subsets of $\cM$.\\
\paragraph{Note:}
Everywhere in the sequel we assume w.l.o.g. that all hypotheses $H_i,\, i=1,...,N,$ are nonempty (i.e., from the start, we reject all empty hypotheses and update accordingly $N$ and the indexes of remaining points $x_i$ and sets $X_j$.
\par
Given a pair {$(i,j)\in\cO$,}
it may happen that there is a simple detector-based {$K$-observation} test $\T_{\{i,j\}}$ as built in Section \ref{sectpairs}, which decides on $H_i$ vs $H_j$ with risk bounded with $\epsilon/(N-1)$; in such case, we say that {pairs $(i,j)$ and $(j,i)$} are {$K$-}good, and say that these pairs are $K$-bad otherwise. We  skip the prefix ``$K$-'' when the value of $K$ is clear from the context.

\item Let for $i\leq N$ $\cJ_i$ be the set of $j\leq N$, $j\neq i$,
such that the pair ${(i,j)}$ is good; note that $j\in\cJ_i$ if and only if $i\in\cJ_j$. For all $i\leq N$ and $j\in \cJ_i$ we run tests $\T_{\{i,j\}}$.
We call index $i$ admissible if hypothesis $H_i$ was never rejected by corresponding tests (i.e., all tests $\T_{\{i,j\}}$ (if any) with  $j\in{\cJ_i}$ accepted $H_i$; in particular, $i$ is admissible, if no pair $(i,j)$ with $j\neq i$ is good). We denote $\cI=\cI(\omega^K)$ the set of all admissible $i$'s.

The output of the procedure---the aggregated estimate $\wh x=\wh x(\omega^K)$---is selected as $x_{\wh i}$ where $\wh i{=\widehat{i}(\omega^K)}$ is the smallest of admissible $i$'s when set $\cI$ is not empty, and selected as, say, $x_1$ otherwise.
\end{enumerate}
We have the following straightforward bound for the error of $\wh x$.
\begin{proposition}\label{ai1}
In the situation described in Section \ref{sect:setting11},  let $\overline{\Omega}^K$ be the set of all $\omega^K$ satisfying the
condition
\begin{quote}
All tests $\T_{\{\ell_*,j\}}$  in good pairs $(\ell_*,j)$, as applied to observation $\omega^K$, accept the hypothesis $H_{\ell_*}$.
 \end{quote}
 Then $(\ell_*,x_*)$-probability\footnote{From now on, for $j\leq N$ and $x\in X_j$ ``$(j,x)$-probability'' of an event is its $p_{\cA_j(x)}^{K}$-probability.}   of $\overline\Omega^K$ is at least $1-\epsilon$, and for all $\omega^K\in \overline\Omega^K$ the set $\bigcup_{i\in \I} B_i$ covers $x_*$. Furthermore, for such $\omega^K$ one has
\be
\|x_*-\wh{x}(\omega^K)\|\leq \|x_*-x_{\ell_*}\|+\max\limits_{j\in \cI^-_{\ell_*}}\|x_j-x_{\ell_*}\|,\quad\quad
\cI^-_{\ell_*}=\{j\in\cI,j<\ell_*\}
\ee{1stb}
(by convention, the maximum over an empty set is zero). Moreover,
\be
\begin{array}{rcl}
\max\limits_{j\in \cI^-_{\ell_*}}\|x_j-x_{\ell_*}\|&\leq& \max\limits_{j\in J^-_{\ell_*}}\|x_j-x_{\ell_*}\|\leq
\myr_{\ell_*}+\max\limits_{j\in J^-_{\ell_*}}(2\delta_{\ell_*j}+\myr_j),\\
J^-_{\ell_*}&=&\{j<\ell_*:\hbox{ $(\ell_*,j)$ is {K-}bad}\},
\end{array}
\ee{onbad-}
 and
\be
\max\limits_{j\in \cI^-_{\ell_*}}\|x_j-x_{\ell_*}\|&\leq& \max\limits_{j\in \cI}\|x_j-x_{\ell_*}\|\leq \max\limits_{j\in J_{\ell_*}}\|x_j-x_{\ell_*}\|
\leq\max\limits_{(i,j)\in \overline{J}}\|x_j-x_i\|\nn &\leq& 2\max\limits_i\myr_{i}+2\max\limits_{(i,j)\in \overline J}\delta_{ij}
\ee{onbad}
\[J_{\ell_*}=\{j\neq \ell_*: \hbox{$(\ell_*, j)$ is {K-}bad}\},\;\;\overline{J}=\{(i,j){\in\cO:} \hbox{$(i,j)$ is {K-}bad}\}.
\]
\end{proposition}

{\subsubsection{Risk analysis}}\label{sec3.2.2}
Given a pair {$(i,j)\in\cO$} and $\varepsilon\in(0,1/2)$ consider the quantity
\be
\myr^K_{ij}(\varepsilon)=\half \max_{x\in X_i,y\in X_j} \left\{\|x-y\|:\,\varrho(\cA_i(x),\cA_j(y))\geq\varepsilon^{1/K}
\right\}
\ee{Deltaij}
where $\varrho(\cdot,\cdot)$ is as defined in \rf{eq1} (here, as before, the maximum over an empty set is 0, by definition). In what follows we refer to $\myr^{K}_{ij}(\varepsilon)$ as {\em separation  {$\varepsilon$-}risk} over $X_i$, $X_j$.
\begin{theorem}\label{thm:genthm}
In the situation described in Section \ref{sect:setting11}, the just built adaptive estimate $\wh x^{(a)}$ (as function of pilot
observation $\omega^\ov{K}$ and independent (secondary) observation $\omega^K$) satisfies
\be
{\Risk^{\{i\}}_{2\epsilon, \ov{K}+K}}[\wh x^{(a)}|X_i]\leq 2\myr_{i}^{\ov{K}}(\epsilon)+\max_{j< i}\left[\myr_j^\ov{K}(\epsilon)+2\myr_{ij}^K\big({\epsilon/(N-1)}\big)\right]\quad\forall i\leq N.
\ee{thm1eq1}
Moreover, whenever $K> \ov\vartheta^{-1}\ov K$
where \[
\ov\vartheta:={\ln(4\epsilon(1-\epsilon))\over 2\ln(\epsilon/(N-1))}\leq 1,
\]
one has
\be
{\Risk^{\{i\}}_{2\epsilon, \ov{K}+K}}[\wh x^{(a)}|X_i]
\leq 2\myr_{i}^{\ov{K}}(\epsilon)+\max_{j< i}\left[\myr_j^\ov{K}(\epsilon)+2\RiskOpt^{\{i,j\}}_{\epsilon,\overline K}[X_i\cup X_j]\right]\quad\forall i\leq N.
\ee{thm1eq2}
\item In addition, in the special case where for every pair $i,j$ there exists $x_{ij}\in X_i\cap X_j$ such that $\cA_i(x_{ij})=\cA_j(x_{ij})$ one has for all $K\geq \ov K$ and $i\leq N$:
\be
{\Risk^{\{i\}}_{2\epsilon, \ov{K}+K}}[\wh x^{(a)}|X_i]
&\leq& 2\myr_{i}^{\ov{K}}(\epsilon)+\max_{j< i}\left[\myr_j^\ov{K}(\epsilon)+2\ov \vartheta^{-1}\RiskOpt^{\{i,j\}}_{\epsilon,\overline K}[X_i\cup X_j]\right].
\ee{thm1eq3}

\end{theorem}
Theorem \ref{thm:genthm} has the following straightforward corollary.
\begin{corollary}\label{cor:if}
{Under} the premise of Theorem \ref{thm:genthm}, suppose that upper bounds $\myr_{i}^{\ov{K}}(\epsilon)$ on partial risks of estimates $\wt{x}_i(\omega^\ov{K}) $ are within factor $\theta$ of the respective $\ov{K}$-observation minimax risks, i.e.,
\[
{\RiskOpt^{\{i\}}_{\epsilon,\ov{K}}}[{X}_i]\leq\myr_i^{\ov{K}}(\epsilon)\leq  \theta{\RiskOpt^{\{i\}}_{\epsilon,\ov{K}}}[{X}_i].
\]
Then the risk  of estimate $\widehat{x}^{(a)}$  is within a moderate factor of  the minimax $\overline K$-observation risk. For instance,
whenever $K\geq \ov \vartheta^{-1} \ov K$
 one has
\begin{equation}\label{eq0012}
{\Risk^{\{i\}}_{2\epsilon,\ov K+K}}[\widehat{x}^{(a)}|{X_i}]\leq (2+3\theta)\max_{j\leq i}\RiskOpt^{\{i,j\}}_{\epsilon,\overline{K}}[X_i\cup X_j]\quad\forall i\leq N,
\end{equation}
and
\be
\RiskOpt^{\overline{1,N}}_{2\epsilon,\overline{K}+K}[X]&\leq& {\Risk^{\overline{1,N}}_{2\epsilon,\ov K+K}}[\widehat{x}^{(a)}|{X}]\leq [2+3\theta]\max_{j,i\leq N}\RiskOpt^{\{i,j\}}_{\epsilon,\overline{K}}[X_i\cup X_j]\nn&\leq&(2+3\theta)\RiskOpt^{\overline{1,N}}_{\epsilon,\overline{K}}[X].
\ee{eqmma}
In the case where for every pair $i,j$ there exists $x_{ij}\in X_i\cap X_j$ such that $\cA_i(x_{ij})=\cA_j(x_{ij})$ one has for all $K\geq \ov K$ and $i\leq N$:
\[
{\Risk^{\{i\}}_{2\epsilon, \ov{K}+K}}[\wh x^{(a)}|X_i]
\leq (3\theta+2\ov\vartheta^{-1})\max_{j\leq i\leq N}\RiskOpt^{\{i,j\}}_{\epsilon,\overline K}[X_i\cup X_j],
\]
so that
\be
{\Risk^{\overline{1,N}}_{2\epsilon, 2\ov{K}}}[\wh x^{(a)}|X]\leq
\max_{i,j\leq N}(3\theta+2\ov\vartheta^{-1})\RiskOpt^{\{i,j\}}_{\epsilon,\overline K}[X_i\cup X_j]\leq (3\theta+2\ov\vartheta^{-1}){\RiskOpt^{\overline{1,N}}_{\epsilon,\overline K}}[X].
\ee{maxun}
\end{corollary}

\paragraph{Discussion.}
Bounds \rf{eq0012}, \rf{eqmma} imply  that  under the premise of the corollary, the minimax risk  $\RiskOpt^{\overline{1,N}}_{2\epsilon,\overline{K}+K}[X]$ of estimation over union $X$ of sets $X_i$, $i=1,...,,N$, is similar, modulo  logarithmic factors, to the maximal ``pairwise'' minimax risk $\max_{j,i\leq N}\RiskOpt^{\{i,j\}}_{\epsilon,\overline{K}}[X_i\cup X_j]$ of estimation over {\em pairwise unions} $X_i\cup X_j$ of sets. Furthermore,
the upper bound \rf{eq0012}  on the maximal over $X_i$ risk ${\Risk^{\{i\}}_{2\epsilon,\ov K+K}}[\widehat{x}^{(a)}|{X_i}]$ of adaptive estimate $\wh x^{(a)}$ is also similar, in the same sense, to the maximal risk
$\max_{j\leq i}\RiskOpt^{\{i,j\}}_{\epsilon,\overline{K}}[X_i\cup X_j]$ of estimation over {\em pairwise unions} $X_j\cup X_i$ with $j\leq i$ and depends on the selected ordering of $X_i$'s. In particular, when this order is chosen so that partial risks of estimation over $X_i$ satisfy
\[
\myr^{\ov K}_1(\epsilon)\leq \myr^{\ov K}_2(\epsilon)\leq ....\leq \myr^{\ov K}_N(\epsilon)
\]
and pairwise separation risks are dominated by partial risks, i.e.,
\be
\myr^K_{ij}\big({\epsilon/(N-1)}\big)\leq C\myr^{\ov K}_i(\epsilon) 
\quad\forall (i,j,\;1\leq j<i\leq N),
\ee{ifadapt}
one has
\[
{\Risk^{\{i\}}_{2\epsilon, \ov{K}+K}}[\wh x^{(a)}|X_i]\leq C'\myr_i^{\ov K}(\epsilon)\quad \forall i\leq N,
\]
and estimate $\wh x^{(a)}$ is adaptive in the sense of \cite{Lepski1991,Lepski1991a}. On the other hand, when relations \rf{ifadapt} do not hold, adaptation in the above sense is impossible what can be seen already when $N=2$. Similar comments are applicable to bound \rf{maxun}.
\subsection{Estimate aggregation, case of Euclidean seminorm}\label{sec:l2est}
We continue to consider  the situation described in Section \ref{sect:setting11}. However, from now on we assume that $\|\cdot\|$ is a Euclidean seminorm such that
$\|x\|=\|Bx\|_2$ where $B\in \bR^{\nu\times n}$ is a given matrix. We build an adaptive estimate of the signal $x_*\in X_{\ell_*}$ underlying our observations: $\omega_k\sim p_{\cA_{\ell_*}(x_*)}$ by aggregating preliminary estimates ---selecting the closest to $x_*$ point among $x_i=\widetilde{x}_i^{\overline{K}}(\widetilde{\omega}^{\overline{K}})$, $1\leq i\leq N$, where, same as before, $\wt\omega^{\ov K}\in{\wt\Omega^{\ov K}}$  is fixed.
\subsubsection{Construction}
We are given the number $K$ of observations and tolerance parameters $\epsilon\in(0,1)$ and  $\underline{\delta}>0$; we put $\overline{N}=2N(N-1)$.
\paragraph{Preliminaries}~\\
$\bullet$ Denote $W_i=B X_i$, $i=1,...,N$, with $W=BX$. Assuming, for the sake of simplicity, that all points $w_i=Bx_i$, $i=1,...,N$, are distinct, we associate with each pair ${(i,j)}\in \cO$ the quantities
$$
r_{ij}=\half\|w_i-w_j\|_2,
$$
vectors $\psi_{ij}={w_j-w_i\over \|w_j-w_i\|_2}$, $w_{ij}=\half (w_i+w_j)$, and for $\delta > 0$ 
consider sets
$$
W^{\ell-}_{ij}=\{v\in W_\ell: \,\psi_{ij}^T(v-w_{ij})\leq 0\},\;\;
W_{ij}^{\ell+}(\delta)=\{v\in W_\ell: \psi_{ij}^T(v-w_{ij})\geq \delta\},\;\ell=1,...,N.
$$
\noindent Observe that $W_{ij}^{\ell-}$ is exactly the set of $v\in W_\ell$ such that $\|v-w_i\|_2\leq \|v-w_j\|_2$, while  $W_{ji}^{\ell+}(\delta)$ is the set of $v\in W_\ell$ such that \[\|v-w_j\|_2^2\leq \|v-w_i\|_2^2-2\delta \|w_i-w_j\|_2.\]
$\bullet$ Let us fix a quadruple $(i,j;\ell,\ell')$, $1\leq i\neq j\leq N$ and $1\leq \ell,\ell'\leq N$.
We denote $H^{\ell-}_{ij}$
(resp., $H^{\ell'+}_{ij}(\delta)$)
the hypothesis stating that observation $\omega^K$ stems from $(\ell, x)$ with $x\in X_\ell$ such that  $w=Bx\in W_{ij}^{\ell-}$ (resp., such that observation $\omega^K$ stems from $(\ell', x)$ with $x\in X_{\ell'}$ and $w\in W_{ij}^{\ell'+}(\delta)$).
 We say that $\delta\in(0, r_{ij}]$ is $(i,j;\ell,\ell')$-good
 if there exists a detector-based test $\T^{\ell\ell'}_{ij}$ deciding on hypothesis $H^{\ell-}_{ij}$ vs $H^{\ell'+}_{ij}(\delta)$ with risk $\leq\varepsilon= {\epsilon {/\ov N}}$. When good $\delta$'s exist, we say that the quadruple $(i,j;\ell,\ell')$ itself is ($\varepsilon$-)good.
  It is obvious that if $\delta'\in[0, r_{ij}]$ is good, so are all $\delta\in [\delta',r_{ij}]$. Note that goodness of $\delta$ can be checked  efficiently, i.e.,  when $(i,j;\ell,\ell')$ is 
 good one can efficiently find, e.g., by bisection, the value  $\delta^{\ell\ell'}_{ij}$ such that $\delta^{\ell\ell'}_{ij}$ is 
 good while $\delta^{\ell\ell'}_{ij}-\underline\delta$ is not.
When $\delta=r_{ij}$ is not $(i,j;\ell,\ell')$-good we say that the corresponding quadruple is bad and set $\delta^{\ell\ell'}_{ij}=r_{ij}$.
\paragraph{Aggregation procedure}~\\
The output of the procedure are two aggregated  estimates $\wh x$ and $\wt x$.
%
\begin{enumerate}
\item For each $(i,j;\ell,\ell')$, $1\leq i\neq j\leq N$ and $1\leq \ell,\ell'\leq N$, we act as follows:\\
$\bullet$ we reject the alternative {$H_{ij}^{\ell'+}(r_{ij})$} if the quadruple in question is bad;\\
$\bullet$ when $(i,j;\ell,\ell')$ is good we apply to observation $\omega^K$ test $\T^{\ell\ell'}_{ij}$ of hypothesis $H_{ij}^{\ell-}$ against $H_{ij}^{\ell\ell'+}={H_{ij}^{\ell'+}}(\delta^{\ell\ell'}_{ij})$.

We say that pair $(i;\ell)$ is {\em admissible} if  corresponding hypotheses $H^{\ell-}_{ij}$ were never rejected by the above procedure. The result of this step is the set $\I=\I(\omega^K)$ of all admissible pairs $(i;\ell)$.
\item If $\I(\omega^K)=\emptyset$  we select the aggregated solution 
as one of $x_i$, e.g., $\wh x=x_1$; when $\I(\omega^K)$ contains pairs corresponding to a unique index $\wh i=\wh i(\omega^K)$, we output $\wh x(\omega^K)=x_{\wh i}$ as aggregated solution. Otherwise,\\
$\bullet$  we select $\wh i=\wh i(\omega^K)$ as (e.g., the smallest) $i$-component corresponding to admissible pairs $(i;\ell)$ with the smallest value of (the second index) $\ell$ and define the estimate $\wh x(\omega^K)=x_{\wh i}$.

$\bullet$ To build the estimate $\wt x$ we find
among $w_i$ corresponding to admissible $i$'s (that is, $i$-components of admissible pairs $(i;\ell)$) points $w_{\bar i},\,w_{\bar j}$ with the maximal length $\|w_{\bar i}-w_{\bar j}\|_2$ of the connecting segment and select as aggregated solution $\wt x(\omega^K)=\half (x_{\bar i}+x_{\bar j})$ (or choose any $\wt x\in \bR^n$ such that $B \wt x=w_{\bar i\bar j}$). 
\end{enumerate}

\begin{proposition}\label{prop:l2aggreg1}
Suppose  that  observation $\omega^K$ stems from the pair $(\ell_*, x_*), \,x_*\in X_{\ell_*}$. Let ${i_*}$ be the index of one of the $\|\cdot\|$-closest to $x$ points among $x_1,...,x_N$, and let $\overline\Omega^K$ be the set of realizations $\omega^K$  such that
as applied to $\omega^K$, all tests $\T^{\ell_*\ell}_{i_*j}$ and $\T^{\ell\ell_*}_{ji_*}$ accept the true, if any, of the hypotheses from the corresponding pair.\footnotemark
Then the $(\ell_*,x_*)$-probability of $\overline\Omega^K$ is at least $1-\epsilon$, and for all $\omega^K \in \overline\Omega^K$
it holds 
\be
\|\wh x-x_*\|\leq\|x_*-x_{i_*}\|+2\wh{\delta}_{i_*}(\omega^K)
\ee{ea0}
where 
$\wh{\delta}^{\ell_*}_{i_*}(\omega^K)=\max\limits_{j\neq i_*,\ell\leq \ell_*,(i;\ell)\in \cI(\omega^K)}\delta^{\ell\ell_*}_{ji_*}$.
Furthermore, one has
\be
\|\wt{x}-{x_*}\|^2\leq\|x_*-x_{i_*}\|^2+4\wt{\delta}_{i_*}^{\ell_*}(\omega^K)^2
\ee{stupid0}
where $\wt{\delta}^{\ell_*}_{i_*}(\omega^K)=\max\limits_{j\neq i_*,(j;\ell)\in \cI(\omega^K)}\delta^{\ell\ell_*}_{ji_*}$.
\end{proposition}
\footnotetext{In other words, as applied to $\omega^K$, test $\T^{\ell_*\ell}_{i_*j}$ accepts $H^{\ell_*-}_{i_*j}$ (recall that $H^{\ell_*-}_{i_*j}$ is the ``true hypothesis'' in this case), while test $\T^{\ell\ell_*}_{ji*}$ rejects $H^{\ell-}_{ji_*}$ and accepts $H^{\ell\ell_*+}_{ji_*}$ if $w\in W^{\ell+}_{ji_*}(\delta^{\ell\ell_*}_{ji_*})$.}

In statistical literature, the bound \rf{stupid0} for prediction loss (in Problem I discussed in the introduction, this corresponds to the seminorm $\|x\|=\|Ax\|_2$) is typically obtained utilizing exponential weights (see, e.g.,  \cite{audibert2004aggregated,bunea2007aggregation,juditsky2008learning,tsybakov2003optimal}) or $Q$-aggregation \cite{dai2012deviation,lecue2014optimal}. When risk of aggregation is measured by a Euclidean seminorm, using the aggregation procedure described above this type of results can be painlessly extended to aggregation problems with convex constraints on unknown signals and aggregation from indirect observations.
\subsubsection{Risk analysis}
\begin{theorem}\label{thm:genthm2}
In the situation of this section, estimate $\wh x^{(a)}(\omega^{\ov K+K})=\wh x(\omega^K)$ satisfies for all $i\leq N$:
\be
{\Risk^{\{i\}}_{2\epsilon, \ov{K}+K}}[\wh x^{(a)}|X_i]\leq \myr_{i}^{\ov{K}}(\epsilon)+2\left[\max_{j< i}\myr_{ij}^K\big({\epsilon{/\ov N}}\big)+\underline \delta\right].
\ee{thm1eq121}
Furthermore, for $\wt x^{(a)}(\omega^{\ov K+K})=\wt x(\omega^K)$ one has
\be
{\left[\Risk^{\{i\}}_{2\epsilon, \ov{K}+K}[\wt x^{(a)}|X_i]\right]^2}\leq [\myr_{i}^{\ov{K}}(\epsilon)]^2+4\left[\max_{j\neq i,j\leq N}\myr_{ij}^K\big({\epsilon{/\ov N}}\big)+\underline \delta\right]^2 \qquad\forall i\leq N
\ee{thm1eq122}
(as before, quantities  $\myr_{ij}^K(\varepsilon)$ are given by  \rf{Deltaij}).
\item
Consequently, when $K> \ov\vartheta^{-1}\ov K$
where \[
\ov\vartheta:={\ln([4\epsilon(1-\epsilon)])\over 2\ln(\epsilon{/\ov N})}\leq 1,
\]
one has
\be
{\Risk^{\{i\}}_{2\epsilon, \ov{K}+K}}[\wh x^{(a)}|X_i]
\leq \myr_{i}^{\ov{K}}(\epsilon)+2\max_{j< i}\RiskOpt^{\{i,j\}}_{\epsilon,\overline K}[X_i\cup X_j]+\underline \delta\quad\forall i\leq N,
\ee{thm1eq222}
and
\be
{\left[\Risk^{\{i\}}_{2\epsilon, \ov{K}+K}[\wt x^{(a)}|X_i]\right]^2}\leq [\myr_{i}^{\ov{K}}(\epsilon)]^2+4\left[\max_{j< i}\RiskOpt^{\{i,j\}}_{\epsilon,\overline K}[X_i\cup X_j]+\underline \delta\right]^2 \qquad\forall i\leq N.
\ee{thm1eq1223}
\item In the case where for every pair $i,j$ there exists $x_{ij}\in X_i\cap X_j$ such that $\cA_i(x_{ij})=\cA_j(x_{ij})$ one has for all $K\geq \ov K$ and $i\leq N$:
\be
{\Risk^{\{i\}}_{2\epsilon, \ov{K}+K}[\wh x^{(a)}|X_i]}
&\leq& \myr_{i}^{\ov{K}}(\epsilon)+2\ov \vartheta^{-1}\max_{j< i}\RiskOpt^{\{i,j\}}_{\epsilon,\overline K}[X_i\cup X_j]+\underline \delta
\ee{thm1eq231}
and
\be
{\left[\Risk^{\{i\}}_{2\epsilon, \ov{K}+K}[\wt x^{(a)}|X_i]\right]^2}\leq [\myr_{i}^{\ov{K}}(\epsilon)]^2+4\left[\ov \vartheta^{-1}\max_{j< i}\RiskOpt^{\{i,j\}}_{\epsilon,\overline K}[X_i\cup X_j]+\underline \delta\right]^2 \qquad\forall i\leq N.
\ee{thm1eq232}
\end{theorem}

\section{"Generic" test-based aggregation}\label{sect:general}
\subsection{Setup}\label{setup} We consider the situation as follows: we are given
\begin{itemize}
\item {\sl observation space} $\Xi$,
\item a compact set $X\subset \bR^n$, with every $x\in X$ associated with a family $\cP_x$ of probability distributions on $\Xi$; we refer to observations distributed according to $P\in\cP_x$ as to observations
 {\sl stemming} from $x$,
\item a seminorm $\|\cdot\|$ on $\bR^n$,
\item $N$ distinct points $x_i\in\bR^n,\,i=1,...,N$.
\end{itemize}
Given observation $\xi\sim P$ stemming from unknown signal $x\in X$,  our objective is to aggregate $x_i$'s---to find among $x_i$'s a point which is ``as close to $x$ as the closest to $x$ point among $x_i$'s.'' Here closeness is measured by the seminorm $\|\cdot\|$. \par
Note that as far as our goal is concerned, we lose nothing when assuming from now on that $\|x_i-x_j\|>0$ whenever $i\neq j$.
\paragraph{Conventions} \begin{itemize}\item
In the sequel we say that an event (a set in the space of observations)
takes place with $x$-probability at most (or at least) $p$ for some $x\in X$ if this is the case for probability w.r.t. any distribution from $\cP_x$.
\item  With a subset $Y$ of $X$ we associate {\sl hypothesis} $H(Y)$ on the distribution of observation $\xi$; the hypothesis states that the observation stems from a signal $x\in Y$. Given $Y^1,Y^2\subset X$, a {\sl test} for the pair of hypotheses $H(Y^1),H(Y^2)$ is a procedure which, given on input an observation $\xi$, {\sl accepts} exactly one of these two hypotheses (informally: claims that $\xi$ is drawn from a distribution obeying the accepted hypothesis)  and rejects the other. We say that such a  test has {\sl risk} $\leq\delta$, if ``the probability to accept the true hypothesis is at least $1-\delta$,'' specifically, for $\chi=1,2$,  when the observation stems from a signal $x\in Y^\chi$, the $x$-probability for the test to accept $H(Y^\chi)$ is at least $1-\delta$. Note that we allow for $Y_1$, or $Y_2$, or both, to be empty; whenever this is the case, the test which always accepts a nonempty hypothesis, if any, and accepts a whatever one of the hypotheses when both $Y_1$ and $Y_2$ are empty, has zero risk.
\par
\item  For $(i,j)\in\cO$ and $\delta\geq 0$ we set
\begin{equation}\label{eqrhos}
\rho_i=\min_{x\in X}\|x-x_i\|,\,1\leq i\leq N,
\end{equation}
and
\begin{equation}\label{eqXij}\\
\begin{array}{rcl}
r_{ij}=\half \|x_i-x_j\|,\;\;
X_{ij}(\delta)=\{z\in X:\,\|z-x_i\|\leq r_{ij}-\delta\}.
\end{array}
\end{equation}
Note that $r_{ij}=r_{ji}$.
\end{itemize}
\subsection{Aggregation in general seminorm}\label{sec:basic}
\subsubsection{The setup}
The setup of the general aggregation scheme is given by
\begin{enumerate}
\item ``reliability tolerance'' $\epsilon\in(0,1)$,
\item a collection $\cC$ of pairs $\{i,j\}\in\cU$ with each pair $\{i,j\}\in \cC$ associated with
    \begin{itemize}
    \item thresholds $\Delta_{ij}=\Delta_{ji}\in[0,r_{ij}]$, giving rise to sets $X_{ij}(\Delta_{ij})$, $X_{ji}(\Delta_{ij})$ and hypotheses $H_{ij}=H(X_{ij}(\Delta_{ij}))$, $H_{ji}=H(X_{ji}(\Delta_{ij}))$, along with
    \item a test $\T_{\{i,j\}}$ deciding on the hypotheses $H_{ij}$ and $H_{ji}$ with risk $\leq \epsilon/(N-1)$.
    \end{itemize}
      When $\{i,j\}\in\cC$, we say that $i$ and $j$ are comparable (same as ``$j$ is comparable to $i$'' and ``$i$ is comparable to $j$'').
      \item For pairs $\{i,j\}\in \cU$ with incomparable to each other $i$ and $j$ (i.e., $\{i,j\}\in\cU\backslash \cC$) we set
$$
\Delta_{ij}=\Delta_{ji}=\max[0,r_{ij}-\max[\rho_i,\rho_j]].
$$

\end{enumerate}
\subsubsection{Aggregation routine}
Aggregation routine associated with the  just described setup is as follows
\begin{enumerate}
\item Given observation $\xi$, for every pair $(i,j)\in\cO$ we ``compare $i$ to $j$'' according to the following rule:
\begin{itemize}
\item when $i,j$ are comparable, we run the test $\T_{\{i,j\}}$ on observation $\xi$ and claim that $i$ looses to $j$ when the test accepts the hypothesis $H_{ji}$, and  claim that $j$ looses to $i$ otherwise
\item when $i,j$ are incomparable, $i$ looses to $j$ whenever $\rho_j\leq\rho_i$, otherwise $j$ looses to $i$.
\end{itemize}
\item For $i\leq N$, we denote by $\cI_i=\cI_i(\xi)$ the set of indices $j\neq i$ such that $i$ looses to $j$, set
\[
d_i(\xi)=\max_{j\in \cI_i(\xi)}\|x_j-x_i\|,
\]
and define the aggregated estimate as $\wh{x}(\xi)=x_{\wh{i}(\xi)}$ where
\[\widehat{i}=\widehat{i}(\xi)\in\Argmin_id_i(\xi).
\]
\end{enumerate}
We have the following simple statement.
\begin{proposition}\label{pro:3p}
Suppose  that the observation stems from a signal $x_*\in X$, and let $x_{i_*}$ be one of the $\|\cdot\|$-closest to ${x_*}$ points among $x_1,...,x_N$.
Let $\overline\Xi$ be the set of realizations $\xi$ satisfying the following condition:
\begin{quote}
 For every $j\neq i_*$  such that $i_*$ and $j$ are comparable and ${x_*}\in X_{i_*j}(\Delta_{i_*j})$,  test
$\T_{\{i_*,j\}}$  as applied to observation $\xi$ accepts  the hypothesis $H_{i_*j}$.
 \end{quote}
 Then the ${x_*}$-probability of $\overline\Xi$ is at least $1-\epsilon$, and for all $\xi\in \overline\Xi$
\be
\|{x_*}-\wh{x}(\xi)\|\leq 3 \|x_{i_*}-{x_*}\|+2\overline{\Delta}_{i_*}(\xi)
\ee{bound0}
where
\[\overline{\Delta}_{i_*}(\xi)={\left\{\begin{array}{ll} 0,&\cI_{i_*}(\xi)=\emptyset,\\
\max_{j\in \cI_{i_*}(\xi)}\Delta_{i_*j},&\cI_{i_*}(\xi)\neq\emptyset.
\end{array}\right.}
\]
\end{proposition}
\paragraph{Remarks.} The above construction is inspired by the aggregation procedure of \cite{goldenshluger2009universal} which itself generalizes the results on density estimation with $\ell_1$-loss from \cite{yatracos1985rates,mahalanabis2007density,devroye2012combinatorial}; it can also be seen as a refinement of the selection procedure of \cite[Section 2.5.3]{juditsky2020statistical}.
\par
The question of (near-)optimality of the accuracy bound \rf{bound0} for the proposed routine is more involved in the considered here general framework than in the direct observation setting of \cite{goldenshluger2009universal}; we postpone the corresponding analysis till Section \ref{sec:genorma}.
Note, however, that there are in fact two questions---that of optimality of the factor ``3'' in front of the minimal loss $\|x_{i_*}-x\|$ which is related to problem geometry (and is independent of the observation scheme), and that of the size of the additive term $\overline{\Delta}$.
It appears that in the problem of aggregation of densities {when} $\|\cdot\|$ is the $\ell_1$-norm the factor 3 in front of the minimal error is (in a certain precise sense, cf. \cite{bousquet2019optimal}) unimprovable  even for problems with $N=2$.
On the other hand, when allowing for ``improper aggregation,'' i.e., when removing  the limitation of the aggregated solution to be one of given points \cite{mahalanabis2007density} {supplies} a randomized  algorithm which attains the factor 2 when $N=2$, and factor 2 is, in a certain sense, unimprovable in the latter setting, see  \cite{chan2014near}. However,  known to us attempts to generalize this kind of result to the case of $N>2$ (cf. \cite{bousquet2019optimal}) result in the inflation of the additive term which is too important in the situation of small minimal loss we are mainly interested here.
There is however a situation where the factor 3 can be removed rather painlessly (at the price of a moderate increase of $\overline{\Delta}$)---this is the case of Euclidean seminorm $\|\cdot\|$, and this is the situation we consider next.
\subsection{Aggregation in Euclidean seminorm}\label{sec:euag} Now assume that $\|\cdot\|$ is a Euclidean seminorm: $\|x\|=\|Bx\|_2$ for a given matrix $B$. For $\delta\geq 0$ and $(i,j)\in \cO$ we define
\begin{equation}\label{scriptXij}
\X_{ij}(\delta)=\{z\in X: \|z-x_j\|\geq \delta+\|z-x_i\|\}.
\end{equation} Aggregation procedures presented below are refined versions of the routine from \cite{juditsky2016hypothesis}.
\subsubsection{The setup}
The setup for the Euclidean aggregation is given by
\begin{enumerate}
\item thresholds $\Delta_{ij}$, $(i,j)\in\cO$, such that
$$
\Delta_{ij}=\Delta_{ji}\geq 0
$$
\item tests $\T_{\{i,j\}}$, $\{i,j\}\in\cU$, with $\T_{\{i,j\}}$ testing the hypothesis $\cH_{ij}=H(\X_{ij}(\Delta_{ij}))$ vs. the alternative $\cH_{ji}=H(\X_{ji}(\Delta_{ji}))$ such that
\begin{itemize}
\item if both hypotheses $\cH_{ij}$ and $\cH_{ji}$ are empty, $\T_{\{i,j\}}$ accepts both hypotheses;
\item if exactly one of the hypotheses $\cH_{ij},$ $\cH_{ji}$ is empty, the test always accepts the nonempty hypothesis and rejects the empty one;
\item if hypotheses  $\cH_{ij}$ and $\cH_{ji}$ are nonempty (in this case, we refer to the
pairs $(i,j)$ and $(j,i)$  as {\sl good}) the test accepts exactly one of these hypotheses,  and the risk of the test does not exceed $\epsilon/\overline{N}$, $\overline{N}={N(N-1)\over 2}$.
\end{itemize}
\end{enumerate}
\subsubsection{Aggregation routine}
Aggregation routine associated with the above setup is  as follows: given observation $\xi$, we run tests $\T_{\{i,j\}}$, $\{i,j\}\in\cU$, and for every $1\leq i\leq N$ record the ``score of $i$''---the number $s_i(\xi)$  of those $j\neq i$, $j\leq N$ for which the test $\T_{\{i,j\}}$ rejects $\cH_{ij}$. We put
\[\widehat{i}(\xi)\in \Argmin_{1\leq i\leq N} s_{i}(\xi)
\]
and define aggregated solution {as} $\wh{x}(\xi)=x_{\wh{i}(\xi)}$.
\begin{proposition}\label{prop:l2aggreg10} Suppose that the observation stems from a signal ${x_*}\in X$.
Let $\overline\Xi$ be the set of realizations of $\xi$ such that
\begin{quote}
as applied to observation $\xi$, each test $\T_{\{i,j\}}$ associated with a good pair $(i,j)$  does not reject the ``true hypothesis,'' if any (i.e., as applied to $\xi$, the test  accepts $\cH_{ij}$ when ${x_*}\in \X_{ij}(\Delta_{ij})$, and accepts $\cH_{ji}$  when ${x_*}\in \X_{ji}(\Delta_{ij})$).
\end{quote}
Then the ${x_*}$-probability of $\overline\Xi$ is at least $1-\epsilon$, and for all $\xi \in \overline\Xi$ one has
\[
\|{x_*}-\wh{x}(\xi)\|\leq\|x_{i_*}-{x_*}\|+2\overline{\Delta}
\]
where $x_{i_*}$ is one of the $\|\cdot\|$-closest to ${x_*}$ points among $x_1,...,x_N$ and $\overline{\Delta}=\max_{i\neq j}\Delta_{ij}$.
\end{proposition}

\section{Test-based aggregation in simple observation schemes}\label{sec:simplebs5}
\subsection{Problem setting}\label{sect:setting}
In the sequel, we deal with the situation as follows. Given are:
\begin{enumerate}
\item simple o.s.  $\mycO=((\Omega,\myP),\{p_\mu(\cdot):\mu\in\cM\},\cF)$,
\item a collection of $J$ convex compact sets $X_\nu\subset \bR^n$, giving rise to the set $X=\cup_{\nu=1}^J X_\nu$,
\item affine mappings $x\mapsto \cA_\nu(x)$ such that $\cA_\nu(X_\nu)\subset\cM$, $\nu=1,...,J$,
\item a seminorm $\|\cdot\|$ on $\bR^n$,
\item $N$  points $x_i\in \bR^n,\,i=1,...,N$.
\end{enumerate}
Our objective is given a stationary  repeated observation $\omega^K=(\omega_1,...,\omega_K)$ with
$$
\omega_k\sim p_{\cA_\nu(x)},\,k=1,2,...,K,
$$
for some {\sl unknown} pair $(\nu,x)$ with $\nu\leq J$  and $x\in X_\nu$, to recover one of the $\|\cdot\|$-closest to  $x$ points among $x_1,...,x_N$.
\par
Note that we are in the situation of Section \ref{setup}, with
\begin{itemize}
\item $\Omega^K$ in the role of observation space $\Xi$, and $K$-repeated observation $\omega^K$ in the role of observation $\xi$,
\item the signal set $X=\bigcup_{\nu=1}^J X_j\subset\bR^n$,
\item
the family {$\cP_x:=\cP^K_x$} of probability distributions of observations stemming from a signal $x\in X$ being comprised of all distributions with densities $p_{\cA_\nu(x)}^K(\omega^K)$ and $\nu$ satisfying $x\in X_\nu$, the densities being taken w.r.t. the reference measure $\Pi^K$.
\end{itemize}
Thus, by convention taken in Section \ref{setup}, claim that an event
takes place with $x$-probability at most (or at least) $p$, $x\in X$, means that this is the case for probability w.r.t. any density $p^K_{\cA_\nu(x)}$ of observation $\omega^K$ with $\nu$ satisfying $x\in X_\nu$.
\par
We are about  to achieve our goal of recovering one of the $\|\cdot\|$-closest to $x$ points
among $x_1,...,x_N$  via techniques developed in Section \ref{sect:general}, and in what follows we use terminology and notation from that section.

\subsection{Aggregation in general seminorm}\label{sec:genorma}
Our current objective is to describe an implementation of the aggregation procedure of Section \ref{sec:basic} in the present setting.
\subsubsection{Preliminaries}\label{sec5.2.1}
Given number $K$ of observations and $\epsilon\in(0,1)$, in order to build for $(i,j)\in\cO$ the quantities $\Delta_{ij}$ and tests $\T_{i,j\}}$, as required  by construction from Section \ref{sec:basic},  we act as follows.
\\
$\bullet$ Let us associate with  $\delta\geq0$ and $(i,j)\in\cO$  sets $X_{ij}(\delta)$, $X_{ji}(\delta)$ (see (\ref{eqXij})) and hypotheses
$$
H_{ij}[\delta]=H(X_{ij}(\delta)),\;\; H_{ji}[\delta]=H(X_{ji}(\delta)).$$
\par
$\bullet$
Let a pair $(i,j)\in\cO$ be fixed. Given  $\delta\geq0$ such that $X_{ij}(\delta)\neq \emptyset$ and $X_{ji}(\delta)\neq \emptyset$, or, which is the same, such that
\begin{equation}\label{deltabarij}
0\leq \delta\leq \overline{\delta}^{ij}:=r_{ij}-\max[\rho_i,\rho_j]
\end{equation}
where $\rho_s$ are defined by (\ref{eqrhos}), observe that $X_{ij}(\delta)$ is a finite union of convex compact sets:
\[
X_{ij}(\delta)=\bigcup_{\nu=1}^{J}\Big\{X_{\nu}\cap \{z:\,\|z-x_i\|\leq r_{ij}-\delta\}\Big\}.
\]
We specify collection $\cR_{ij}(\delta)=\{R_{ij}^s(\delta):1\leq s\leq J_{ij}^{\delta\r}\}$ of ``red'' nonempty convex compact sets as the collection of all {\sl nonempty} sets of the form
\[
R_{ij\nu}(\delta)=\{\cA_\nu(x):\,x\in X_\nu,\,\|x-x_i\|\leq r_{ij}-\delta\}, \;1\leq \nu\leq J.
\]
Similarly, we specify the collection $\cB_{ij}(\delta)=\{B_{ij}^s(\delta):1\leq s\leq J_{ij}^{\delta\b}\}$  of ``blue''  nonempty convex compact sets as the collection of all nonempty sets of the form
\[B_{ij\mu}(\delta)=\{\cA_\mu(x):\,x\in X_\mu,\,\|x-x_j\|\leq r_{ij}-\delta\},\;1\leq \mu\leq J.
\]
When applying to $\cR_{ij}(\delta)$ and $\cB_{ij}(\delta)$ the color inferring test from Section \ref{sectinfcol}, depending on $\delta$ it may happen that the risk bound
$\epsilon_K$ of the inference as defined in Section \ref{sectinfcol} satisfies  $\epsilon_K\leq \epsilon/(N-1)$.
Let us refer to $\delta$ as {\sl $(i,j)$-appropriate}, if $0\leq \delta\leq \overline{\delta}^{ij}$ and $\epsilon_K\leq\epsilon/(N-1)$.
\par
$\bullet$ Given $i,j$ and $\delta$ satisfying (\ref{deltabarij}), we can check efficiently whether $\delta$ is $(i,j)$-appropriate---to this end we should compute the spectral norm of a {$[J^{\delta\r}_{ij}+J^{\delta\b}_{ij}]\times[J^{\delta\r}_{ij}+J^{\delta\b}_{ij}]$}  symmetric matrix filled with optimal values of  {$J^{\delta\r}_{ij}J^{\delta\b}_{ij}$} explicit convex optimization problems.
Clearly, if $\delta$ is $(i,j)$-appropriate and $\delta'\in[\delta,\overline{\delta}^{ij}]$, then $\delta'$ is $(i,j)$-appropriate
along with $\delta$. \par
$\bullet$
Let us call $(i,j)$ {\sl appropriate}, if $\overline{\delta}^{ij}$ is nonnegative and $(i,j)$-appropriate.
In this case the infimum $\underline{\delta}^{ij}$ of $(i,j)$-appropriate $\delta\in [0,\overline{\delta}^{ij}]$
is well defined,
and bisection in $\delta$ allows to obtain rapidly $(i,j)$-appropriate upper bounds  on $\underline{\delta}^{ij}$ to whatever high accuracy.
The bottom line is that one can efficiently check whether the pair $(i,j)\in\cO$ is appropriate, and  whenever it is the case, the quantity
$$
\underline{\delta}^{ij}\in\left[0,\overline{\delta}^{ij}=r_{ij}-\max[\rho_i,\rho_j]\right]
$$
is efficiently computable, and whenever
$$
\Delta_{ij}=\Delta_{ij}=\widetilde{\delta}_{ij}\in[0,\overline{\delta}_{ij}]
$$
with an $(i,j)$-appropriate $\widetilde{\delta}_{ij}$, we can point out $K$-observation test $\T_{ij}$ which decides on the hypotheses $H_{ij}(\Delta_{ij})$, $H_{ji}(\Delta_{ji})$ with risk $\leq\epsilon/(N-1)$. Under the circumstances, the latter means that as applied to observation $\omega^K$, test $\T_{ij}$ accepts at most one of the hypotheses $H_{ij}(\Delta_{ij})$, $H_{ji}(\Delta_{ij})$, and
whenever $\omega^k\sim p^K_{A_\nu(x)}$ for $\nu$ and $x$ such that $x\in X_\nu$, the probability for $\T_{ij}$ to accept $H_{ij}(\Delta_{ij})$ is at least $1-\epsilon/(N-1)$ when $\|x-x_i\|\leq r_{ij}-\Delta_{ij}$, and
the probability for $\T_{ij}$ to accept $H_{ji}(\Delta_{ij})$ is at least $1-\epsilon/(N-1)$ when $\|x-x_j\|\leq r_{ij}-\Delta_{ij}$.
{Note that for an appropriate pair $(i,j)$, the above $\widetilde{\delta}_{ij}$ can be made arbitrarily close to $\underline{\delta}_{ij}$. }
\par
$\bullet$
As is immediately seen, a pair $(i,j)$ is appropriate if and only if so is the pair ${j,i}$, and because
$\overline{\delta}_{ij}$ and $\underline{\delta}_{ij}$ are symmetric in $i,j$, $\delta$ is $(i,j)$-appropriate if and only if $\delta$ is ${j,i}$-appropriate, which allows to restrict ourselves to $\widetilde{\delta}_{ij}$ which are symmetric in $i,j$ as well.  Consequently, appropriateness of a pair, appropriateness of a $\delta$ for this pair, and the parameters  $\overline{\delta}_{ij}$, $\underline{\delta}_{ij}$, $\widetilde{\delta}_{ij}$ are attributes of {\sl unordered} pair $\{i,j\}$ rather than of the ordered pair $(i,j)$.  For appropriate $\{i,j\}$, let us set $\bar{i}=\min[i,j]$, $\bar{j}=\max[i,j]$ and $\T_{\{i,j\}}=\T_{\bar{i}\bar{j}}$, so that $\T_{\{i,j\}}$ decides on $H_{ij}$ vs. $H_{ji}$ with risk $\leq\epsilon/(N-1)$.
\subsubsection{Aggregation routine}
Consider the following procedure.
\begin{itemize}
\item We specify the set $\cG$ of all appropriate pairs  $\{i,j\}\in\cU$ along with the related quantities $\widetilde{\delta}_{ij}$ (the smaller the better) and tests $\T_{\{i,j\}}$. Next, we declare a whatever subset $\cC$ of the set $\cG$ to be the  set of comparable pairs of indices as defined in Section \ref{sec:basic}, and set
    $$
    \Delta_{ij}=\left\{\begin{array}{ll}\widetilde{\delta}_{ij},&\{i,j\}\in\cC\\
    \max[0;r_{ij}-\max[\rho_i,\rho_j]],&\hbox{otherwise}\\
    \end{array}\right.
    $$
With the thresholds $\Delta_{ij}$ just defined, $K$-repeated observation $\omega^K$ in the role of $\xi$, and with $\cC$ in the role of the set of comparable pairs and associated tests, we arrive at the aggregation setup as described in Section \ref{sec:basic}, satisfying all the requirements from that section.
\item Given observation $\omega^K$, we apply the aggregation procedure associated with the above setup, resulting in the aggregated estimate $\wh x(\omega^K)$.
\end{itemize}
Results of Propositions \ref{pro:3p} and \ref{basicprop}  imply the following property of the resulting estimate.
\begin{proposition}\label{prop:consbas}
In the situation of this section, suppose that the just described routine is applied to observation $\omega^K$ stemming from ${x_*}\in X$, so that $\omega^K\sim p^K_{\cA_\nu({x_*})}$ for some $\nu\leq J$ such that  ${x_*}\in X_\nu$. Let also $i_*$ be the index of one of the $\|\cdot\|$-closest to ${x_*}$ points among $x_1,...,x_N$. Finally,  let $\overline{\Omega}$ be the set of all $\omega^K$ satisfying the
condition (cf. Proposition \ref{pro:3p})
\begin{quote}
 For every $j\neq i_*$  such that $i_*$ and $j$ are comparable and ${x_*}\in X_{i_*j}(\Delta_{i_*j})$, test $\T_{\{i_*,j\}}$ as applied to observation $\omega^K$ accepts  the hypothesis $H_{i_*j}(\Delta_{i_*j})$
 \end{quote}
Then the $p_{\cA_\nu({x_*})}^K$-probability of $\overline{\Omega}$ is at least $1-\epsilon$, and
the aggregated solution $\wh x(\omega^K)$ satisfies
\[
\omega^K\in\overline{\Omega}\Rightarrow
\|{x_*}-\wh{x}(\omega^K)\|\leq 3 \|x_{i_*}-{x_*}\|+2\overline{\Delta}_{i_*},
\]
where
\[\overline{\Delta}_{i_*}=\left\{\begin{array}{ll}0,&\cI_{i_*}(\omega^K)=\emptyset,\\
\max_{j\in \cI_{i_*}(\omega^K)} \Delta_{i_*j},&\cI_{i_*}(\omega^K)\neq\emptyset.
\end{array}\right.\]
(for notation, see the description of aggregation in Section \ref{sec:basic}).
\end{proposition}

\subsubsection{Characterizing performance}
\begin{theorem}\label{thm:lowerbas} In the setting described in Section \ref{sect:setting}, assume that $x_i\in X$, $1\leq i\leq N$, and that
 for some positive  integer $\overline{K}$, $\epsilon\in(0,1/2)$ and $(\gamma,\delta)\geq0$ there exists inference $\omega^{\overline{K}}\mapsto \bar x(\omega^{\overline{K}})\in \bR^n$ such that
\[
\Prob_{\omega^{\overline{K}}\sim p^{\overline{K}}_{\cA_\nu(x)}}\Big\{\|x-\bar x(\omega^{\overline{K}})\|\leq \gamma\|x-x_{i_*}\|+\delta\Big\}\geq 1-\epsilon\;\;\forall(\nu\leq J,\,x\in X_\nu)
\]
where $x_{i_*}$ is one of the $\|\cdot\|$-closest  to $x$ point among $x_1,...,x_N$. Now let
{$$
\gamma'>\gamma
$$
and let} $K$ satisfy the relation
\begin{equation}\label{eqKl}
K\geq \left\lceil {2\ln\left(J(N-1)/\epsilon\right)\over \ln\left( [4\epsilon(1-\epsilon)]^{-1}\right)} \overline{K}\right\rceil.
\end{equation}
Then, with $K$-repeated observations $\omega^K$, all pairs {$\{i,j\}$} with $r_{ij}>\delta$ are appropriate, and specifying these pairs as comparable, the aggregation procedure described in {this section} {with properly selected $\widetilde{\delta}_{ij}$} ensures that
the resulting aggregated estimate $\wh x(\omega^K)$ satisfies
\be
\Prob_{\omega^{K}\sim p^{K}_{\cA_\nu(x)}}\Big\{\|\wh x(\omega^K)-x\|>(3+2\gamma')\|x-x_{i_*}\|+2\delta\Big\}\leq \epsilon\;\;\forall(\nu\in J,\,x\in X_\nu).
\ee{finalbasic}
\end{theorem}

\subsection{Aggregation in Euclidean seminorm}\label{sec:eunags}
We now consider the special case of situation described in Section \ref{sect:setting} where $\|\cdot\|$ is an Euclidean seminorm: $\|x\|=\|Bx\|_2$ where $B\in \bR^{q\times n}$ is a given matrix.
In this case, the sets $\X_{ij}(\delta)$ as defined in (\ref{scriptXij})
are finite unions of convex compact sets, and we can apply the ``near-optimal'' inferring color machinery to build the tests
required by aggregation scheme from Section \ref{sec:euag}.

We assume to be given the number of observations $K$ along with tolerance parameters $\epsilon\in(0,1)$ and a ``negligibly small'' $\underline{\delta}>0$ (say, $10^{-100}$); we put $\overline{N}=\half N(N-1)$.
\subsubsection{Preliminaries}\label{euclprelim}
Given a pair $(i,j)\in \cO$ and $\delta>0$, it may happen that one or both of the sets $\X_{ij}(\delta)$ and $\X_{ji}(\delta)$ as defined in (\ref{scriptXij}) is/are empty, in which case we qualify $\delta$ as $(i,j)$-good. Now let $i,j,\delta$ be such that both of the sets $\X_{ij}(\delta)$ and $\X_{ji}(\delta)$ are nonempty. In this case we build the collection $\cR_{ij}(\delta)=\{R_{ij}^s(\delta):1\leq s\leq J_{ij}^{\delta\r}\}$ of nonempty convex compact ``red'' sets  comprised of all nonempty sets of the form
\[
R_{ij\nu}(\delta)=\{\cA_\nu(x):\,x\in X_\nu,\,\|x-x_i\|\leq \|x-x_j\|-\delta\}, \;1\leq \nu\leq J.
\]
Similarly, we build the collection
$\cB_{ij}(\delta)=\{B_{ij}^s(\delta):1\leq s\leq J_{ij}^{\delta\b}\}$ of nonempty convex compact ``blue'' sets  comprised of all nonempty sets of the form
\[
B_{ij\nu}(\delta)=\{\cA_\nu(x):\,x\in X_\nu,\,\|x-x_j\|\leq \|x-x_i\|-\delta\}, \;1\leq \nu\leq J.
\]
Applying to the collections $\cR_{ij}$, $\cB_{ij}$  the $K$-observation color inferring procedure from in Section \ref{sectinfcol}, 
depending on $\delta$
it may happen that the resulting risk bound
$\epsilon_K$ satisfies $\epsilon_K\leq \epsilon/\overline{N}$.
In this case we say that $\delta$ is {\sl $(i,j)$-good}, and that it is $(i,j)$-bad otherwise.
\par
Clearly, whenever $\delta$ is $(i,j)$-good, so is $\delta'\geq\delta$. Similarly to the case of general seminorm, given $(i,j)\in\cO$ and $\delta>0$, we can check efficiently whether $\delta$ is or is not $(i,j)$-good. Given $(i,j)\in\cO$, large enough $\delta$ definitely are $(i,j)$-good, since the corresponding sets $\X_{ij}(\delta)$ are empty. Applying Bisection, we can rapidly find the value $\Delta_{ij}$ of $\delta$ such that $\Delta_{ij}$ is $(i,j)$-good, and either $\Delta_{ij}\leq\underline{\delta}$, or $\Delta_{ij}{-\underline{\delta}}$ is not $(i,j)$-good.
\par
Same as in the case of general seminorm,  it is immediately seen that $\delta$ is $(i,j)$-good if and only if $\delta$ is $(j,i)$-good. As a result, we can select the above $\Delta_{ij}$ to be symmetric: $\Delta_{ij}=\Delta_{ji}$. Note that as a result, every pair $\{i,j\}\in \cU$ is assigned threshold $\Delta_{ij}=\Delta_{ji}$ which is $(i,j)$-good. Besides this, we can equip this pair with $K$-observation test
$\T_{\{i,j\}}$ deciding on the hypotheses $\cH_{ij}(\Delta_{ij}):=H(\cX_{ij}(\Delta_{ij}))$ and $\cH_{ji}(\Delta_{ij}):=H(\cX_{ji}(\Delta_{ij}))$, specifically, the test as follows:
\begin{itemize}
\item[]--- when both hypotheses are empty, the test accepts both hypotheses,
\item[]--- when exactly one of the hypotheses is nonempty, the test accepts this nonempty hypothesis and rejects the empty one,
\item[]--- when both hypotheses are nonempty, $\T_{\{i,j\}}$ is the above color inferring test associated with ($(i,j)$-good!) $\Delta_{ij}$, so that it accepts exactly one of the hypotheses, and its risk does not exceed
$\epsilon/\overline{N}$.
\end{itemize}
\subsubsection{Aggregation routine} Aggregation routine is  the procedure from Section \ref{sec:euag} as applied to the $K$-repeated observation $\omega^K$ in the role of $\xi$ and the just defined $\Delta_{ij}=\Delta_{ij}$, $\T_{\{i,j\}}$; as we have seen,
these entities meet all the requirements of the setup of Section \ref{sec:euag}.
Denoting by $\widehat{i}(\omega^K)\in\{1,...,N\}$ the output of our aggregation, the observation being $\omega^K$, and applying Proposition \ref{prop:l2aggreg10}, we arrive at the following result:
\begin{proposition}\label{prop:consbaseu}
In the situation of this section, suppose that the just described aggregation routine is applied to observation $\omega^K$ stemming from ${x_*}\in X$.
Then
\be
\Prob_{\omega^K\sim P}\left\{
\|{x_*}-x_{\widehat{i}(\omega^K)}\|\leq\|x_{i_*}-{x_*}\|+2\overline{\Delta}\right\}\geq 1-\epsilon\qquad\forall P\in\cP_{x_*}^K
\ee{uppereu}
where $x_{i_*}$ is one of the closest to ${x_*}$ points among $x_1,...,x_N$ and, same as in Proposition \ref{prop:l2aggreg10},
$\overline{\Delta}=\max_{j\neq i}\Delta_{ij}$.
\end{proposition}

\subsubsection{Characterizing performance}
\begin{theorem}\label{thm:lowerbaseu} In the situation under consideration,
assume that for some positive  integer $\overline{K}$, $\epsilon\in(0,1/2)$ and real $\bar{\delta}>{0}$, {\em for every pair $(i,j)\in \cO$} there exists inference {$\omega^{\overline{K}}\mapsto \iota_{ij}(\omega^{\overline{K}})\in\{i,j\}$ such that for every $x_*\in X$ and $P\in\cP^{\overline{K}}_{x_*}$ one has
\begin{equation}\label{haim}
\Prob_{\omega^{\overline{K}}\sim P}\Big\{\|x_*-x_{\iota_{ij}(\omega^{\overline{K}})}\|<\min[\|x_*-x_i\|,\|x_*-x_j\|]+\bar{\delta}\Big\}\geq1-\epsilon.
\end{equation}}

Then whenever
\begin{equation}\label{KLarge}
K\geq \left\lceil {2\ln\left(J\overline{N}/\epsilon\right)\over \ln\left( [4\epsilon(1-\epsilon)]^{-1}\right)} \overline{K}\right\rceil \end{equation}
the aggregated estimate $x_{\widehat{i}(\omega^K)}$ yielded by the above aggregation procedure as applied to $K$-repeated observation $\omega^K$ for every $x_*\in X$ satisfies
\be
 \Prob_{\omega^{K}\sim P}\left\{\|x_{\widehat{i}(\omega^K)}-{x_*}\|\geq\|{x_*}-x_{i_*}\|+2(\bar{\delta}+\underline\delta)\right\}\leq \epsilon\qquad \forall P\in\cP_{x_*}^K,
\ee{finalbasiceu}
$x_{i_*}$ being one of the $\|\cdot\|$-closest  to $x_*$ point among $x_1,...,x_N$.
\end{theorem}
\subsection{Application: adaptive estimation over unions of convex sets}\label{sec:aggunions}
It is clear that just developed aggregation routines may be applied to the problem of adaptive estimation over unions of convex sets defined in Section \ref{sect:setting11}. Our next objective is to discuss this application in more detail and derive corresponding accuracy bounds. From now on, notation and entities such as reliability tolerance $\epsilon$, number $\overline{K}$ of pilot observations, pilot $\overline{K}$-observation estimates $\widetilde{x}_i(\omega^{\overline{K}})$, risks $\Risk^{\cJ}_{\epsilon,M}[\widehat{x}|Y]$, and upper bounds $\myr_j=\myr_j^{\ov K}(\epsilon)$ on $\Risk^{\{j\}}_{\epsilon,\overline{K}}[\widetilde{x}_j|X_j]$, are as defined in that section.

\subsubsection{Estimation over unions using point aggregation}\label{EstPointAggr}
{The quantities} $J=N$,  $X=\cup_jX_j$ and points $x_i=\widetilde{x}_i(\omega^{\overline{K}})$ taken together with the mappings $\cA_j(\cdot)$ and the seminorm $\|\cdot\|$ form the data meeting the requirements of the setup of Section \ref{sect:setting}. Given (post-pilot) $K$-repeated observation $\omega^K$ with $\omega_k\sim p_{\cA_{j_*}}(x_*)$, $k=1,...,K$, with $x_*\in X_{j_*}$,  we can use the routines in Sections \ref{sec:genorma} and \ref{sec:eunags} to
aggregate  points $x_i$ into an estimate $\wh x$ of $x_*$.
\paragraph{Case of general seminorm.}
Let us start with the aggregation procedure described in Section \ref{sec:genorma}. In our present setting its implementation
is as follows. For $(i,j)\in \cO$ we set
\be
r_{ij}=\half\|x_i-x_j\|,\;X_{ij}(\delta)=X\cap \{z:\,\|z-x_i\|\leq r_{ij}-\delta\},
\ee{X_{ij}}
and consider hypotheses $H_{ij}(\delta)$ and $H_{ji}(\delta)$ stating, respectively, that observations stem from a signal
$x\in X_{ij}(\delta)$ and $x\in X_{ji}(\delta)$. Same as before, we say that $\delta$ is $(i,j)$-appropriate, if the risk of the $K$-observation test $\T_{\{i,j\}}$, yielded by the machinery from Section \ref{sectinfcol},  deciding on $H_{ij}(\delta)$ vs $H_{ji}(\delta)$ does not exceed  $\epsilon/(N-1)$. We define parameters $\Delta_{ij}$ and tests $\T_{\{i,j\}}$ as prescribed by the construction in  Section \ref{sec:genorma} and utilize the resulting entities in the aggregation procedure from Section \ref{sec:basic} thus arriving at the aggregated estimate $\wh{x}^{(a)}(\omega^{\ov K+K})=\wh{x}(\omega^{K})$ of $x_*$.
\par
By  Proposition \ref{prop:consbas}, for all $j_*\leq N$ and $x_*\in X_{j_*}$ aggregation $\wh{x}(\omega^{K})$
satisfies
\be
\Prob_{\omega^K\sim p_{\cA_{j_*}(x_*)}^K}\left\{
\|{x_*}-\wh{x}(\omega^{K})\|\leq 3 \|x_{i_*}-x_*\|+2\overline{\Delta}_{i_*}\right\}\geq 1-\epsilon
\ee{eq0001}
where $x_{i_*}$ is one of the $\|\cdot\|$-closest to $x_*$ points among $x_1,...,x_N$, and $\overline{\Delta}_{i_*}\leq\max_{j\neq i_*}\Delta_{i_* j}$ is defined in Proposition \ref{prop:consbas}.
Note that due to ${\omega}^{\ov K}\in\wt\Omega^{\ov K}$ we also have $\|x_*-x_{j_*}\|\leq \myr_{j_*}=\myr_{j_*}^{\ov K}(\epsilon)$, which combines with (\ref{eq0001}) and $\|x_*-x_{i_*}\|\leq \|x_*-x_{j_*}\|$ to imply that the $x_*$-probability for $\omega^K$ to satisfy
\begin{equation}\label{newnewnew}
\|x_*-\wh{x}(\omega^{K})\|\leq  3\myr_{j_*}+2\overline{\Delta}_{i_*}\leq 3\ov\myr+2\overline{\Delta},\;\ov\myr=\max_i\myr_i,\;\overline{\Delta}=\max_{i\leq N}\overline{\Delta}_{i},
\end{equation}
is at least $1-\epsilon$.

\begin{proposition}\label{prop:louna} In the situation described in Section \ref{sect:setting11}, suppose that we are given a positive integer $\overline{K}$, tolerances $\epsilon\in(0,1/2)$ and $\kappa>0$, and $K$ such that
\be
K\geq \left\lceil {2\ln\left(N(N-1)/\epsilon\right)\over \ln\left( [4\epsilon(1-\epsilon)]^{-1}\right)} \overline{K}\right\rceil.
\ee{thesameK}
Then estimate $\wh{x}(\omega^{K})$ yielded by  the procedure described above with properly selected $\Delta_{ij}$ as applied to observation $\omega^K$ satisfies
\[
\Prob_{\omega^{K}\sim p^{K}_{\cA_{j_*}(x_*)}}\Big\{\|\wh{x}(\omega^{K})-x_*\|>3\max_i\myr^{\ov K}_i(\epsilon)+2\RiskOpt^{\overline{1,N}}_{\epsilon,\ov K}[{X}]+\kappa\Big\}\leq \epsilon.
\]
In particular, when the upper bounds $\myr^{\ov K}_i(\epsilon)$ on the risks $\Risk^{\{i\}}_{\epsilon,\overline{K}}[\widetilde{x}_i|{X}_i]$ of estimates $\wt{x}_i(\omega^{\ov K}) $ are within factor $\theta$ of the respective $\overline{K}$-observation minimax risks, i.e.,
\[
{\RiskOpt^{\{i\}}_{\epsilon,\ov K}}[{X}_i]\leq\myr^{\ov K}_i(\epsilon)\leq  \theta\RiskOpt^{\{i\}}_{\epsilon,\ov K}[{X}_i]
\]
the risk ${\Risk^{\overline{1,N}}_{2\epsilon,\ov K + K}}[\widehat{x}^{(a)}|{X}]$ of the estimate $\widehat{x}^{(a)}(\omega^{\ov K+K})=\wh{x}(\omega^{K})$ (as function of pilot
observation $\omega^{\ov K}$ and independent observation $\omega^K$) is within a moderate factor from  the minimax $\overline K$-observation risk ${\RiskOpt^{\overline{1,K}}_{\epsilon,\ov K}}[{X}]$:
\[
\Risk^{\overline{1,N}}_{2\epsilon,\ov K + K}[\widehat{x}^{(a)}|{X}]\leq [2+3\theta]\RiskOpt^{\overline{1,N}}_{\epsilon,\ov K}[{X}]+\kappa.
\]
\end{proposition}

\paragraph{Case of Euclidean seminorm.}
When $\|\cdot\|$ is a Euclidean seminorm,  we can utilize the aggregation procedure described in Section \ref{sec:eunags} to build the ``two-stage'' estimate $\wh{x}^{(a)}(\omega^{\ov K+K})=\wh{x}(\omega^K)$. Specifically, given $\epsilon\in(0,\half)$, ``negligibly small'' $\underline{\delta}>0$,  and $\delta\geq 0$, consider
 sets
\[
\X_{ij}(\delta)=\{z\in X: \|z-x_j\|\geq \delta+\|z-x_i\|\},\;\X_{ji}(\delta)=\{z\in X: \|z-x_i\|\geq \delta+\|z-x_j\|\}.
\]
We apply the construction of
Section \ref{sec:eunags} to compute for every $(i,j)\in\cO$ $(i,j)$-good quantities $\Delta_{ij}=\Delta_{ji}$ such that either $\Delta_{ij}\leq\underline{\delta}$, or $\Delta_{ij}>\underline{\delta}$ and $\Delta_{ij}{-\underline{\delta}}$ is not $(i,j)$-good, and proceed as explained in that section, ending up with the aggregated estimate $\wh{x}(\omega^K)$. Invoking Proposition
\ref{prop:consbaseu},  we have
\[
\Prob_{\omega^K\sim p_{\cA_{j_*}(x_*)}^K}\left\{
\|x_*-\wh{x}(\omega^K)\|\leq\|x_{i_*}-x_*\|+2\overline{\Delta}\right\}\geq 1-\epsilon, \;\overline{\Delta}=\max_{i,j}\Delta_{ij},
\]
where $x_{i_*}$ is one of the closest to $x_*$ points among $x_1,x_2,...,x_N$.
We have the following analog of Proposition \ref{prop:louna} in this case.
\begin{proposition}\label{prop:lounae} Let $\|\cdot\|$ be a Euclidean seminorm. In the situation described in Section \ref{sect:setting11},  suppose that we are given a positive integer $\overline{K}$, tolerances $\epsilon\in(0,1/2)$ and $\varkappa>0$, and $K$ satisfying
\be
K\geq \left\lceil {2\ln\left(N^2(N-1)/(2\epsilon)\right)\over \ln\left( [4\epsilon(1-\epsilon)]^{-1}\right)} \overline{K}\right\rceil.
\ee{thesameKe}
Then estimate $\wh{x}(\omega^{K})$ yielded by  the above procedure with properly selected parameters as applied to observation $\omega^K$ satisfies
\[
\Prob_{\omega^{K}\sim p^{K}_{\cA_{j_*}(x_*)}}\Big\{\|\wh x(\omega^K)-x_*\|>\max_i\myr^{\ov K}_i(\epsilon)+4\RiskOpt^{\overline{1,N}}_{\epsilon,\ov K}[{X}]+{\varkappa}\Big\}\leq \epsilon.
\]
In particular, when the upper bounds $\myr^{\ov K}_i(\epsilon)$ on the partial risks $\Risk_{\epsilon,\ov K}^{\{i\}}[\wt x_i|{X}_i]$ of ${\ov K}$-observation estimates $\wt x_i(\cdot) $ are within a factor $\theta$ of the respective $\overline K$-observation minimax risks, i.e.,
\[
\RiskOpt_{\epsilon,\ov K}^{\{i\}}[X_i]\leq\myr_i^{\ov K}(\epsilon)\leq  \theta\RiskOpt_{\epsilon,\ov K}^{\{i\}}[X_i],
\]
the maximal risk $\Risk_{2\epsilon,\overline{K}+K}^{\overline{1,N}}[\widehat{x}^{(a)}|X]$ of aggregated estimate $\wh x^{(a)}(\omega^{{\ov K}+K}):=\wh x(\omega^{K})$ (considered as function of the pilot observation $\omega^{\ov K}$ and independent observation $\omega^K$) is within a moderate factor from  the minimax $\overline K$-observation risk $\RiskOpt_{\epsilon,\ov K}^{\overline{1,N}}[{X}]$:
\[
\Risk_{2\epsilon,K+\ov K}^{\overline{1,N}}[\widehat{x}|{X}]\leq [4+\theta]\RiskOpt_{\epsilon,\ov K}^{\overline{1,N}}[{X}]+\varkappa.
\]
\end{proposition}

\section{Adaptive estimation over unions of ellitopes}\label{sec:ellall}
\subsection{Ellitopic setup}\label{sec:elli}
Ellitopes, as introduced in \cite{l2estimation,juditsky2020statistical}, are symmetric w.r.t. the origin convex and compact sets.  In this section we consider the special case of the estimation problem described in Section \ref{sect:setting11} in which
\begin{enumerate}
\item observation scheme is Gaussian, i.e., observations $\omega_k\in \bR^m$ stemming from $(j,x)$, $x\in X_j$, are normal,
$\omega_k\sim \cN(\cA_j(x),I_m)$ where $\cA_j(x)$ are {\sl linear}, rather than affine, mappings: $\cA_j(x)=A_jx$, where $A_j\in\bR^{m\times n}$, $j=1,...,N$, are given matrices;
\item sets $X_j$ $j=1,...,N$, are {\em basic ellitopes}:
\[
X_j=\{x\in\bR^n:\, \exists r\in\cR_j: x^TR_{j\tau} x\leq r_\ell,\,\tau\leq L\}, \quad j=1,...,N,
\]
\item seminorm $\|\cdot\|$ is of the form $\|x\|=\pi(Bx)$ where $B$ is a $q\times n$ matrix and the unit ball $\cB_*$ of the conjugate to $\pi(\cdot)$ norm $\pi_*(\cdot)$ is an ellitope 
\[
\cB_*=\{y\in\bR^q: y\in\bR^q:\,\exists s\in\cS: x^TS_\tau x\leq r_\tau,\,\tau\leq L'\}.
\]
Here
\begin{itemize}
\item $\cR_j\subset\bR^{L}_+$, $j=1,...,N$, and $\cS\subset\bR^{L'}_+$ are computationally tractable convex compact sets intersecting with $\inter \bR^{L}_+$ which are {\em monotone}.\footnote{Here monotonicity of $V\subset \bR^k_+$ means that if  $0\leq v'\leq v$ and $v\in V$ then also $v'\in V$.}
\item $R_{j{\tau}}$, $1\leq j\leq N$, $1\leq {\tau}\leq L$, are $n\times n$ matrices with $R_{j{\tau}}\succeq 0$ and $\sum_{\tau} R_{j{\tau}}\succ 0$; $S_{\tau}$ are
     $q\times q$ matrices such that  $S_{\tau}\succeq 0$ and $\sum_{\tau} S_{\tau}\succ 0$.
\end{itemize}
We refer to $L$ and $L'$ as {\em sizes} of corresponding ellitopes.
\end{enumerate}
Particular choices of sets $X_j$ and seminorm $\|\cdot\|$ encompass a variety of situations.
\begin{itemize}
	\item When $L = 1 $, $\cR_j = [0,1]$ and $ R^{(j)}_1 \succ 0 $, $ X_j $ is an ellipsoid.
	\item When $L \geq 1 $, $\cR_j = [0,1]^L$, $X_j$ is an  intersection of ellipsoids and elliptic cylinders centered at the origin, $\bigcap_{{\tau}\leq L} \{ z:  z^TR_{j{\tau}} z\leq  1 \}.$
	\item When $U=[u_1,...,u_L] \in \bR^{n\times L}$, $\rank[U] = n$, $\cR_j = [0,1]^L$,  and $R_{j{\tau}}=u_{\tau} u_{\tau}^{T}$,  $X_j$ is a symmetric w.r.t. the origin polytope $\{z: \|U^T z\|_{\infty} \leq 1\}$.	
\item
When for $ p \geq 2$, $ \cS = \Big\{s \in \bR^{L'}_+ : \sum_{\tau} [s]_{\tau}^{p/2} \leq 1  \Big\} $ and as, in the previous example $U=[u_1,...,u_{L'}]\in \bR^{q\times L'}$, $\rank[U] = q$, and
$S_{\tau} = u_{\tau} u_{\tau}^{T}$, $\cB_*$ is the set $ \{ y: \| U^T y \|_{p} \leq 1\}$ and the seminorm $\|\cdot\|$ is $\|w\|=\|UBw\|_{p/(p-1)}$.
\end{itemize}
The family of ellitopes admits simple and fully algorithmic ``calculus'' demonstrating that this family is closed w.r.t. nearly all operations preserving convexity and symmetry w.r.t. the origin (e.g., taking finite intersections, direct products, linear images, and inverse images under linear embeddings; for details, see \cite[Section 4.6]{juditsky2020statistical}).

We are about to show that in the present situation, {\sl estimates yielded by the approach described in Section \ref{sec:genest} are nearly optimal in the minimax sense}.\footnote{An analog of the results below in the special case where $\|\cdot\|$ is a Euclidean seminorm can be obtained by applying construction of Section \ref{sec:l2est}.} Moreover, in this case we are able to provide ``reasonably good'' bounding of minimax risks of recovery over pairwise unions $X_i\cup X_j$ of ellitopes implying that tight bounds for the minimax risk of estimation over $X=\cup_{i=1}^NX_i$ can be efficiently computed.
\subsection{Near-optimality of the aggregated estimate}
Let $\ov K$ and $K$ be positive integers, and let us assume that in the situation described in Section \ref{sec:elli} we are given $\epsilon\in(0,1/8)$ and $M=\ov K+K\geq 2$ independent observations $\omega_k$ stemming from unknown pair $(\ell_*,x_*),\,x_*\in X_{\ell_*}$, $1\leq \ell_*\leq N$. To build an $M$-observation estimate  $\wh x^{(a)}(\omega^M)$ of $x_*$ we proceed as explained in Section \ref{sec:genest}:
\begin{itemize}\item we split the observation sample into two observations: a $\ov K$-repeated observation  $(\omega_1,...,\omega_{\ov K})$ (preliminary observation) and  $\omega^K=(\omega_{\ov K+1},...,\omega_{\ov K+K})$ (secondary observation).
\item Preliminary observation is averaged to build observation
\[
\wh\omega={1\over \ov K}\sum_{k=1}^{\ov K} \omega_k
\]
which is then used to compute $N$ {\em polyhedral} estimates $\wt x_i({\widehat{\omega}})$ following the recipe in \cite{juditsky2020polyhedral} and \cite[Section 5.1.5]{juditsky2020statistical}.
\item Finally, we apply the aggregation routine from Section \ref{sec:genest} to assemble points $x_i= \wt x_i({\widehat{\omega}})$ into estimate $\wh x(\omega^K)$ obtaining as a result adaptive estimate
 $\wh x^{(a)}(\omega^{M})=\wh x(\omega^{K})$ of $x_*$.
\end{itemize}
Recall  (cf. \cite[Proposition 5.10]{juditsky2020statistical}) that polyhedral estimates $\wt x_i$ satisfy $\wt x_i(\wh \omega)\in X_i$ and
\be
{\Risk^{\{i\}}_{\epsilon,1}}[\wt x_i|X_i]\leq \myr_{i}(\epsilon)\leq \myC_1 \ln (L+L')\sqrt{\ln\left[{m/\epsilon}\right]}\,
{\overline{\RiskOpt}^{\{i\}}_{{1\over 8},1}}[X_i]
\ee{barrisk}
where bound $\myr_i(\epsilon)$ for {maximal risk of estimation under the $i$th observation model} is  efficiently computable and ${\overline{\RiskOpt}^{\{i\}}_{\epsilon,1}}[X_i]$ is the minimax $\epsilon$-risk of recovering $x\in X_{i}$ from single ``averaged'' observation $\wh \omega$ stemming from i.i.d. $\omega_k\sim\cN(A_jx,I_m)$; From now on, $\myC_i$ stand for appropriate {\sl absolute} constants.

Given $(i,j)\in\cO$,  positive integer $K$, and ${\delta}\in (0,1/2)$ we (re-)define the notion of $\delta$-separation risk (cf. \rf{Deltaij}) in the present situation according to
\be
\myrg^{K}_{ij}({\delta})=\half \max_{x\in X_i,y\in X_j} \left\{\|x-y\|:\,\|\cA_i(x)-\cA_j(y)\|_2\leq{2\over \sqrt{K}} q_{\N}(1-{\delta})
\right\}
\ee{DeltaijN}
where $q_{\N}(p)$ is the
$p$-quantile of the standard normal distribution: ${1\over\sqrt{2\pi}}\int_{-\infty}^{q_{\N}(p)}{\rm e}^{-s^2/2}ds=p$. Note that (\ref{DeltaijN}) is feasible, and therefore solvable, due to $0\in X_j$ and $\cA_j(0)=0$ for all $j$.
\par
For the sake of simplicity, from now on we restrict  ourselves to the case
\begin{equation}\label{restrict}
K\leq\overline{K}
\end{equation}
The next statement provides a refined version of results of Section \ref{sec3.2.2}  in the present setting:

\begin{theorem}\label{thm:genthmG}
In the situation of this section, assuming (\ref{restrict}) and $0<\epsilon<1/16$, the just built estimate $\wh x^{(a)}$ (as function of pilot
observation $\omega^{\ov K}$ and secondary observation $\omega^K$) satisfies
\begin{equation}\label{eq0999}
\Risk^{\{i\}}_{2\epsilon,\ov{K}+K}[\wh x^{(a)}|X_i]\leq 2\myr_{i}(\epsilon)+\max_{j< i}\left[\myr_j(\epsilon)+2\myrg_{ij}^{K}\big(\varepsilon\big)\right]\,\,\forall i\leq N,
\end{equation}
with $\varepsilon={\epsilon\over N-1}$.
Moreover, setting
\[
\ov \vartheta:={\sqrt{K}q_{\N}(1-\epsilon)\over\sqrt{\ov K} q_{\N}(1-\varepsilon)}
\]
one has $\ov \vartheta\leq1$ and
\be
{\Risk^{\{i\}}_{2\epsilon, \ov{K}+K}}[\wh x^{(a)}|X_i]
\leq 2\myr_{i}(\epsilon)+\max_{j< i}\left[\myr_j(\epsilon)+
2\ov\vartheta^{-1}\RiskOpt^{\{i,j\}}_{\epsilon,\overline K}[X_i\cup X_j]\right]\quad\forall i\leq N,
\ee{barthe21}
whence, in particular,
\[
{\Risk^{\{i\}}_{2\epsilon, 2\ov{K}}}[\wh x^{(a)}|X_i]
\leq \max_{1\leq j\leq  i}\left[3\myr_{j}(\epsilon)
+2{q_{\N}\left(1-{\varepsilon}\right)\over q_{\N}(1-\epsilon)}\RiskOpt^{\{i,j\}}_{\epsilon,\overline K}[X_i\cup X_j]\right]
\qquad\forall i\leq N.\]
\item Besides this, one has
\bse
{\Risk^{\{i\}}_{2\epsilon, 2\ov{K}}[\wh x^{(a)}|X_i]}
\leq \myC_2\Big(\ln{(L+L')}\sqrt{\ln\left[{m/\epsilon}\right]}+{\sqrt{\ln[N/\epsilon]}}\Big)\max_{1\leq j\leq i}\RiskOpt^{\{i,j\}}_{{1\over 16},\overline K}[X_i\cup X_j] \qquad\forall i\leq N,
\ese
\item so that
\be
{\Risk^{\overline{1,N}}_{2\epsilon, 2\ov{K}}}[\wh x^{(a)}|X]\leq
\myC_3\Big(\ln {(L+L')}\sqrt{\ln\left[{m/\epsilon}\right]}+{\sqrt{\ln[N/\epsilon]}}\Big){\RiskOpt^{\overline{1,N}}_{{1\over 16},\overline K}}[X].
\ee{maxung}
\end{theorem}\noindent

\subsection{Bounding the maximal risk of estimation}
Our current objective is to provide efficient bounding for separation risks $\myrg^K_{ij}(\varepsilon)$; taken together with bounds $\myr_j(\epsilon)$ for partial risks this would allow to bound the minimax risk of estimation over $X$.
Under the premise of Theorem \ref{thm:genthmG}, let ${(i,j)}$, $1\leq i, j\leq N$, and $K\leq\ov{K}$ be fixed, let, same as in (\ref{eq0999}), $\varepsilon=\epsilon/(N-1)$,  and let $\delta=2K^{-1/2}q_{\N}(1-\varepsilon)$. Observe that
\be
\myrg^{K}_{ij}(\varepsilon)&=&\half \max_{x\in X_i,y\in X_j} \left\{\|x-y\|:\,\|A_ix-A_jy\|_2\leq \delta
\right\}\nn
&=&\half \max_{x\in X_i,y\in X_j} \left\{\pi(B(x-y)):\,\|A_ix-A_jy\|^2_2\leq \delta^2 
\right\}\nn
&=&\half \max_{[x;y;u]} \left\{u^TB(x-y):\,u\in \cB_*,\,x\in X_i,\,y\in X_j,\,\|A_ix-A_jy\|_2^2
\leq \delta^2
\right\}.
\ee{DeltaijNN}
Because the direct product of ellitopes $\cB_*\times X_i\times X_j$ is an ellitope of the size not exceeding $L'+2L$ (cf. \cite[Section 4.6]{juditsky2020statistical}), when writing
$u^TB(x-y)= [u;x;y]^T\overline B[u;x;y]$ and $\|A_ix-A_jy\|_2^2=[u;x;y]^TQ_{ij}[u;x;y]$
with \[
\overline B=\left[\begin{array}{r|c|c} 0_{q\times q}&\half B&-\half B\\\hline
\half B^T& 0_{n\times n}&0_{n\times n}\\\hline
-\half B^T&0_{n\times n}&0_{n\times n}
\end{array}
\right],\;\;Q_{ij}=\left[\begin{array}{r|r|r} 0_{q\times q}&0_{q\times n}&0_{q\times n}\\\hline
0_{n\times q}& A_i^TA_i&-A_i^TA_j\\\hline
0_{n\times q}&-A_j^TA_i&A_j^TA_j
\end{array}
\right]
\] we conclude that the quantity $\myrg^{K}_{ij}(\varepsilon)$ is the maximum of a homogeneous quadratic form over an ellitope of
size at most $\overline D=L'+2L+1$. Therefore, it 
can be upper-bounded by an efficiently computable quantity $\overline \myrg^{K}_{ij}(\varepsilon)$ within factor $2\ln \overline D+2\sqrt{\ln \overline D}+1$ (see, e.g., \cite[Proposition 3.3]{l2estimation}) using semidefinite relaxation.

As a result, given a pair ${(i,j),i\neq j}$, we can upper-bound the $2\epsilon$-minimax risk
$\RiskOpt^{\{i,j\}}_{2\epsilon,2\ov K} [X_i\cup X_j]$
with efficiently computable quantity
\[
\overline \myr_{ij}(\epsilon)=3\max[\myr_i(\epsilon),\myr_j(\epsilon)]+2\overline\myrg^{\ov K}_{ij}(\varepsilon)
\]
 such that
\bse
\overline \myr_{ij}(\epsilon)&\leq& \myC_6\ln[\overline D]\sqrt{\ln[m/\epsilon]}
\max\big({\RiskOpt^{\{i\}}_{{1\over 16},\ov K}}[X_i],{\RiskOpt^{\{j\}}_{{1\over 16},\ov K}}[X_j]\big)
+\myC_7\ln [\overline D]\myrg^{\ov K}_{ij}(\varepsilon)\\
&\leq&
\myC_8\ln[\overline D][\sqrt{\ln[m/\epsilon]}+\sqrt{\ln[N/\epsilon]}]\RiskOpt^{\{i,j\}}_{{1\over 16},\ov K}[X_i\cup X_j]
\ese
(we have used \rf{eq12222}).
Similarly, $2\epsilon$-minimax risk ${\RiskOpt^{\overline{1,N}}_{2\epsilon,2\ov K}} [X]$ of estimation over $X$ can be bounded with efficiently computable quantity
\[
\overline \myr(\epsilon)=\max_{i,j\leq N}\left[3\myr_i(\epsilon)+2\overline\myrg^{\ov K}_{ij}\big(\varepsilon\big)\right]
\]
such that
\bse
\overline \myr(\epsilon)&\leq&
\myC_9
\max_{i,j\leq N}\left[\ln[\overline D]\sqrt{\ln[m/\epsilon]}{\RiskOpt^{\{i\}}_{{1\over 16},\ov K}}[X_i]+\myC_{10}\ln [\overline D]\myrg^{\ov K}_{ij}\big(\varepsilon\big)\right]\\&\leq&
\myC_{11}\ln[\overline D]\left(\sqrt{\ln[m/\epsilon]}+\sqrt{\ln[N/\epsilon]}\right){\RiskOpt^{\overline{1,N}}_{{1\over 16},\ov K}}[X].
\ese
\subsection{Numerical illustration: application to estimation in the single-index model}
\begin{figure}
\vspace{-0.5cm}
\centering
\includegraphics[scale=0.8]{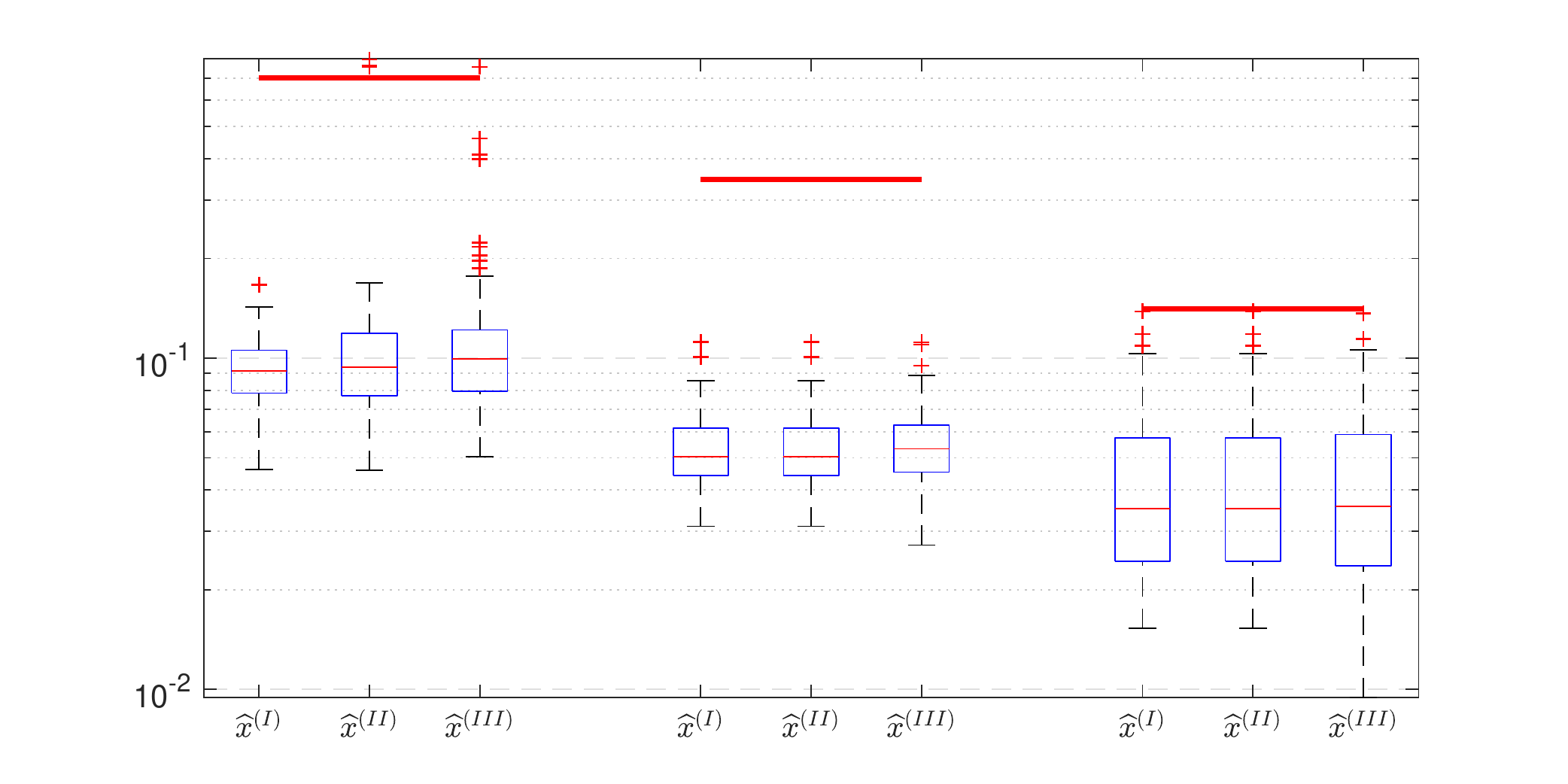}
\caption{\label{fig:1} Error distribution of recoveries $\wh x^{(I)}$, $\wh x^{(II)}$,  and $\wh x^{(III)}$ for different values of noise variance $\sigma^2$: from left to right, box plots for $\sigma=0.1$, $ \sigma=0.05$, and $\sigma=0.02$.}
\end{figure}

In this section, we apply the proposed adaptive estimate to a toy problem of estimation  in the simple single index model in which
\begin{itemize}
\item ``Unknown signal'' $x$ is a vector of coefficients of one-dimensional spline $s(t)$ on $[-1,1]$ split into 10 equal segments. In each segment, $s$ is a quadratic polynomial, and its derivative $s'(t)$ is continuous on the entire $[-1,1]$, making the number of degrees of freedom in the spline---dimension of the parameter vector $x$---equal to 12. Signal vector $x$ is restricted to have $\|\cdot\|_2$-norm not exceeding 1, thus, the signal set $X$ is the unit Euclidean ball in $\bR^{12}$.
\item We consider the situation in which all signal sets $X_j$, $j=1,...,N=64$, are equal to $X$, but there are $N$ different encodings $\cA_j(\cdot)=A_j\in \bR^{1024\times 12}$ built as follows: for $j=1,...,J=64$, we specify $e_j$ as unit vector in $\bR^2$ at angle $2\pi(j-1)/N$ with the first basis vector. Specifying $\Gamma$ as a set of $1024$ points sampled from a uniform distribution on
    $\{u\in\bR^2:\|u\|_\infty\leq 1\}$,  $A_jx$ is the restriction onto $\Gamma$ of the function $f_{j,x}(u)=s(e_j^Tu)$. \\
    Note: for $u\in\Gamma$, $e_j^Tu$ can be outside of $[-1,1]$, and when defining $s(e_j^Tu)$,  we extend $s$ from $[-1,1]$ onto the entire real axis in such a way that the extended function is continuously differentiable and is affine to the left of $-1$ and to the right of $1$.
\item Observations $A_jx$ are corrupted by white Gaussian noise $\xi\sim \N(0,\sigma^2 I)$.
\item We deal with $\overline{K}=K=1$ and split our actual observation into two independent unbiased Gaussian observations---pilot $\wt{\omega}$ and secondary $\omega$---of variance $2\sigma^2$ each.
\end{itemize}
It is worth mentioning that the considered situation differs from the ``classical'' setting of the single index estimation problem: here our objective is neither to estimate the index---unit vector $e$ corresponding to the ``orientation'' in $\bR^2$ of the univariate function underlying observations, nor to estimate the bivariate regression function $f_{i,x}(\cdot)$,\footnote{For ``state of art'' adaptive estimates of regression function $f$ in a general $d$-dimensional single index model under $L_2([-1,1]^d)$-losses  see, e.g.,  \cite{golubev1992asymptotic}; see also \cite{lepski2014adaptive} for adaptation w.r.t pointwise and general $L_p$-risks.}   but to recover vector $x$ of spline coefficients of $s(\cdot)$, the norm $\|\cdot\|$ being the Euclidean norm. As such, the problem we consider is that of recovery from noisy indirect observations, the latter being equivalent to estimating univariate function $s(\cdot)$, estimation error being measured in the $L_2$-norm on $[-1,1]$.
We consider two implementations of the recovery procedure; in both implementations we utilize polyhedral estimate of \cite{juditsky2020polyhedral} to build pilot estimates $\wt x_i(\wt\omega)$, $i=1,...,N$. The first recovery, we denote it $\wh x^{(I)}$, utilizes the aggregated estimate described in Sections \ref{sec:eunags}, \ref{sec:aggunions};    $\wh x^{(II)}$ is the adaptive estimate of Section \ref{sec:l2est}; finally, estimate $\wh x^{(III)}$ is the slightly modified adaptive estimate of Section \ref{sec:genest} in which, when the set $\cI(\omega)$ of admissible estimates contains more than 1 point, instead of selecting the admissible estimate with the smallest index $i$, adaptive estimate $\wh x$ is obtained by aggregating admissible points $\wt x_i,\,i\in \cI(\omega)$, as the optimal solution to the optimization problem
\[
\wh x=\argmin_u\max_{i\in\cI(\omega)}\|u-\widetilde{x}_i\|_2.
\]
To see how the error of recovery depends on the noise variance $\sigma^2$, for each value of the variance we sample $K=100$ realizations of the signal $x_k$ from the uniform distribution on the unit sphere along with directions $e_k$ from the uniform distribution on the unit circle. Results of these experiments are presented in Figure \ref{fig:1} (note the logarithmic scale of the $y$-axis); the red bar over each box plot represents the upper bound $\max_j\myr^{1}_{j}(\epsilon)$ of partial $\epsilon$-risks of preliminary estimates  $\wt x_j(\wt\omega)$.

\begin{figure}
\vspace{-0.5cm}
\begin{tabular}{cc}
\includegraphics[scale=0.55]{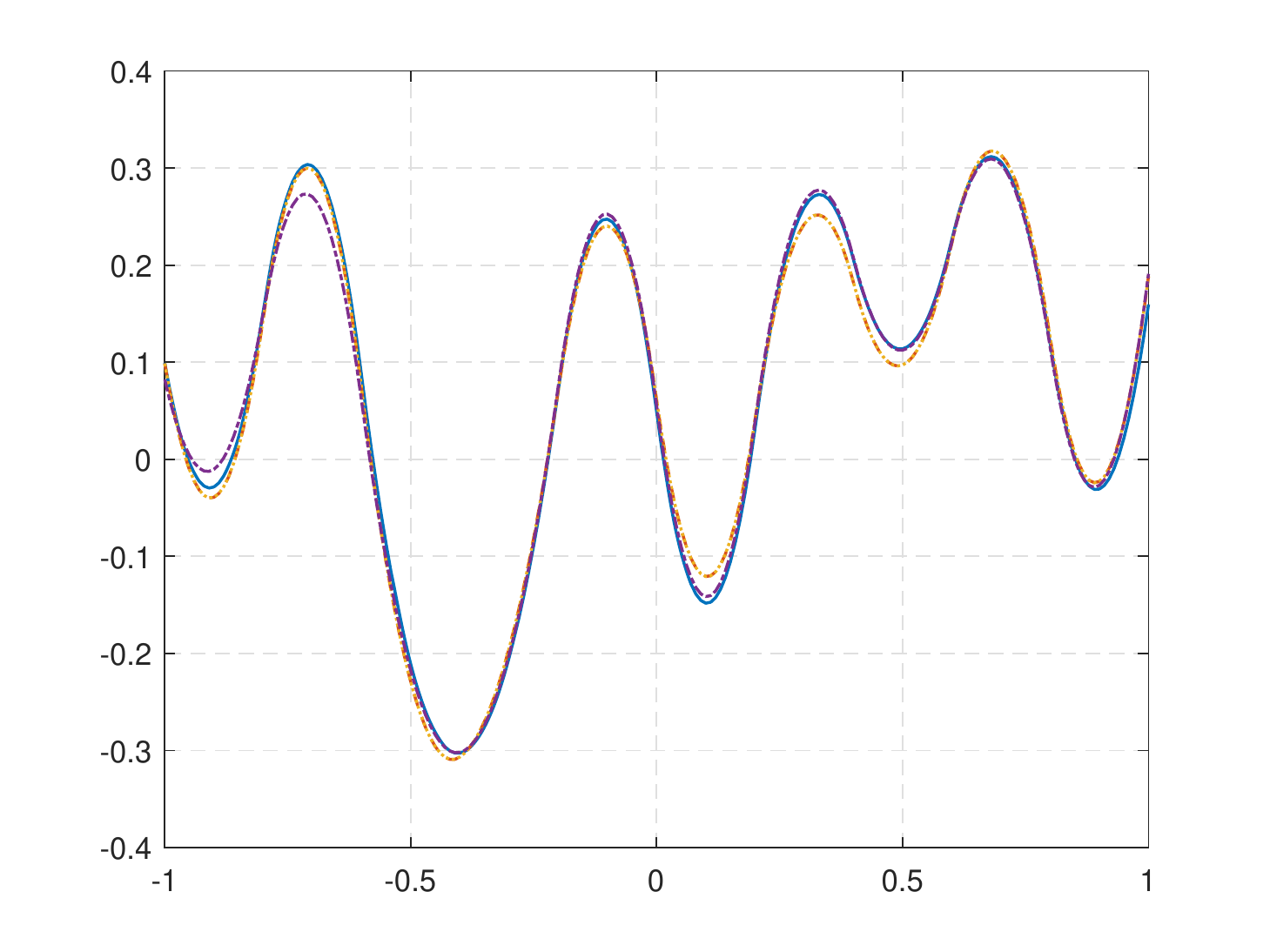}&\includegraphics[scale=0.55]{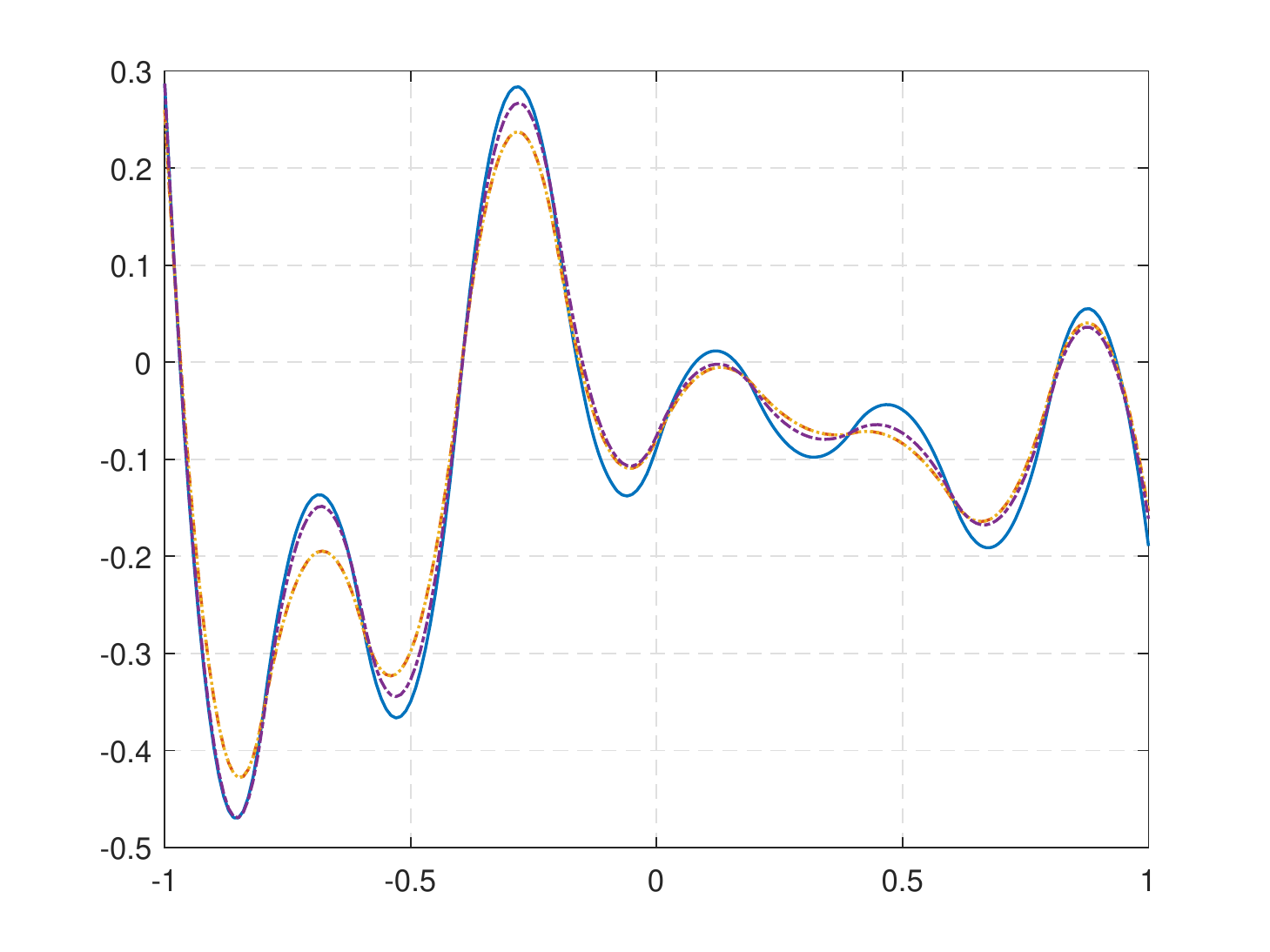}\\(a)&(b)
\end{tabular}
\caption{\label{fig:0} Typical graphs of the true function $s(\cdot)$ (solid line) and its recoveries utilizing estimate $\wh x^{(I)}$ (dotted line), estimate $\wh x^{(II)}$ (dash-dot line), and estimate $\wh x^{(III)}$ (dashed line).}
\end{figure}

\begin{figure}
\vspace{-0.5cm}
\centering
\includegraphics[scale=0.8]{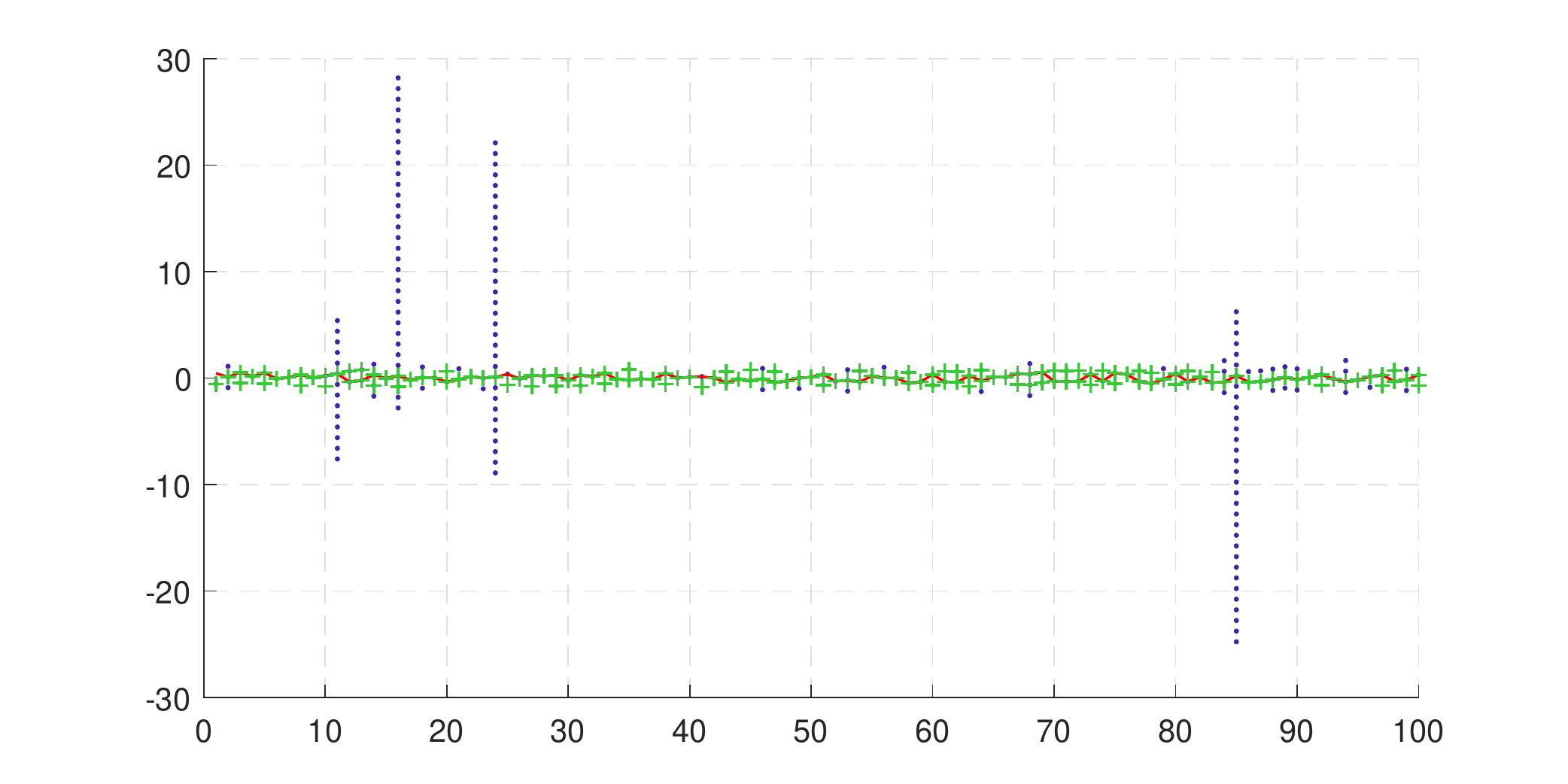}
\caption{\label{fig:2} Admissible sets $\cI(\omega)$ of  recovery $\wh x^{(III)}$ (blue dots), those of recovery $\wh x^{(II)}$ (green crosses),  and selected $\wh i$'s of recovery $\wh x^{(I)}$ (solid red line).}
\end{figure}

\begin{figure}
\vspace{-0.5cm}
\centering
\includegraphics[scale=0.8]{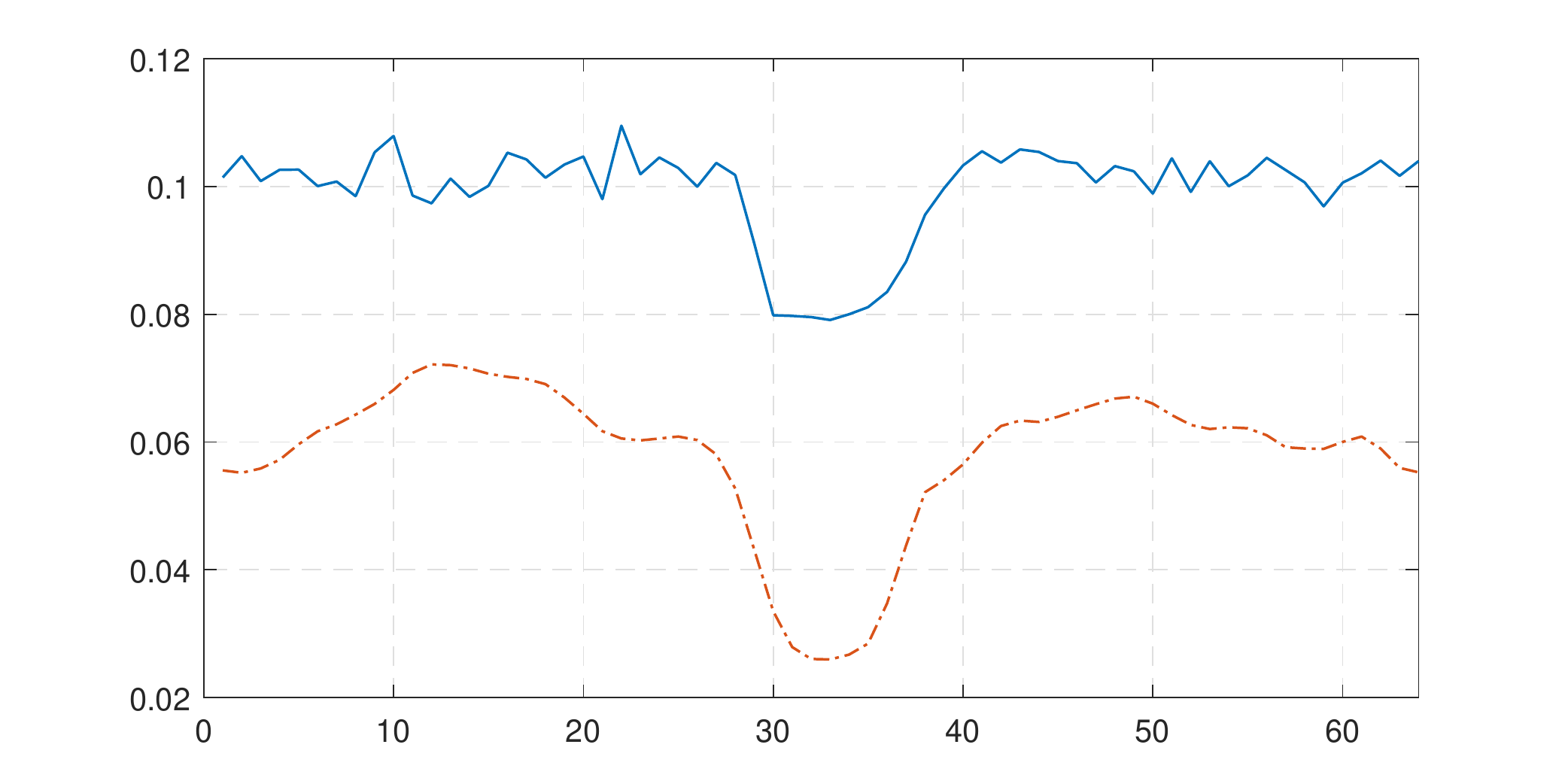}
\caption{\label{fig:3} Bounding $\myrg^{1}_{ij}(\varepsilon)$: solid line---upper bound $\ov\myrg^{1}_{ij}(\varepsilon)$ by semidefinite relaxation,
dash-dot line---lower bound by Monte-Carlo simulations.}
\end{figure}

The reliability parameter of the recoveries being set to 95\% (i.e., $\epsilon=0.05$), upper bounds $\myr^1_i(\eps)$ exceed 1 for $\sigma>0.15$.
We present in Figure \ref{fig:0}, for $\sigma=0.1$,  typical graphs of the true signal $s(\cdot)$ and its recoveries.
\begin{itemize}
\item Plot (a):  set $\cI(\omega)$ of admissible estimates for recoveries $\wh x^{(II)}$ and $\wh x^{(III)}$ is a singleton (in this case, $\|x_*-\wh x^{(I)}\|_2=0.0949$, $\|x_*-\wh x^{(II)}\|_2=0.0616$,
and $\|x_*-\wh x^{(III)}\|_2=0.0710$).
\item Plot (b): cardinality $|\cI(\omega)|=3$ of the set of  admissible estimates for recovery $\wh x^{(III)}$,  $\cI(\omega)$  is a singleton for recovery $\wh x^{(II)}$ (in this case, $\|x_*-\wh x^{(I)}\|_2=0.0752$, $\|x_*-\wh x^{(II)}\|_2=0.0846$,
and $\|x_*-\wh x^{(III)}\|_2=0.1508$).
\end{itemize}
For $\sigma\leq 0.05$ in all simulations the set $\cI(\omega)$ of admissible estimates was a singleton for all recoveries. Moreover, in these simulations, selected indices  $j$ of encodings $A_j$ were the same for both recoveries, corresponding to the closest to the ``true direction $e$''
element $e_j$ of the ``grid of directions.'' When $\sigma=0.1$, corresponding admissible sets for recovery $\wh x^{(I)}$ were singletons, with corresponding direction  being the closest to true $e$ in $93/100$ simulations (and second close in remaining $7/100$); admissible set for recovery
$\wh x^{(II)}$ was a singleton corresponding to the closest direction in $56/100$ experiments, in the remaining $44/100$ the admissible set contained two closest to $e$ directions. ``Population'' of admissible sets of recovery $\wh x^{(III)}$ is represented in Figure \ref{fig:2}; admissible $i$'s
obtained in each simulation are ``centered'' w.r.t. the ``index'' $j_e={N\theta \over 2\pi}+1$ where $\theta$ is the angle between random vector $e$
underlying the observation and the first basis vector of $\bR^2$.
For $\sigma=0.1$ we also present in Figure \ref{fig:3} typical plot of the bound $\ov\myrg^{1}_{ij}(\varepsilon)$, $\varepsilon={\epsilon/(N-1)}$,
$i=33$, for separation risk $\myrg^{1}_{ij}(\varepsilon)$ along with the lower bound computed by Monte-Carlo simulations.

\appendix
\section{Proofs}
\subsection{Proofs for Section \ref{sec:adest}}
\subsubsection{Proof of Proposition \ref{ai1}}
The fact that the ${(\ell_*,x_*)}$-probability of $\overline\Omega^K$ is at least $1-\epsilon$  is readily given by the union bound and the fact that for good pairs ${(\ell_*,j)}$ the ${(\ell_*,x_*)}$-probability for test $\T_{\{\ell_*,j\}}$ to accept $H_{\ell_*}$ is at least $1-\epsilon/(N-1)$. Furthermore, because the preliminary observation belongs to $\widetilde{\Omega}^{\ov K}$ we have
\begin{equation}\label{neweq999}
 \|x_*-x_{\ell_*}\|\leq \myr_{\ell_*}.
 \end{equation}
Let now $\omega^K\in \overline \Omega^K$ be fixed; then the set $\cI=\cI(\omega^K)$ is not empty because  $\ell_*\in\cI$; indeed, when $\omega^K\in \overline \Omega^K$ ``true'' hypothesis $H_{\ell_*}$ is never rejected. Consequently, if $i'\in\cI$ differs from $\ell_*$, then $(i',\ell_*)$ is bad, since otherwise $i'\in \cJ_{\ell_*}$ and the test $\T_{\{i',\ell_*\}}$ would reject the true hypothesis $H_{\ell_*}$ (otherwise $i'\in\cI$ would be  impossible), which contradicts  $\omega^K\in \overline{\Omega}^K$. As a byproduct of the just made observation, $\cI_{\ell_*}^-\subset J_{\ell_*}^-$. Now, since we are in the case of $\ell_*\in\cI$, either $\widehat{i}=\ell_*$, or $\cI\ni \widehat{i}<\ell_*$. In the first case, (\ref{1stb}) is evident, in the second $\widehat{i}\in\cI_{\ell_*}^-$, and therefore \rf{1stb} holds true as well. \rf{1stb} is proved.
Next, the first inequality in \rf{onbad-} is trivially true due to already proved inclusion $\cI_{\ell_*}^-\subset J_{\ell_*}^-$; the second inequality is evident from the definitions of $\mr_i$'s and $\delta_{ij}$'s (recall that we have assumed all $B_i$ to be nonempty). Finally, \rf{onbad} is an immediate consequence of inclusions $\cI_{\ell_*}^-\subset \cI\backslash\{\ell_*\}\subset J_{\ell_*}\subset \overline{J}$ (the first and the third are evident, the second has been proved) and the definition of $\mr$'s and $\delta_{ij}$'s.  \qed

\subsubsection{Proof of Theorem \ref{thm:genthm}}
\paragraph{1$^o$.}
Let $\varepsilon=\epsilon/(N-1)$. Consider a pair $(i,j)\in\cO$ which is {bad}. In this case, one has $\delta_{ij}\leq \myr^{K}_{ij}(\varepsilon)$.
Indeed,  consider optimization problem
\be
\max_{x\in B_i,y\in B_j}\varrho(\cA_i(x),\cA_j(y)).
\ee{1rhoprob}
Observe that problem \rf{1rhoprob} is solvable, and its optimal solution $x'\in B_i,\,y'\in B_j$ satisfies
$\|x'-y'\|\geq 2\delta_{ij}$. On the other hand, the optimal value $\bar\rho$ of \rf{1rhoprob} is greater than
$\varepsilon^{1/K}$ because, otherwise, the risk of a $K$-observation test $\T_{\{i,j\}}$ deciding on hypotheses $H_i$ and $H_j$, as discussed in Section \ref{sectpairs}, would be bounded by
$\bar\rho^K=\varepsilon$, implying that pair ${(i,j)}$ is good what is not the case. We conclude that $\myr^{K}_{ij}(\varepsilon)$, as defined in \rf{Deltaij}, satisfies $\myr^{K}_{ij}(\varepsilon)\geq \half \|x'-y'\|\geq \delta_{ij}.$ Combined with (\ref{neweq999}) and the bounds \rf{1stb} and \rf{onbad-}, the latter relation implies \rf{thm1eq1}.

\paragraph{2$^o$.} To show \rf{thm1eq2} we need the following statement.
\begin{lemma}\label{lem:Dij}
Given $(i,j)$ with $1\leq j<i\leq N$, let   $\mR_{ij}=\RiskOpt^{\{i,j\}}_{\epsilon,\overline K}[X_i\cup X_j]$ be the minimax $\overline K$-observation $\epsilon$-risk of estimation over $X_i\cup X_j$. Suppose that $\varepsilon\in(0,1/2)$ is such that
\[
{\wt\vartheta}={\ln\left(4\epsilon(1-\epsilon)\right)\over 2\ln(\varepsilon)}\leq 1
\]
\item {\rm (i)} {Assume that}
$ K>  {\wt\vartheta^{-1}}\overline{K}.
$ 
Then $\myr^{K}_{ij}(\varepsilon)\leq \mR_{ij}$. 
\item {\rm (ii)}
In addition, if $\cA_i(x_{ij})=\cA_j(x_{ij})$ for some $x_{ij}\in X_i\cap X_j$, one has
\[
\myr^{K}_{ij}(\varepsilon)\leq \wt\vartheta^{-1} \mR_{ij}
\]
whenever $K\geq \overline{K}$.
\end{lemma}
\noindent {\bf Proof.}
When problem (\ref{Deltaij}) is infeasible, we have $\myr_{ij}^K(\varepsilon)=0$, and the claims in Lemma are trivially true. Now let \rf{Deltaij} be feasible. Then the problem is solvable; let $(\bar x,\bar y)$, $\bar x\in X_i, \bar y\in X_j$, be an optimal solution.
Suppose that $\mR_{ij}< \half \|\bar x-\bar y\|$. This would imply existence of a $\overline K$-observation estimate $\tilde x(\cdot)$ with maximal $\epsilon$-risk over $X_i\cup X_j$ which is smaller than $\half \|\bar x-\bar y\|$, meaning that there is a simple {$\overline{K}$-observation} test deciding on hypothesis $H_{\bar x}: \mbox{``observation $\omega^{\overline{K}}$ stems from $(i,\bar x)$''}$ against $H_{\bar y}:\mbox{``observation $\omega^{\overline{K}}$ stems from $(j,\bar y)$''}$ with risk bounded with $\epsilon$, namely, the test which accepts $H_{\bar x}$ whenever $\|\tilde x-\bar x\|\leq \|\tilde x-\bar y\|$ and accepts $H_{\bar y}$ otherwise. By what we know about testing in simple observation schemes, this means that Hellinger affinity $\varrho(\cA_i(\bar x),\cA_j(\bar y))$ between the corresponding distributions of observations satisfies (cf. \rf{epsilonstar})
$ \varrho(\cA_i(\bar x),\cA_j(\bar y))\leq [4\epsilon(1-\epsilon)]^{1/(2\overline K)}< \varepsilon^{1/K}$ contradicting the fact that, by construction of $\bar x$ and $\bar y$, $\varrho(\cA_i(\bar x),\cA_j(\bar y)\geq  \varepsilon^{1/K}$. {(i) is proved.}
%

Next, to prove (ii), for $\vartheta\in[0,1]$, let $x(\vartheta)=x_{ij}+\vartheta(\bar x-x_{ij})$
and $y(\vartheta)=x_{ij}+\vartheta(\bar y-x_{ij})$; observe that $\tilde\rho(\vartheta):=\ln \varrho(\cA_i(x(\vartheta)),\cA_j(y(\vartheta)))$ is a concave function of $\vartheta$ with $\tilde\rho(1)\geq
K^{-1}\ln \varepsilon$ and $\tilde\rho(0)=0$. Thus, for any $\vartheta<\wt\vartheta$,
\[
\tilde\rho(\vartheta)\geq\vartheta{K^{-1}}\ln \varepsilon>{\wt\vartheta}{K^{-1}}\ln \varepsilon
\geq\half{K^{-1}}\ln[4\epsilon(1-\epsilon)].
\]
As we already know, this means that there is no {$K$-observation} test capable of deciding between hypotheses $H_{x(\vartheta)}$ and $H_{y(\vartheta)}$ with risk bounded with $\epsilon$, implying in its turn that
$$
\mR_{ij}>\half\|x(\vartheta)-y(\vartheta)\|=\half\vartheta\|\bar x-\bar y\|=\vartheta\myr^{K}_{ij}(\varepsilon).\eqno{\mbox{\qed}}
$$
Setting $\varepsilon=\epsilon/(N-1)$ (which results in $\wt\vartheta=\ov\vartheta$) and substituting into \rf{thm1eq1} bounds of  Lemma \ref{lem:Dij} we arrive at \rf{thm1eq2} and \rf{thm1eq3}. \qed

\subsubsection{Proof of Proposition \ref{prop:l2aggreg1}}
 {\bf 1$^o$.} Observe first that the ``true hypothesis'' $H_{i_*j}^{\ell_*-}$ in quadruple $(i_*,j;\ell_*,\ell)$ is never empty because $w_*=Bx_*\in W^{\ell_*-}_{i_*j}$ for all $j\neq i_*$. Furthermore, whenever one of the hypotheses  $H^{\ell-}_{ij}$, $H_{ij}^{\ell\ell'+}$ is true in a good quadruple $(i,j;\ell,\ell')$,  test $\T^{\ell\ell'}_{ij}$ will accept it with probability at least $1-{\epsilon{/\ov N}}$. Indeed, we have assumed that this is the case if both hypotheses are nonempty; since a hypothesis, when true, cannot be empty, the only other case to be considered is that of the other hypothesis in the pair being empty. It remains to recall that in this case the test always accepts the nonempty hypothesis.
Thus, the $(\ell_*,x_*)$-probability of $\overline\Omega^K$  is $\geq 1-\epsilon$ due to the union bound.
\\{\bf 2$^o$.}
From now on, let $\omega^K\in \overline \Omega^K$; in this case we have $(i_*;\ell_*)\in \I(\omega^K)$, implying that $\I(\omega^K)\neq\emptyset$; thus, if all pairs ${(i,\ell)}\in \I(\omega^K)$ share the same $i$-component $\wh i=\wh i(\omega^K)$ we clearly have $\wh i(\omega^K)=i_*$. Next, suppose that $(i';\ell')\in \I(\omega^K)$ with $i'\neq i_*$.
Observe that for all $j\neq i_*$ one has
\[
\|x_*-x_j\|\leq \|x_*-x_{i_*}\|+\|x_{i_*}-x_j\|\leq  \|x_*-x_{i_*}\|+2r_{ij}.
\]
We conclude that whenever quadruple $(i',i_*;\ell',\ell_*)$ is bad one has
\[
\|x_*-{x_{i'}}\|\leq \|x_*-x_{i_*}\|+2\delta_{i'i_*}^{\ell'\ell_*}.
\]
Let us now fix a good quadruple $(i',i_*;\ell',\ell_*)$.
We have
\[0\leq \psi_{i'i_*}^T(w_*-w_{i'i_*})<\delta^{\ell'\ell_*}_{i'i_*}\]
where the first inequality is due to
\[
\|w_*-w_{i'}\|_2=\|x_*-x_{i'}\|\geq \|x_*-x_{i_*}\|=\|w_*-w_{i_*}\|_2,
\] while the second one is due to the fact that were it false, the hypothesis $H^{\ell'+}_{i'i_*}(\delta^{\ell'\ell_*}_{i'i_*})$ would be true and thus $H^{\ell'-}_{i'i_*}$ would be rejected by  test $\T^{\ell'\ell_*}_{i'i_*}$ (recall that $\omega^K\in\overline{\Omega}^K$), which is not the case because $(i';\ell')\in \I(\omega^K)$. Denoting by $\pi_{i_*i'}$ the projection of $w_*$ onto the line 
passing through $w_{i'}$ and $w_{i_*}$, let
 $\tau_*=\tau(w_{i_*})$, $\tau_{\pi}=\tau(\pi_{i_*i'})$, and $\tau'=\tau(w_{i'})$ be coordinates of
$w_{i_*}$, $\pi_{i_*i'}$, and $w_{i'}$ on this line, the origin on the line being the midpoint $w_{i_*i'}$ of the segment $[w_{i_*}, w_{i'}]$, its orientation given by $\psi_{i'i_*}$. One has $\tau_*={r_{i_*i'}}$, $\tau'={-r_{i_*i'}}$, and $\tau_\pi{\leq} \delta^{\ell'\ell_*}_{i'i_*}$, and so
\bse
\|x_*-x_{i'}\|^2-\|x_*-x_{i_*}\|^2&=&\|w_*-w_{i'}\|^2_2-\|w_*-w_{i_*}\|_2^2=\|\pi_{i_*i'}-w_{i'}\|^2_2-\|\pi_{i_*i'}-w_{i_*}\|_2^2
\\&=&(\tau_{\pi}+r_{i_*i'})^2-(\tau_{\pi}-r_{i_*i'})^2=4r_{i_*i'}\tau_\pi\leq 4r_{i_*i'}\delta_{i'i_*}^{\ell'\ell_*}.
\ese
We conclude that
\[\|{x_*}-x_{i'}\|-\|{x_*}-x_{i_*}\|={\|{x_*}-x_{i'}\|^2-\|{x_*}-x_{i_*}\|^2\over \|{x_*}-x_{i'}\|+\|{x_*}-x_{i_*}\|}\leq {\|{x_*}-x_{i'}\|^2-\|{x_*}-x_{i_*}\|^2\over 2r_{i_*i'}}\leq 2\delta_{i'i_*}^{\ell'\ell_*}
\]
implying \rf{ea0}. 
\paragraph{3$^o$.} Denote
\[
R=\half \max_{(i;\ell),(j;\ell')\in \I(\omega^K)} \|w_i-w_j\|_2,
\]
and let $w_{\bar i}$ and $w_{\bar j}$ be the endpoints of a maximizing segment with $w_{\bar{i}\bar{j}}=\half (w_{\bar i}+w_{\bar j})$ being its midpoint and $\wh x(\omega^K)=\half (x_{\bar i}+x_{\bar j})$ being the aggregated solution yielded by our algorithm.
W.l.o.g. assume that
$\|w_{\bar{i}}-w_*\|_2\leq \|w_{\bar{j}}-w_*\|_2$, implying that $(w_*-w_{\bar i\bar j})^T(w_{\bar{j}}-w_{\bar{i}\bar j})\leq 0$.  We have $\|w_{\bar j}-w_{i_*}\|_2\leq 2R$, whence, as we have just established,
\[
\|w_{\bar{j}}-w_*\|_2^2-\|{w_{i_*}}-w_*\|_2^2\leq 2\|w_{\bar j}-w_{i_*}\|_2\max_{\ell:(j;\ell)\in \I(\omega^K)}\delta^{\ell'\ell_*}_{ ji_*}\leq 4R\wt \delta^{\ell_*}_{i_*}(\omega^K).
\]
On the other hand,
\[
\|w_{\bar{j}}-w_*\|_2^2-\|w_{\bar{i}\bar{j}}-w_*\|_2^2=2(w_{\bar i\bar j}-w_*)^T(w_{\bar{j}}-w_{\bar{i}\bar{j}})+\|w_{\bar{j}}-w_{\bar{i}\bar{j}}\|_2^2\geq \|w_{\bar{j}}-w_{\bar{i}\bar{j}}\|_2^2=R^2,
\]
 and we conclude that
\bse
\|\wt x(\omega^K) -x_*\|^2&=&\|w_{\bar{i}\bar{j}}-w_*\|_2^2\leq \|w_{\bar{j}}-w_*\|_2^2-R^2\leq \|{w_{i_*}}-w_*\|_2^2+4R\wt\delta^{\ell_*}_{i_*}(\omega^K)-R^2\\&\leq& \|{w_{i_*}}-w\|_2^2+4\wt\delta^{\ell_*}_{i_*}(\omega^K)^2.
\ese
what is \rf{stupid0}.\qed
\subsubsection{Proof of Theorem \ref{thm:genthm2}}
{\bf 1$^o$.} Suppose that quadruple $(i,j;\ell,\ell')$ is bad. Let us verify that in this case one has $r_{ij}\leq \myr^{K}_{\ell\ell'}(\varepsilon)$ where  (cf. \rf{Deltaij}, \rf{thm1eq121})
\begin{equation}\label{arik}
\varepsilon={\epsilon{/\ov N}},\;\;\myr^{K}_{\ell\ell'}(\varepsilon)=\half \max_{x\in X_\ell,y\in X_{\ell'}}\left\{\|x-y\|:\,\varrho(\cA_\ell(x),\cA_{\ell'}(y))\geq\varepsilon^{1/K}
\right\}.
\end{equation}
To this end, consider optimization problem
\be
\begin{array}{c}
{
\max_{x\in X^{\ell-}_{ij},y\in X^{\ell'+}_{ij}(\delta)}\varrho(\cA_{\ell}(x),\cA_{\ell'}(y))},\\
{X^{\ell-}_{ij}=\{x\in X_\ell:\,Bx\in W^{\ell-}_{ij}\},\;\;X^{\ell'+}_{ij}(\delta)=\{x\in X_{\ell'}:\,Bx\in W^{\ell'+}_{ij}(\delta)\}}\\
\end{array}
\ee{2rhoprob}
for $\delta=r_{ij}$.
Note that ${X^{\ell-}_{ij}}$ and ${X^{\ell'+}_{ij}(r_{ij})}$ are nonempty (otherwise the corresponding quadruple would be $\varepsilon$-good) convex and compact sets. Thus, problem \rf{2rhoprob} is solvable, and its optimal solution $x'\in X^{\ell-}_{ij},y'\in X^{\ell'+}_{ij}(\delta)$ satisfies
{$\|x'-y'\|=\|Bx'-By'\|_2\geq \delta=r_{ij}$}. On the other hand, optimal value $\bar\rho$ of \rf{2rhoprob} is greater than $\varepsilon^{1/K}$ because, otherwise, the risk of a $K$-observation test $\T^{\ell\ell'}_{ij}$ deciding on hypothesis $H^{\ell-}_{ij}$ against
$H^{\ell'+}_{ij}(r_{ij})$, as built in Section \ref{sectpairs}, would be bounded by
$\bar\rho^K=\varepsilon$, so the quadruple $(i,j;\ell,\ell')$ would be $\varepsilon$-good what is not the case. In other words, $x',y'$ is a feasible solution to the maximization problem in (\ref{arik}) with the value of the objective $\geq r_{ij}$, implying the desired inequality $r_{ij}\leq \myr^{K}_{\ell\ell'}(\varepsilon)$.

Next, assume that quadruple $(i,j;\ell,\ell')$ is good and that $\delta^{\ell\ell'}_{ij}>\underline\delta$. In this case set ${X^{\ell'+}_{ij}(\delta^{\ell\ell'}_{ij}-\underline\delta)}$ is not empty because $\delta=\delta^{\ell\ell'}_{ij}-\underline\delta$ would be $\varepsilon$-good otherwise, and we know it is not.
Same as above, we conclude that in this case $\delta^{\ell\ell'}_{ij}-\underline\delta\leq \myr^{K}_{\ell\ell'}(\varepsilon)$, implying that whether quadruple $(i_*,j;\ell_*,\ell)$ is good or bad, one has
\[
\delta^{\ell_*\ell}_{i_*j}\leq \myr^{K}_{\ell_*\ell}(\varepsilon)+\underline\delta,
\]
  so that
  \begin{subequations}
  \begin{align}
  \wh{\delta}^{\ell_*}_{i_*}(\omega^K)&\leq \max_{\ell\leq \ell_*}\myr^{K}_{\ell_*\ell}(\varepsilon)+\underline\delta,\label{whdel.a}\\
  \wt{\delta}^{\ell_*}_{i_*}(\omega^K)&\leq \max_{\ell}\myr^{K}_{\ell_*\ell}(\varepsilon)+\underline\delta.\label{whdel.b}
\end{align}
  \label{whdel}
\end{subequations}
 \\{\bf 2$^o$.} Now 
 \rf{whdel.a} combined with bound \rf{ea0} imply that whenever $\omega^K\in \overline \Omega^K$
 \[
 \|x_*-\wh x(\omega^K)\|\leq \|x_*-x_{i_*}\|+2\left[\max_{\ell\leq \ell_*}\myr^{K}_{\ell_*\ell}(\varepsilon)+\underline\delta\right]
 \]
 (recall that $\|x_*-x_{i_*}\|\leq \|x_*-x_{\ell_*}\|$ by construction and $\|x_*-x_{\ell_*}\|=\|x_*-\wt x_{\ell_*}(\wt\omega^{\ov K})\|\leq \myr^{\ov K}_{\ell_*}(\epsilon)$ due to $\wt\omega^{\ov K}\in \wt \Omega^{\ov K}$).
 Utilizing the bound in \rf{whdel.b} we conclude that
 $$
 \|x_*-\wt x(\omega^K)\|^2\leq \|x_*-x_{i_*}\|^2+4\left[\max_{\ell}\myr^{K}_{\ell_*\ell}(\varepsilon)+\underline\delta\right]^2
 .
 $$
Finally, the second part of the statement of the theorem (starting with ``Consequently...'') is readily given by \rf{thm1eq121} and \rf{thm1eq122} combined with the result of Lemma \ref{lem:Dij} applied with $\varepsilon={\epsilon{/\ov N}}$.\qed

\subsection{Proofs for Sections \ref{sect:general} and \ref{sec:simplebs5}}
\subsubsection{Proof of Proposition \ref{pro:3p}}
{\bf 1$^o$.}
The fact that the ${x_*}$-probability of $\overline\Xi$ is at least $1-\epsilon$  is readily given by the union bound and the fact that when $i,j$ are comparable and ${x_*}\in X_{ij}(\Delta_{ij})$, the ${x_*}$-probability for $\T_{ij}$ (when $i<j$) or $\T_{ji}$ (when $i>j$) to accept $H_{ij}$ is at least $1-\epsilon/(N-1)$.
\paragraph{2$^o$.}
Let us fix $\xi\in\overline\Xi$ and set $\wh{x}=\wh{x}(\xi)$, $\wh{i}=\wh{i}(\xi)$, $\rho_*=\|{x_*}-x_{i_*}\|$. We claim that whenever $i_*$ {loses} to $j\neq i_*$, we have
\begin{equation}\label{eq4444}
\rho_*\geq r_{i_*j}-\Delta_{i_*j}
\end{equation}
Indeed, let $j\neq i_*$ be such that  $i_*$ loses to $j$. If $j$ is comparable to $i_*$, we have ${x_*}\not\in X_{i_*j}(\Delta_{i_*j})$. Indeed, otherwise the  test $\T_{\{i_*,j\}}$  would accept $H_{i_*j}$ due to $\xi\in\overline{\Xi}$ and $j$ would loose to $i_*$, which is not the case. On the other hand,  ${x_*}\not\in X_{i_*j}(\Delta_{i_*j})$ is exactly the same as
$\rho_*=\|{x_*}-x_{i_*}\|>r_{i_*j}-\Delta_{i_*j}$. Now, let $j$ and $i_*$ be incomparable; in this case  $i_*$ loosing to $j$ means that $\rho_j\leq\rho_{i_*}$, that is, $\Delta_{i_*j}=\max[0,r_{i_*j}-\rho_{i_*}]$, implying that \[r_{i_*j}-\Delta_{i_*j}\leq \rho_{i_*}\leq\rho_*\] (the concluding $\leq$ being given by $\rho_{i_*}=\min_{x'\in X}\|x'-x_{i_*}\|$ combined with $\rho_*=\|{x_*}-x_{i_*}\|$  and $x\in X$). 
\paragraph{3$^o$.}
Note that if $\widehat{i}=i_*$ (\ref{bound0}) clearly is true. Let now $i_*\neq\widehat{i}$. Then, if $i_*$ {loses} to no $j$ we would have $d_{i_*}=-\infty$, and since every $i\neq i_*$ loses to $i_*$, $d_i\geq0$ for all $i\neq i_*$, resulting in $\widehat{i}=i_*$ which is not the case.
Let us assume that $\widehat{i}\neq i_*$ and that $i_*$ loses to some $j$'s; let also $j_*=j_*(\xi)\in \cI_{i_*}(\xi)$ be such that $\|x_{i_*}-x_{j_*}\|=d_{i_*}(\xi)$. There are two possibilities:
\begin{itemize}
\item $i_*$ loses to $\widehat{i}$; when it is the case, (\ref{eq4444}) says that $\rho_*\geq r_{i_*\widehat{i}}-\Delta_{i_*\widehat{i}}$, whence
$$\|{x_*}-x_{\widehat{i}}\|\leq \|{x_*}-x_{i_*}\|+\|x_{i_*}-x_{\widehat{i}}\|=\rho_*+2r_{i_*\widehat{i}}\leq \rho_*+2[\rho_*+\Delta_{i_*\widehat{i}}]=3\rho_*+2\Delta_{i_*\widehat{i}},$$
and (\ref{bound0}) follows.
\item $\widehat{i}$ loses to $i_*$, implying that
$$\|x_{\widehat{i}}-x_{i_*}\|\leq d_{\widehat{i}}(\xi)\leq d_{i_*}(\xi)= \|x_{i_*}-x_{j_*}\|.
$$
Since $i_*$ loses to $j_*$, we have  $\rho_*\geq r_{i_*j_*}-\Delta_{i_*j_*}$ due to (\ref{eq4444}), resulting in
\bse
\|{x_*}-x_{\widehat{i}}\|&\leq& \|{x_*}-x_{i_*}\|+\|x_{i_*}-x_{\widehat{i}}\|\leq \|{x_*}-x_{i_*}\|+\|x_{i_*}-x_{j_*}\|
\\
&=&\rho_*+2r_{i_*j_*}\leq \rho_*+2[\rho_*+\Delta_{i_*j_*}],
\ese\end{itemize}
and (\ref{bound0}) follows.\qed
\subsubsection{Proof of Proposition \ref{prop:l2aggreg10}}
The fact that ${x_*}$-probability of $\overline\Xi$ is at least $1-\epsilon$ is obvious 
(cf. the proof of Proposition \ref{pro:3p}). Now let us fix $\xi\in\overline\Xi$ and let
 $\wh{x}=\wh{x}(\xi) $, $\wh{i}=\wh{i}(\xi)$, and $\rho_*=\|{x_*}-x_{i_*}\|$. Consider the following ``coloring'' of indices $1\leq j\leq N$:
\begin{itemize}
\item $j$ is white if $\|{x_*}-x_j\|\leq \rho_*+\overline{\Delta}$;
\item $j$ is gray if ${\rho_*}+\overline{\Delta}<\|{x_*}-x_j\|\leq \rho_*+2\overline{\Delta}$;
\item $j$ is black if $\|{x_*}-x_j\|>\rho_*+2\overline{\Delta}$.
\end{itemize}
Let $k_w$, $k_g$, $k_b$ be the numbers of white, gray, and black indices, respectively.
Recalling that  $\xi\in\overline \Xi$, observe that
\begin{itemize}\item
 When $j$ is gray or black,
 \[\|{x_*}-x_j\|>\rho_*+\overline{\Delta}=\|{x_*}-x_{i_*}\|+\overline{\Delta}\geq \|{x_*}-x_{i_*}\|+\Delta_{i_*j},\] that is, ${x_*}\in \X_{i_*j}(\Delta_{i_*j})$. It follows that as applied to observation $\xi$, the test $\T_{\{i_*,j\}}$
of hypotheses $\cH_{i_*j}$ and $\cH_{ji_*}$ accepts $\cH_{i_*j}$, that is, $s_{i_*}(\xi)\leq k_w-1$.
\\
\item
When index $i$ is black and $j$ is a white index, we have \[
\|{x_*}-x_i\|>\rho_*+2\overline{\Delta} \geq \overline{\Delta}+\|{x_*}-x_j\|,
\] that is, ${x_*}\in\X_{ji}(\Delta_{ij})$. As a consequence, as applied to observation $\xi$, the test $\T_{\{i,j\}}$ of hypotheses $\cH_{ij}$ and $\cH_{ji}$ accepts the second hypothesis, implying that $s_i(\xi)\geq k_w$.\end{itemize}
Taken together, the above observations say that when $\xi\in\overline\Xi$ stems from ${x_*}$, index $\widehat{i}(\xi)$ is either white or gray, but definitely is not black, implying that
$$
\|\wh x-{x_*}\|\leq \rho_*+2\overline{\Delta}.
\eqno{\mbox{\qed}}
$$
\subsubsection{Proof of Theorem \ref{thm:lowerbas}}
{\bf 1$^o$.} Given a pair $(i,j)\in{\cal O}$ such that $r_{ij}=\half \|x_i-x_j\|> \delta$, let us set
\[
X_{ij}(\lambda,\delta)=\Big\{z\in X:\,\|z-x_i\|\leq  \underbrace{\lambda(r_{ij}-\delta)}_{d_{ij}}\Big\}
\]
where $0<\lambda<(1+\gamma)^{-1}$.
Under the premise of Theorem, for any such pair $i,j$ there exists a $\overline{K}$-observation test deciding with risk $\leq\epsilon$ on a pair of hypotheses ${\overline{H}}_{ij}$ and ${\overline{H}}_{ji}$ stating, respectively, that
$\omega^{\overline{K}}$ stems from signal  belonging to $X_{ij}(\lambda,\delta)$ and  $X_{ji}(\lambda,\delta)$, and both these sets are nonempty (recall that we are in the case where $x_s\in X$, $1\leq s\leq N$).
 The desired test $\T$ is as follows: given observation $\omega^{\overline{K}}$ we compute $\bar x(\omega^{\overline{K}})$ and accept
${\overline{H}}_{ij}$ when
$\|\bar x(\omega^{\overline{K}})-x_{i}\|\leq \|\bar x(\omega^{\overline{K}})-x_{j}\|$, and accept ${\overline{H}}_{ji}$ otherwise.
\par Let us verify that the risk of this test is indeed at most  $\epsilon$. Suppose, first, that ${\overline{H}}_{ij}$ takes place, so that $\omega^{\overline{K}}\sim p_{\cA_\nu({x_*})}^{\overline{K}}$ for some $\nu\leq J$ and ${x_*}\in X_{ij}(\lambda(r_{ij}-\delta))$. Then, if $x_{i_*}$ is the closest to $x_*$ point among $x_1,...,x_N$, we have
\[
\|x_*-x_{i_*}\|\leq \|x_*-x_i\|\leq \lambda(r_{ij}-\delta),
\]
and so
$p_{\cA_\nu(x)}^{\overline{K}}$-probability of the event
\[\E=\{\omega^{\ov K}:\,\|\bar x(\omega^{\overline{K}})-x\|\leq \gamma d_{ij}+\delta\}
\]
 is at least $1-\epsilon$ due to the origin of $\bar x(\cdot)$. But if $\E$ takes place,
\[
\|\bar x(\omega^{\overline{K}})-x_i\|\leq\|\bar x(\omega^{\overline{K}})-{x_*}\|+\|x_i-{x_*}\|\leq (\gamma+1) d_{ij}+\delta< r_{ij},
\]
so that
\[
\|\bar x(\omega^{\overline{K}})-x_j\|\geq \|x_i-x_j\|-\|\bar x(\omega^{\overline{K}})-x_i\|> r_{ij}.
\]
We conclude that
the $p^{\overline{K}}_{\cA_\nu({x_*})}$-probability for $\T$ not to accept ${\overline{H}}_{ij}$ is $\leq\epsilon$. By ``symmetric reasoning,'' when ${\overline{H}}_{ji}$ holds true, so that
$\omega^K\sim p^{\overline{K}}_{\cA_\nu({x_*})}$ for some $\nu\leq J$ and ${x_*}\in X_{ji}(\lambda(r_{ij}-\delta))$,  $p^{\overline{K}}_{\cA_\nu({x_*})}$-probability to reject ${\overline{H}}_{ji}$ is at most $\epsilon$.
\par
Now, testing ${\overline{H}}_{ji}$ against ${\overline{H}}_{ij}$ is equivalent to deciding between ``red'' set $R_{ij}(\lambda,\delta)$ and ``blue'' set  $B_{ij}(\lambda,\delta)$ in the space $\cM$ of parameters of distribution $p^{\overline{K}}_\mu$ of $\omega^{\overline{K}}$, each set being  a union of at most $J$ convex and compact sets:
\[
R_{ij}(\lambda,\delta)=\bigcup_{\nu=1}^JR_{ij\nu}(\lambda,\delta),\;R_{ij\nu}(\lambda,\delta)=\{\cA_\nu(x):\,x\in X_\nu,\,\|x-x_i\|\leq \lambda(r_{ij}-\delta)\},
\]
and
\[
B_{ij}(\lambda,\delta)=\bigcup_{\nu=1}^JB_{ij\nu}(\lambda,\delta),\;B_{ij\nu}(\lambda,\delta)=\{\cA_\nu(x):\,x\in X_\nu,\,\|x-x_j\|\leq \lambda(r_{ij}-\delta)\},\;\nu=1,...,J.
\]
From what we know about color inferring test in simple observation schemes, the fact that the hypotheses ${\overline{H}}_{ij}$ and ${\overline{H}}_{ji}$ can be decided upon via $\overline{K}$-repeated observation $\omega^{\overline{K}}\sim p^{\overline{K}}_{\cA_\nu(x)}$ with risk $\epsilon\in(0,1/2)$ implies (cf. Proposition \ref{col:basic}) that
when $K$ satisfies (\ref{eqKl}),
we have at our disposal test ${\T}_{ij}$ utilising $K$-repeated observation $\omega^K$ which decides {with maximal risk not exceeding $\epsilon/(N-1)$} upon hypotheses ${H}_{ij}$ and ${H}_{ji}$ stating that $\omega^K$ stems from a signal  $x$ such that $x\in X_{ij}(\lambda(r_{ij}-\delta))$ (for ${H}_{ij}$) and $x\in X_{ji}(\lambda(r_{ij}-\delta)$ (for ${H}_{ji}$).
\par\noindent
{\bf 2$^o$.} Now let us apply the aggregation procedure described in Section \ref{sec:genorma} to $K$-repeated observations, with $K$ satisfying (\ref{eqKl}).  From what we have just seen, 
in this case, all pairs $(i,j)$  such that $r_{ij}> \delta$  are appropriate, and  the quantity $(1-\lambda)r_{ij}+\lambda\delta$ with $\lambda\in D:=(0,(1+\gamma)^{-1})$ is $(i,j)$-appropriate. Let us set
\[
\widetilde{\lambda}=(1+\gamma')^{-1}\in D, \,
\widetilde{d}_{ij}=\widetilde{\lambda}(r_{ij}-\delta), \, \widetilde{\delta}_{ij}=r_{ij}-\widetilde{d}_{ij}=(1-\widetilde{\lambda})r_{ij}+\widetilde{\lambda}\delta,
\]
so that $\Delta_{ij}:=\widetilde{\delta}_{ij}$ is $(i,j)$-appropriate, as required in the construction we are implementing.  Let us define the set of comparable pairs to be exactly the set of pairs $\{i,j\}\in\cU$ with $r_{ij}>\delta$ and equip these pairs with the tests $\T_{\{i,j\}}=\T_{\min[i,j],\max[i,j]}$.  For these pairs the quantities $\Delta_{ij}$  satisfy the relations
\begin{equation}\label{Deltaijsmall}
\Delta_{ij}=(1-\widetilde{\lambda})r_{ij}+\widetilde{\lambda}\delta= {\gamma'\over1+\gamma'}r_{ij}+{1\over 1+\gamma'}\delta.
\end{equation}
 Note that for incomparable pairs $(i,j)\in{\cO}$ we have
 $\Delta_{ij}=\max[0,r_{ij}-\max[\rho_i,\rho_j]]\leq r_{ij}\leq\delta$.
\par\noindent
{\bf 3$^o$.}
Now let observation $\omega^K$ stem from signal ${x_*}\in X$, so that $\omega^K\sim p^K_{\cA_\nu({x_*})}$ for some $\nu$ such that ${x_*}\in X_\nu$, and let $i_*$ be the index of one of the $\|\cdot\|$-closest to ${x_*}$ points $x_1,...,x_N$. Finally, let $\overline\Omega$ be the set of $\omega^K$ satisfying the condition
\begin{quote}
whenever $j\neq i_*$ is such that $r_{i_*j}>\delta$ (i.e., $i_*$ and $j$ are comparable) and \[
\|{x_*}-x_{i_*}\|\leq r_{i_*j}-\Delta_{i_*j}\] (i.e., hypothesis $H_{i_*j}$ holds true), test $\T_{\{i_*,j\}}$ as applied to observation $\omega^K$ accepts $H_{i_*j}$.
\end{quote}
By Proposition \ref{pro:3p}, the $p^K_{\cA_\nu({x_*})}$-probability of $\overline{\Omega}$ is at least $1-\epsilon$, and
\begin{equation}\label{andrgh}
\|{x_*}-\widehat{x}(\omega^K)\|\leq 3\|{x_*}-x_{i_*}\|+2\overline{\Delta}_{i_*}(\omega^K).
\end{equation}
Next, let $\omega^K\in \overline{\Omega}$, and let $j\in \cI_{i_*}(\omega^K)$.
It may happen that $i_*$ and $j$ are comparable; in this case $H_{i_*j}$ cannot be true due to $\omega^K\in\overline{\Omega}$, that is, $\|{x_*}-x_{i_*}\|>r_{i_*j}-\Delta_{i_*j}$, and besides this,
\[
\Delta_{i_*j}\leq{\gamma'\over 1+\gamma'}r_{i_*j}+{1\over 1+\gamma'}\delta
\] due to (\ref{Deltaijsmall}), implying that $r_{i_*j}\leq (1+\gamma')\|{x_*}-x_{i_*}\|+\delta$. Hence,
\[
\Delta_{i_*j}\leq {\gamma'\over 1+\gamma'}r_{i_*j}+{1\over 1+\gamma'}\delta\leq \gamma'\|{x_*}-x_{i_*}\|+\delta.
\]
When $i_*$ and $j$ are incomparable, we have $\Delta_{i_*j}= r_{i_*j}\leq\delta$. We see that when $\omega^K\in \overline{\Omega}$ (what happens with $p^K_{\cA_\nu({x_*})}$-probability at least $1-\epsilon$), we have $\overline{\Delta}_{i_*}(\omega^K)\leq \gamma'\|{x_*}-x_{i_*}\|+\delta$. This combines with (\ref{andrgh}) to imply that when $\omega^K\in\overline{\Omega}$, we also have  $\|\widehat{x}-x_{i_*}\|\leq (3+2\gamma')\|{x_*}-x_{i_*}\|+2\delta$.
\qed
\subsubsection{Proof of Theorem \ref{thm:lowerbaseu}}
{\bf 1$^o$.} Let  $\ov \I$ be the set of pairs $\{i,j\}\in\cU$ such that both hypotheses $\cH_{ij}(\bar \delta)=H(\X_{ij}(\bar{\delta}))$ and $\cH_{ji}(\bar \delta)=H(\cX_{ji}(\bar{\delta}))$ are nonempty. Let $\{i,j\}\in \ov \I$ be fixed, and let us show that under the premise of the theorem
the simple test which, given $\omega^{\overline{K}}$, accepts $\cH_{ij}$ when  $\iota_{ij}(\omega^{\overline{K}})=i$, and accepts $\cH_{ji}$ otherwise has its risk bounded with $\epsilon$. Indeed, let $\cH_{ij}$ be true, that is, the distribution $P$ of observation $\omega^{\overline{K}}$ satisfies $P\in\cP_x^{\overline{K}}$ for some  $x\in \X_{ij}(\bar{\delta})$, so that
$\|x-x_j\|\geq\|x-x_i\|+\bar{\delta}$, whence $\|x-x_j\|\geq \min[\|x-x_i\|,\|x-x_j\|]+\bar{\delta}$. By (\ref{haim}) the $P$-probability of the event $\iota_{ij}=j$, that is, the probability  of the test in question rejecting $\cH_{ij}$, is $\leq\epsilon$.  By ``symmetric'' reasoning, the $P$-probability to reject $\cH_{ji}(\bar{\delta})$ when the hypothesis is true is $\leq\epsilon$ as well.
 \par
Now recall that testing $\cH_{ij}(\bar \delta)$ vs $\cH_{ji}(\bar \delta)$ via a $\overline{K}$-repeated observation is equivalent to deciding via this observation on ``red'' set $R_{ij}(\bar \delta)$ vs. ``blue'' set  $B_{ij}(\bar \delta)$ in the space $\cM$ of parameters of distribution $p_\mu$ of $\omega_k$, and each set is a union of at most $J$ convex and compact sets:
\[
R_{ij}(\bar \delta)=\bigcup_{\nu=1}^J R_{ij\nu}(\bar \delta),\; R_{ij\nu}(\bar \delta)=\{\cA_\nu(x):\,x\in X_\nu,\,\|x-x_i\|\leq \|x-x_j\|-\bar \delta)\},
\]
and
\[
B_{ij}(\bar \delta)=\bigcup_{\nu=1}^J B_{ij\nu}(\bar \delta),\;B_{ij\nu}(\bar \delta)=\{\cA_\nu(x):\,x\in X_\nu,\,\|x-x_j\|\leq \|x-x_i\|-\bar \delta\}.
\]
The fact that hypotheses $\cH_{ij}(\bar \delta)$ and $\cH_{ji}(\bar \delta)$ can be decided upon via $\overline{K}$-repeated observation  with risk $0\leq\epsilon<1/2$ implies by Proposition \ref{col:basic} that whenever
\[
K\geq \left\lceil {2\ln\left(J\overline{N}/\epsilon\right)\over \ln\left( [4\epsilon(1-\epsilon)]^{-1}\right)} \overline{K}\right\rceil,
\]
$\overline{\delta}$ is $(i,j)$-good in the sense of Section \ref{euclprelim}.
\par\noindent
{\bf 2$^o$.}
Now let $K$ satisfy (\ref{KLarge}) and $\{i,j\}\in \overline{\cI}$, that is, both  $\X_{ij}(\bar \delta)$ and  $\X_{ji}(\bar \delta)$ are nonempty. Recall that in this case in our aggregation procedure $\Delta_{ij}$ is selected to be $(i,j)$-good (that is, with $K$ observations, the test yielded by the machinery from Section \ref{sectinfcol} decides on the  hypothesis $\cH_{ij}(\Delta_{ij})$  vs. the alternative $\cH_{ji}(\Delta_{ij})$ with risk not exceeding $\epsilon/\overline{N}$)  and either $\Delta_{ij}\leq\underline{\delta}$, or $\Delta_{ij}-\underline{\delta}$ is not $(i,j)$-good.
By item 1$^o$, for our $i,j,K$ $\delta=\bar{\delta}$ is $(i,j)$-good, so that the second option implies that $\Delta_{ij}-\underline{\delta}\leq\bar{\delta}$ and one always has $\Delta_{ij}\leq \bar\delta+\underline{\delta}$.

\par
On the other hand, if $i\neq j$ and $\{i,j\}\not\in \ov \I$, at least one of  the sets $\X_{ij}(\bar \delta)$, $\X_{ji}(\bar \delta)$ is empty, implying that $\bar \delta$ is $(i,j)$-good. Consequently, in our aggregation procedure, same as in the case of $\{i,j\}\in\cI$, one has $\Delta_{ij}\leq \overline{\delta}+\underline{\delta}$. Thus, $\Delta_{ij}\leq\overline{\delta}+\underline{\delta}$ for all $i\neq j$, and \rf{finalbasiceu} is given by Proposition \ref{prop:l2aggreg10}.
\qed

\subsubsection{Proof of Proposition \ref{prop:louna}}
We start with the following observation.
\begin{lemma}\label{lem:lounr}
Under the premise of the proposition, let $\mR=\RiskOpt^{\overline{1,N}}_{\epsilon,\overline{K}}[{X}]$ be the minimax risk of $\overline{K}$-observation estimation over ${X}$. Let also $K$ satisfy \rf{thesameK} and $\{i,j\}\in \cU$ be such that $\mR<\overline{\delta}^{ij}$ (cf. (\ref{deltabarij})). Then  any $\delta$ such that $\mR<\delta\leq\overline{\delta}^{ij}$ is $(i,j)$-appropriate.
\end{lemma}
\noindent{\bf Proof.} Under the lemma's premise, for any $\rho>\mR$ there exists an estimate $\bar x=\bar x(\omega^{\overline K})$ such that for every $x\in {X}$, the $x$-probability of the event
$\|\bar x-x\|\leq \rho$ is at least $1-\epsilon$. As a result, for any $i\neq j$ and $\delta>\rho$ there exists a $\overline{K}$-observation test deciding on hypotheses  $H_{ij}(\delta)$ and $H_{ji}(\delta)$ with risk bounded with $\epsilon$, namely, test $\overline{T}_{\{i,j\}}$ accepting $H_{ij}(\delta)$ if
$\|\bar x-x_i\|\leq r_{ij}$ and accepting $H_{ji}(\delta)$ otherwise. Indeed, assuming that $H_{ij}(\delta)$ takes place, the distribution $P^{\ov K}$ of observation $\omega^{\ov K}$ stems from some $x\in X$ satisfying $\|x-x_i\|\leq r_{ij}-\delta<r_{ij}-\rho$, so that when the event $\|\bar{x}-x\|\leq\rho$ takes place
(which happens with $P^{\ov K}$-probability $\geq1-\epsilon$), we have
\[
\|\bar{x}-x_i\|\leq \|\bar{x}-x\|+\|x-x_i\|<r_{ij},
\]
and test $\overline{\T}_{\{i,j\}}$ accepts $H_{ij}(\delta)$. Similarly, when $H_{ji}(\delta)$ takes place, the distribution $P^{\overline{K}}$ of $\omega^{\ov K}$ stems from some $x\in X$ satisfying $\|x-x_j\|\leq r_{ij}-\delta<r_{ij}-\rho$, so that when the event $\|\bar{x}-x\|\leq\rho$ takes place (which happens with $P^{\overline{K}}$-probability $\geq1-\epsilon$), we have $\|\bar{x}-x_j\|\leq\|\bar{x}-x\|+\|x-x_j\|< \rho+r_{ij}-\rho=r_{ij}$, whence $\|\bar{x}-x_i\|>2r_{ij}-r_{ij}>r_{ij}$, and $\overline{\T}_{\{i,j\}}$ accepts $H_{ji}(\delta)$.
\par
Recalling that $X$ is the union of at most $N$ convex and compact sets, we conclude that when $K$ satisfies (\ref{thesameK}), the risk of the $K$-observation test deciding on $H_{ij}(\delta)$ vs $H_{ji}(\delta)$ constructed in Section \ref{sec:genorma} does not exceed $\epsilon/(N-1)$. \qed
The claim of the proposition is readily given by combining the bound (\ref{newnewnew}) with the fact that by Lemma \ref{lem:lounr} the quantity $\overline{\Delta}$, as is immediately seen, can be bounded by a quantity arbitrarily close to $\mR$.\qed

\subsubsection{Proof of Proposition \ref{prop:lounae}}
The statement of the proposition is readily implied by the following analog of Lemma \ref{lem:lounr}.
\begin{lemma}\label{lem:lounre}
Given a positive integer $\overline{K}$ and $\epsilon\in(0,1/2)$, let $\mR=\Risk_{\epsilon,\ov K}^{\overline{1,N}}[{X}]$ be the minimax $\epsilon$-risk of estimation over ${X}$, and let $K$ satisfy \rf{thesameKe}.
Then $\overline{\Delta}\leq  2\mR+\underline{\delta}$.
\end{lemma}
\paragraph{Proof.}
Let $\wt{\delta}>2\mR$ and let $(i,j)\in\cO$; let us show that $\wt{\delta}$ is $(i,j)$-good (for the definition of $(i,j)$-goodness, see Section \ref{sec:eunags}).
There is nothing to prove when at least one of the sets $\X_{ij}(\wt{\delta})$, $\X_{ji}(\wt{\delta})$ is empty.
Assuming these sets nonempty, let $\rho>\mR$ be such that $2\rho<\wt{\delta}$. Then there is an estimate $\bar{x}(\omega^{\overline{K}})$ such that for every $x\in X$, the $x$-probability of the event $\|x-\bar{x}(\omega^{\overline{K}})\|\leq \rho$ is $\geq1-\epsilon$.
We immediately convert this estimate into a $\overline{K}$-observation test deciding on the hypothesis $\cH_{ij}=H(\X_{ij}(\wt{\delta}))$ vs the alternative $\cH_{ji}=H(\X_{ji}(\wt{\delta}))$: given $\omega^{\overline{K}}$ and setting $\bar{x}=\bar{x}(\omega^{\overline{K}})$, this test accepts $\cH_{ij}$ (and rejects $\cH_{ji}$) when $\|\bar{x}-x_i\|\leq \|\bar{x}-x_j\|$, and accepts $\cH_{ji}$ (and rejects $\cH_{ij}$) otherwise.
Observe that the risk of this test is $\leq\epsilon$. Indeed, when $\cH_{ij}$ takes place, the distribution $P^{\ov K}$ of $\omega^{\ov K}$ stems from some $x\in \X_{ij}(\wt\delta)$, that is, $x\in X$ and $\|x-x_i\|\leq\|x-x_j\|-\wt\delta$. Therefore when  $\|x-\bar{x}\|\leq\rho$ (the latter happens with $P^{\ov K}$-probability $\geq1-\epsilon$), we have
\bse
\|\bar{x}-x_i\|&\leq&\|x-\bar{x}\|+\|x-x_i\|\leq \rho+\|x-x_j\|-\wt\delta\leq \rho+\|\bar{x}-x_j\|+\|\bar{x}-x\|-\wt\delta\\&\leq& \|\bar{x}-x_j\|+2\rho-\wt\delta<\|\bar{x}-x_j\|,
\ese
and the test accepts $\cH_{ij}$.
 ``Symmetric'' reasoning shows that when $\cH_{ji}$ takes place, the test accepts $\cH_{ji}$ and rejects $\cH_{ij}$  when $\|x-\bar{x}\|\leq\rho$, which happens with $x$-probability $\geq1-\epsilon$, implying that the risk of the test is $\leq\epsilon$.\par
Because $X$ is the union of $N$ convex compact sets, existence of pairwise $\overline{K}$-observation tests deciding with risk $\leq\epsilon$ on all pairs $\cH_{ij}$, $\cH_{ji}$ of nonempty hypotheses with $(i,j)\in\cO$ implies, by the results of Section \ref{sectinfcol}, that with $K$ as in \rf{thesameKe} $\wt{\delta}$ indeed is $(i,j)$-good for all $(i,j)\in\cO$.
\par
Now, for every $(i,j)\in\cO$ by construction the quantity $\Delta_{ij}$ is $(i,j)$-good and is either  $\leq\underline{\delta}$, or is such that $\Delta_{ij}-\underline{\delta}$ is not $(i,j)$-good.
By the above, the second option implies that $\Delta_{ij}<\wt{\delta}+\underline{\delta}$ for all $(i,j)\in\cO$, so that $\overline{\Delta}< \wt{\delta}+\underline{\delta}$.
The latter inequality holds true whenever $\wt{\delta}>2\mR$, and the conclusion of the lemma follows. \qed
\subsection{Proof of Theorem \ref{thm:genthmG}}
 Let us verify that in a {$K$-}bad pair ${(i,j)}$ $\delta_{ij}$, as defined in \rf{emeq101a}, satisfies $$\delta_{ij}\leq \myrg^{K}_{ij}(\varepsilon).$$
Indeed,  consider optimization problem
\be
\min_{x\in B_i,y\in B_j}\|\cA_i(x)-\cA_j(y)\|_2;
\ee{1rhoprobg}
observe that \rf{1rhoprobg} is solvable, and its optimal solution $x'\in B_i,\,y'\in B_j$ satisfies
$\|x'-y'\|\geq 2\delta_{ij}$. On the other hand, the optimal value of \rf{1rhoprob} is less than
${2\over \sqrt{K}}q_{\N}(1-\varepsilon)$ because, otherwise, the risk of {$K$}-observation test $\T_{\{i,j\}}$  deciding on hypotheses $H_i$ and $H_j$, see \rf{emeq101a}, as yielded in Gaussian case by the machinery from Section \ref{sectpairs}, would be bounded by
$\varepsilon$, implying that pair ${(i,j)}$ is {$K$-}good what is not the case. We conclude that $\myrg^{K}_{ij}(\varepsilon)$, as defined in \rf{DeltaijN}, satisfies $\myrg^{K}_{ij}(\varepsilon)\geq \half \|x'-y'\|\geq \delta_{ij}$; along with the result of Proposition \ref{ai1} (see\rf{1stb} and \rf{onbad-}) this implies relation (\ref{eq0999}).
\paragraph{2$^o$.}
Let us fix $(i,j)\in\cO$.
Let for $\upsilon\in(0,1)$ $\mR_{ij}(\upsilon)=\RiskOpt^{\{i,j\}}_{\upsilon,\overline K}[X_i\cup X_j]$; let also $(\bar x,\bar y)$, $\bar x\in X_i, \bar y\in X_j$, be an optimal solution to \rf{DeltaijN} with $\delta=\varepsilon$. Note that for $\vartheta\in[0,1]$, $x(\vartheta)=\vartheta\bar x\in X_i$
and $y(\vartheta)=\vartheta\bar y\in X_j$, while
$\rho(\vartheta)=\|\cA_i(x(\vartheta))-\cA_j(y(\vartheta))\|_2$ is a linear function of $\vartheta$ with
$\rho(1)\leq {2\over \sqrt{K}}q_{\N}(1-\varepsilon)$ and $\rho(0)=0$. Let now
\[
\wt\vartheta={\sqrt{ K} q_{\N}(1-\upsilon)\over \sqrt{\overline K} q_{\N}\Big(1-{\varepsilon}\Big)}\leq 1;
\]
then for $\vartheta<\wt \vartheta$ one has
\[
\rho(\vartheta)\leq \vartheta\rho(1)< {2\over \sqrt{\overline K}} q_{\N}(1-\upsilon).
\]
The latter relation means that there is no $\overline{K}$-observation test capable of deciding between hypotheses $
H_{x(\vartheta)}:\omega_k\sim\cN(A_ix(\vartheta),I_n)$ and $H_{y(\vartheta)}:\omega_k\sim\cN(A_jy(\vartheta),I_m)$ with risk bounded with $\upsilon$, implying in its turn that
$$
\mR_{ij}(\upsilon)>\half\|x(\vartheta)-y(\vartheta)\|=\half\vartheta\|\bar x-\bar y\|=\vartheta\myrg^{K}_{ij}(\varepsilon).
$$
Applying the latter bound to $\upsilon={\epsilon}$ (recall that $\wt\vartheta\leq1$ for $\upsilon=\epsilon$  due to (\ref{restrict})), we obtain
\[
\myrg^{K}_{ij}\left(\varepsilon\right)\leq 
\ov\vartheta^{-1} \mR_{ij}(\epsilon)
\]
which combines with (\ref{eq0999}) to imply \rf{barthe21}.

The same bound as applied with $\upsilon=1/16$ and $K=\ov K$ (this again is possible due to $\varepsilon\leq \epsilon<1/16$) implies that
\be
\myrg^{\ov K}_{ij}\left(\varepsilon\right)\leq {q_{\N}(1-\varepsilon)\over q_{\N}(\mbox{\small ${15\over 16}$})}
 \mR_{ij}(\mbox{\small ${1\over 16}$})\leq \myC_4\sqrt{\ln[N/\epsilon]}\RiskOpt^{\{i,j\}}_{{1\over 16},\overline K}[X_i\cup X_j].
 \ee{eq12222}
\paragraph{3$^o$.}
Thus, all we need to show the last statement of the theorem, is to bound the quantity ${\ov{\RiskOpt}^{\{j\}}_{{1\over 8},\ov K}}[X_j]$---the minimax $1/8$-risk of recovering $x\in X_j$ from single ``averaged'' observation
\[
\wh\omega\sim \N(\cA_j(x),\ov K^{-1}I_m).
\]
Common sense says that ${\ov{\RiskOpt}^{\{j\}}_{{1\over 8},\ov K}}[X_j]$ is exactly the same as ${\RiskOpt^{\{j\}}_{1/8,\ov K}}[X_j]$, but we do not know why this would be the case.\footnote{Recall that $\wh \omega$ is sufficient statistics when estimating functions of the mean of the Gaussian distribution---conditional distributions of $\omega^{\overline{K}}$ given $\wh \omega$ is Gaussian and does not depend on $x$. Were the considered loss convex, the corresponding result would be readily given by the Rao-Blackwell theorem.}
Instead, we  are about to establish a slightly weaker fact which is sufficient for our purposes.
\def\Riskopt{\mathrm{Riskopt}}
\begin{lemma}\label{lem:rand}
Suppose that for a positive integer $M$, $Y\subset \bR^n$, and $\upsilon\in (0,1/4)$ $\Riskopt_{\upsilon,M}[Y]$ is the minimax over $Y$ $\upsilon$-risk of estimation given an $M$-repeated observation $(\omega_1,...,\omega_M)$, $\omega_k\sim \N(\cA(x), I_m)$, $x\in Y$. Then minimax over $Y$ $2\upsilon$-risk $\overline{\Riskopt}_{2\upsilon,M}[Y]$ of estimation given single observation
\[\ov \omega={1\over M}
\sum_{i=1}^M\omega_k
\]
satisfies\be
\ov \Riskopt_{2\upsilon,M}[Y]\leq 2\Riskopt_{{\upsilon},M}[Y].
\ee{eq22122}
\end{lemma}
\noindent{\bf Proof.} Note that $\omega_k$, $k=1,...,M$ can be represented as $\omega_k=\eta_k+\ov\omega$ where $\eta_k\sim \N\Big(0,{M-1\over M}I\Big)$ are independent of $\ov\omega$. This observation implies that if  $\ov\Riskopt^{\cR}_{\upsilon,M}[X_j]$ is defined in the same fashion as $\ov{\Riskopt}_{\upsilon,M}[X_j]$ but with candidate estimates which may be {\em randomized} then
\[
\ov \Riskopt^\cR_{\upsilon,M}[X_j]\leq \Riskopt_{\upsilon,M}[X_j].
\]
We claim that
\be
\overline{\Riskopt}_{2\upsilon,M}[X_j]\leq 2\overline{\Riskopt}^\cR_{\upsilon,M}[X_j]
\ee{eq221}
what obviously implies the lemma.
Indeed, let $\rho>\overline{\Riskopt}^{\cR}_{\upsilon,M}$, so that there exists a deterministic function $\phi(\omega,\zeta)$ taking values in $\bR^n$ such that
for every $x\in Y$ it holds
$$
\Prob_{(\ov\omega,\zeta)\sim P}\{\|\phi(\ov\omega,\zeta)-x\|>\rho\}\leq \upsilon,
$$
where $P$ is the distribution of $(\ov\omega,\zeta)$ with independent of each other $\ov\omega\sim\cN(\cA(x),M^{-1}I_m)$ and $\zeta\sim
U$, $U$ being the uniform distribution over $[0,1]$. Let
$$
\ov\Omega=\big\{\ov \omega\in\bR^m: \exists y\in Y: \Prob_{\zeta\sim U}\{\zeta:\,\|\phi(\ov \omega,\zeta)-y\|\leq \rho\}>1/2\big\}.
$$
For every $\ov\omega\in\ov\Omega$, we can specify $\psi(\ov\omega)\in Y$ in such a way that
$$
\Prob_{\zeta\sim U}\{\zeta:\|\phi(\ov\omega,\zeta)-\psi(\ov\omega)\|\leq \rho\}>1/2,
$$
and define $\psi(\ov\omega)$ as once for ever fixed point of $Y$ when $\ov\omega\notin\ov\Omega$. For $x\in Y$, let also
$$
\wt\Omega[x]=\big\{\ov\omega\in\bR^m: \Prob_{\zeta\sim U}\{\zeta:\|\phi(\bar\omega,\zeta)-x\|\leq \rho\}>1/2\big\},
$$
note that $\wt{\Omega}[x]\subset\ov\Omega$. Let now $\wt{\Omega}^c[x]$ be the complement of $\wt{\Omega}[x]$; due to the origin $\rho$, we have for every $x\in Y$:
$$
\Prob_{\ov\omega\sim\cN(\cA(x),M^{-1}I_m)}\wt{\Omega}^c[x]\leq 2\upsilon.
$$
On the other hand, whenever $\ov\omega\in \wt{\Omega}[x]$, both sets
$$
\{\zeta:\|\phi(\ov\omega,\zeta)-x\|\leq\rho\}\;\;\mbox{and}\;\; \{\zeta:\|\phi(\ov\omega,\zeta)-\psi(\ov\omega)\|\leq\rho\}
$$
are subsets of $[0,1]$ of measure $>1/2$ and thus intersect, implying that $\|x-\psi(\ov\omega)\|\leq 2\rho.$ We conclude that for every $x\in Y$ the stemming from $x$ probability of the event $\|\psi(\ov\omega)-x\|> 2\rho$ is at most $2\upsilon$, that is,
$$
\overline{\RiskOpt}_{2\upsilon,M}[X_j]\leq 2\rho.
$$
Because $\rho$ may be arbitrary  $>\RiskOpt^\cR_{\upsilon,M}[Y]$, (\ref{eq221}) follows.\qed
When combining \rf{barrisk} and (\ref{eq22122}) we conclude that because $\epsilon\leq1/16$
one has
\[
\myr_j(\epsilon)\leq \myC_5\ln(L+L')\sqrt{\ln(m/\epsilon)}\RiskOpt_{{1\over16},M}[X_j]\quad\forall j
\leq N.
\]
Taken together with \rf{eq12222} the latter bound implies the last statement of the theorem. \qed

\end{document}